\documentclass{amsart}
\usepackage{amsmath}
\usepackage{float}
\usepackage{amssymb}
\usepackage{tikz}
\usepackage[all]{xy}
\usepackage{graphicx}
\usepackage{tikzsymbols}
\usepackage{textcomp}
\usepackage{parskip}
\usepackage{mathrsfs}
\usepackage[breaklinks=true]{hyperref}
\usepackage{epigraph}

\newtheorem{theorem}{Theorem}[section]

\newtheorem{question}{Question}
\newtheorem*{song*}{The song of despair}
\theoremstyle{definition}
\newtheorem{definition}[theorem]{Definition}

\newtheorem{claim}[theorem]{Claim}

\newenvironment{psmallmatrix}
  {\left(\begin{smallmatrix}}
  {\end{smallmatrix}\right)}

\theoremstyle{remark}
\newtheorem{remark}[theorem]{Remark}
\newcommand{\com}[1]{}

\newcommand*{\TakeFourierOrnament}[1]{{%
\fontencoding{U}\fontfamily{futs}\selectfont\char#1}}
\newcommand*{\danger}{\TakeFourierOrnament{66}}

\title{IntroSurvey of representation theory}

 \author{Nicolas Libedinsky }

\begin{document}


\begin{abstract}
There could be thousands of Introductions/Surveys of representation theory, given that it is an enormous field. This is just one of them, quite personal and informal. It has an increasing level of difficulty; the first part is intended for final year undergrads. We explain some basics of representation theory, notably Schur-Weyl duality and representations of the symmetric group. We then do the quantum version, introduce Kazhdan-Lusztig theory, quantum groups and their categorical versions. We then proceed to a survey of some recent advances in modular representation theory. We finish with twenty open problems and a song of despair.
\end{abstract}

\maketitle

\tableofcontents

\subsection*{Acknowledgements} You would not be reading this if it was not for Apoorva Khare, who was the engine for this work and also, with his warmth, kept me going in times of despair (my baby Gael was 1-5 months old when this paper was written). I would like to thank warmly Jorge Soto-Andrade for his many comments that improved an earlier version of this paper. Also to Juan Camilo Arias and Karina Batistelli for carefully correcting earlier versions of this manuscript. I would also like to thank  Stephen Griffeth, Giancarlo Lucchini, David Plaza, Gonzalo Jimenez and Felipe Gambardella for helpful comments. Finally, enormous thanks to Geordie Williamson who corrected the paper in super detail, finding some errors in a previous version. This project was funded by ANID project Fondecyt regular 1200061.

\part{Classical level}
\section{Quotes}

Groups are absolutely magnificent! They are magical beasts, beautiful monsters teeming with wonders. They are the breathing soul of mathematics. Although,  throughout history some people have dared to disagree, sometimes quite acidly. For example, the philosopher Simone Weil,  described by Albert Camus as ``the only great spirit of our times'', said: 

\vspace{0.2cm}

\textbf{``\textit{Money, machinery, algebra; the three monsters of current civilization\footnote{``Cahier I'', page I.100.}.''}}

\vspace{0.2cm}

Some people could be offended by this phrase, but as an algebraist, I am amazed and honored that algebra is included in her shortlist of monsters of the whole of civilization! In the past I have also sometimes referred to groups as monsters (see, for example, the first line of the present paper).  The physicist Sir James Jeans, was more humble in his criticism: 

\vspace{0.2cm}

\textbf{``\textit{We may as well cut out group theory. That is a subject that will never be of any use in physics}}\footnote{Discussing mathematics curriculum reform at Princeton University (1910), as quoted in Abraham P. Hillman, Gerald L. Alexanderson, ``Abstract Algebra: A First Undergraduate Course'' (1994).}

\vspace{0.2cm}
A funnier (and slightly more disrespectful) critique about how group theory actually works, comes from James Newman:

\vspace{0.2cm}

\textbf{\textit{``The theory of groups is a branch of mathematics in which one does something to something and then compares the result with the result obtained from doing the same thing to something else, or something else to the same thing.}}\footnote{ ``The World of Mathematics'' (1956) p.1534.}''

\vspace{0.2cm}

I will say nothing (due to outrage) about the first and third of these horrifying judgments. I will just do a little mosaic of quotes to answer Sir James Jeans. In order, they are from Irving Adler, George Whitelaw Mackey and Freeman Dyson: 

\vspace{0.2cm}

\textbf{\textit{ ``The importance of group theory was emphasized very recently when some physicists using group theory predicted the existence of a particle that had never been observed before, and described the properties it should have. Later experiments proved that this particle really exists and has those properties.}}\footnote{Quoted in ``Out of the Mouths of Mathematicians'' (1993) by R. Schmalz.}\textbf{\textit{'' ``Nowadays group theoretical methods, especially those involving characters and representations, pervade all branches of quantum mechanics.}}\footnote{George Whitelaw Mackey
``Group Theory and its Significance'', Proceedings American Philosophical Society (1973), 117, No. 5, 380.}\textbf{\textit{''
``We have seen particle physics emerge as the playground of group theory.}}\footnote{Joseph A. Gallian ``Contemporary Abstract Algebra'' (1994) p. 55}\textbf{\textit{''} }

\vspace{0.1cm}

\textbf{ BAM!!}



\section{Inside groups}

\epigraph{\textit{The introduction of the cipher $0$ or the group concept was general nonsense too, and mathematics was more or less stagnating for thousands of years because nobody was around to take such childish steps...}}{{Alexander Grothendieck}}

\vspace{0.2cm}

There are essentially two questions that you can ask  a group that you meet in the hood:

\vspace{0.2cm}

\begin{enumerate}
    \item What are you made of?     (looking inside)
    \item  How do you behave?     i.e.  where and how do you act?   (looking outside)
\end{enumerate}

\subsection{What are you made of?}\label{2.1}

Of course, like any group psychologist knows, both questions are oh so very interrelated. In this paper we care about Question (2). But let me say a few words about Question (1). That question  has totally obsessed mathematicians for over a century. For finite groups, they invented the following program to solve it: (a) Find all finite simple groups (the building blocks that form every finite group), and (b) Understand how are they glued together to produce  the specific group that you have met. 

\subsubsection{} The quest (a) of finding all simple groups was probably the biggest collaborative effort in the history of mathematics. Hundreds of mathematicians worked tirelessly from 1832 until 2004 and finished $\mathcal{THE\ LIST}$. It goes like this: $\mathbb{Z}/p\mathbb{Z}$ (the easiest groups), the alternating groups $A_n$ with $ n\geq 5$, the groups of Lie type (projective special linear groups $PSL_n(q)$, projective symplectic groups $PSp_{2n}(q),$ etc.) and finally $26$ sporadic groups (Monster, Baby Monster, etc.).

Some comments are in order. 
\begin{itemize}
    \item Already Galois in $1832$ knew that $\mathbb{Z}/p\mathbb{Z}$, $A_n$ with $ n\geq 5$ and $PSL_2(p)$ with $ p\geq 5$ are simple. Fantastic!
    \item The \textit{monstrous moonshine} is the unexpected connection between the monster group $M$ (the largest sporadic group) and modular functions (something in number theory). It was proved by Richard Borcherds who obtained the Fields Medal for it. I am saying this just because I want to give a mathematical version of Simone Weil's quote: 
    \vspace{0.3cm}
   \begin{center}
       \textbf{``\textit{Monster group, Baby monster, Monstrous moonshine; the three monsters of current mathematics}''}
   \end{center} 
\end{itemize}
\vspace{0.3cm}
\subsubsection{}\label{two}
Now let us speak about part (b) of the program: understand how simple finite groups can be glued together to form any group. Some people (see for example Section \ref{2.1}) use the metaphor that simple groups are ``building blocks'' that produce groups, in the same way as prime numbers produce all natural numbers. 

\begin{figure}[htp]\label{lego1}
    \centering
    \includegraphics[width=8cm]{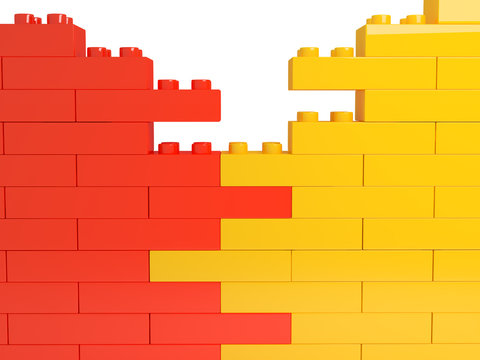}
    \caption{Mental image of a group using this metaphor}
    \label{bricks}
\end{figure}

\begin{figure}[htp]\label{lego2}
    \centering
    \includegraphics[width=8cm]{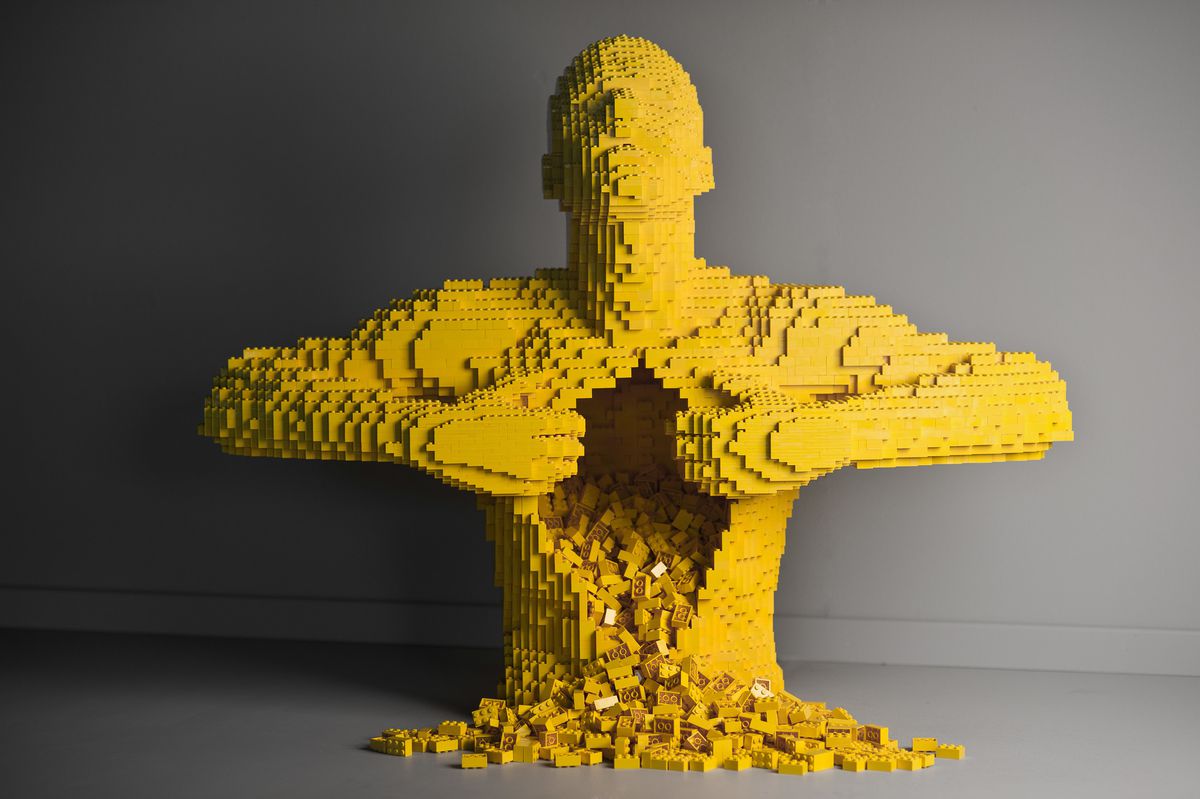}
    \caption{What groups really are}
    \label{lego}
\end{figure}

In hindsight, finding all finite simple groups was an incredibly hard problem for human beings, but part (b) of the project (understanding how  they can be glued together) seems several orders of magnitude harder. Probably not a problem for human beings.


\section{The two most beautiful groups in the world}

\epigraph{\textit{(On the concept of group:) ... what a wealth, what a grandeur of thought may spring from what slight beginnings.}}{{Henry Baker}}

The groups mentioned in the title of this section are, of course,  the symmetric group and the general linear group. The \emph{symmetric group} $S_n$ is the group of permutations of $\{1,2,\ldots, n\}$ and the \emph{general linear group} $GL_n(k)$ is the group of all $n\times n$ invertible matrices with coefficients in a field $k$. They are so important that even their subgroups have special names. The subgroups of $GL_n$ are called ``linear groups'' and the subgroups of $S_n$ are called ``finite groups". In this section we will propose three very different ways in which they are intimately related. This section could be summarized like this:

\ref{phantom} $GL_n(\mathbb{F}_q)$ is a  deformation of $S_n$. 

\ref{Bruhat} From $GL_n$ to $S_n$ and from $S_n$ to $GL_n$.

\ref{Duality} Do something to $GL_n$ and obtain $S_n$. Do it to $S_n$ and obtain back $GL_n$. 

\subsection{A mathematical phantom}\label{phantom} 
\epigraph{\textit{Nothing is more fruitful—all mathematicians know
it—than those obscure analogies, those disturbing reflections of one theory in another; those furtive caresses, those inexplicable discords; nothing also gives
more pleasure to the researcher.
The day comes when
the illusion dissolves; the yoked theories reveal their
common source before disappearing. As the Gita
teaches, one achieves knowledge and indifference at
the same time.}}{André Weil}
%
It was already noticed  in 1951 by Robert Steinberg \cite{St51} that $S_n$  and $GL_n(\mathbb{F}_q)$ have remarkable similarities (here $\mathbb{F}_q$ is the unique field with $q$ elements, where $q$ is a power of a prime number). He even used this analogy along with the representation theory of $S_n$ (see Section \ref{Sn}) to understand an important part of the representation theory of $GL_n(\mathbb{F}_q)$. 

In 1957, Jacques Tits \cite{Tits57} proposed the most incredibly bold version of the analogy found by Steinberg. He said $S_n=GL_n(\mathbb{F}_1)$, where $\mathbb{F}_1$ is the nonexistent field with one element. 

Let me stop right here and say a few words about this ``field''. First things first,  $\mathbb{F}_1$ does not officially exist, as the cipher $0$ was also nonexistent at some point, or the complex number $i$ was, for a long time as well, a mathematical phantom. 
It is quite probable that at some point the field with one element will have the kind of existence that we mathematicians love and respect (some kind of wide consensus in the definition and satisfying the properties one expects it to have), but for now it seems that it has some sort of existence, or, if you wish, it is just a highly suggestive idea. 

This astonishingly simple yet deep idea, the field of characteristic one, was not really taken seriously at the time Tits proposed it. More than 35 years elapsed until the field with one element appeared in the literature again, by Mikhail Kapranov and Alexander Smirnov (unpublished work). But in \cite{Manin}  Yuri Manin proposed (following ideas by  Christoph Deninger \cite{Deninger} and Nobushige Kurokawa \cite{Kurokawa})  that the Riemann hypothesis might be solved if one fully understands the geometry over $\mathbb{F}_1$, or, more precisely, if one is able to reproduce over $\mathbb{F}_1$ Andr\'{e} Weil's proof for the case of arithmetic curves over $\mathbb{F}_q$. Since then, the Riemann hypothesis is the main motivation for many people (like Cristophe Soul\'{e} \cite{Soule} or Alain Connes and Catarina Consani \cite{Connes}) to produce a geometry over $\mathbb{F}_1$. For an overview of the existing approaches to these geometries, see \cite{Lopez}.

Anyways, the point is that this is serious mathematics. Of course, one could solve the problem by  saying that in $\mathbb{F}_1$ we have $0=1.$ But if we do so, then a vector space  over $\mathbb{F}_1$ would be the zero vector space because every vector $v=1\cdot v=0\cdot v=0$. And that is far  from what we want. 

Let us get into more detail on the similarities observed by Steinberg.  Define the $q$-\emph{integer}:
$$[n]_q:=1+q+q^2+\cdots +q^{n-1}\in \mathbb{N}[q].$$
Of course, when $q=1$ one obtains $[n]_q=n.$ We will say that the $q$-deformation of the set $\mathbf{n}:=\{1,2,\ldots, n\}$ is the projective space $\mathbf{n}_q:=\mathbb{P}(\mathbb{F}_q^n)$ (i.e. the set of lines in the vector space $\mathbb{F}_q^n$) because the cardinality of $\mathbf{n}$ is $n$ and the cardinality of  $\mathbf{n}_q$ is $$[n]_q=\frac{q^n-1}{q-1}.$$ Likewise, the $q$-deformation of the set $$\mathbf{n}^k:=\{S\subseteq \mathbf{n}\, :\, \mathrm{card}(S)=k\}$$  is the set $$\mathbf{n}^k_q:=\{V\ \mathrm{subspace\ of\ } \mathbf{n}_q\, :\, \mathrm{dim}(V)=k\}$$ because (with some effort) one can prove that
$$ \mathrm{card}(\mathbf{n}^k)=\frac{n!}{(n-k)!k!} \ \  \ \  \mathrm{and }\ \ \ \  \mathrm{card}(\mathbf{n}^k_q)=\frac{[n]_q!}{[n-k]_q![k]_q!},$$
where $[n]_q!:=[n]_q\cdots [2]_q[1]_q.$ In other words, if  $\mathbb{F}_1$ were to exist, a vector space over $\mathbb{F}_1$ is just a set.

Finally, we can see the symmetric group as the set of all functions $f:\mathbf{n}\rightarrow \mathbf{n}$ sending $\mathbf{n}^k$ to itself, for all $k$. The $q$-deformation should be the set of maps  $f: \mathbf{n}_q\rightarrow \mathbf{n}_q$   sending $\mathbf{n}^k_q$ to itself, for all $k$. This is the group $PGL_n(\mathbb{F}_q)$, i.e. the quotient of $GL_n(\mathbb{F}_q)$ by the set of scalar matrices $\{\lambda\, \mathrm{I}_n : \lambda \in \mathbb{R}\}$ (these are the matrices that fix every point in $\mathbf{n}^k_q$)\footnote{This fits well with the fact that in classical mechanics there is a set describing the states of a system, while in quantum physics the states are vector lines in a Hilbert space $H$, thus in the projective Hilbert space $P(H)$.}.

You might complain that the deformation of $S_n$ that we have obtained is $PGL_n(\mathbb{F}_q)$ and not $GL_n(\mathbb{F}_q)$ as I promised, but one should not take this too seriously (we could even have obtained $SL_n(\mathbb{F}_q)$, as explained in the next remark), because these groups are almost the correct candidates for being deformations of $S_n$ but not quite: the cardinality of  $GL_n(\mathbb{F}_q)$ (or that of $PGL_n(\mathbb{F}_q)$)  does not tend to the cardinality of $S_n$ when $q \rightarrow 1$. But almost. This is the meat of the next remark and the mammoth footnote in it.

\begin{remark}
\textbf{For the more advanced reader}. In \cite{Tits57} Tits explains that if $G$ is a Chevalley group scheme ($PGL_n, GL_n, SL_n, Sp_{2n}, SO_n,$ etc.)  over $\mathbb{Z}$, then the cardinality of  the group of points of $G$ in $\mathbb{F}_q$ satisfies the formula: 
\begin{equation}\label{cheva}
    \mathbb{N}[q]\ni \frac{\mathrm{card}(G(\mathbb{F}_q))}{(q-1)^r} \xrightarrow{q\rightarrow 1} \mathrm{card}(W) 
\end{equation}
 where $r$ is the rank of $G$ and $W$ is its Weyl group. By ignoring\footnote{ This ``ignoring'' business bothers me. I would like some slight twist in the analogy for the numbers to really fit. I find intriguing that $(q-1)^r$ is the number of points in a maximal torus $T$, so one could replace $G(\mathbb{F}_q)$ by $G(\mathbb{F}_q)/T(\mathbb{F}_q)$ or even by $G(\mathbb{F}_q)/B(\mathbb{F}_q)$ (with $B$ a Borel) for the numbers to have the correct limit, but this only give homogeneous spaces and not groups, and part of the magic of this analogy is that it also extends to the level of representations (see Section \ref{another}). One may counterargue this by saying that each homogeneous space has a corresponding groupoid (à la Connes), and maybe the groupoid representation theory of $G(\mathbb{F}_q)/B(\mathbb{F}_q)$ is similar to the group representation theory of the whole group. Another option, staying in the realm of groups, is that maybe there is something like an $\mathbb{F}_q$-version (instead of the usual $\mathbb{Z}$ version) of the affine Weyl group of cardinality $(q-1)^r\times \mathrm{card}(W), $ that would make the numbers agree. } the $(q-1)^r$ in the formula, he defines $G(\mathbb{F}_1):=W$.
\end{remark}
\begin{remark}
One striking feature of this analogy is that if $PGL_n(\mathbb{F}_q)$ is the quantum deformation of $S_n$, the determinant is the quantum deformation of the sign function $$\mathrm{sgn}:S_n\rightarrow \{+1,-1\}.$$ Indeed the sign of the determinant is positive (resp. negative) if a basis of your vector space is sent to another basis with the same orientation (resp. opposite). When your map is between $\mathbb{F}_1$-vector spaces, a.k.a. finite sets, the orientation is given by the sign. This means that $A_n$ (the kernel of the sign homomorphism) corresponds to $PSL_n(\mathbb{F}_q)$ (the kernel of the determinant homomorphism). The striking part is that $A_n$ and $PSL_n(\mathbb{F}_q)$ are both simple groups (except for very small $n$ and $q$).
\end{remark}

\begin{remark}
Yuri Manin observed that the equation $SL_n(\mathbb{F}_1)=S_n$ implies (using the Barratt-Priddy-Quillen theorem) that $K_i(\mathbb{F}_1)$, the higher $K$-theory of $\mathbb{F}_1$ must be the stable homotopy groups of spheres $\lim_{m\to \infty}\pi_{m+i}(S^m)$. What I find so fascinating about this, is that Geordie Williamson has realized that the $p$-canonical basis (allegedly the most important mystery in modern modular representation theory, see Section \ref{pcan}) could potentially be understood using some categories obtained from stable homotopy theory. 
\end{remark}

\subsection{François Bruhat and the Han Dynasty}\label{Bruhat}

There is a  beautiful observation: the symmetric group $S_n$ is a subgroup of $GL_n(\mathbb{F}_p)$. Each permutation $w\in S_n$ corresponds to $ \dot{w}$ a \textit{permutation matrix} (that is, a matrix that has exactly one entry of $1$ in each row and each column and $0$'s elsewhere) obtained in the following way. Start with the $n\times n$ identity matrix $I_n$ and permute \footnote{This is a particular example (for $GL_n$ and $S_n$) of a general construction that associates to a connected reductive algebraic group over $k$ its Weyl group $W:=N_G(T)/T$, where $T$ is a maximal torus.} the columns according  to $w$. 

It is implicit in the paper \cite{St51} by Robert Steinberg  that, if $B$ is the set of upper triangular matrices, then

\begin{equation}\label{bd}
GL_n(\mathbb{F}_p)=\coprod_{w\in S_n}B\dot{w}B,
\end{equation}
i.e. $GL_n(\mathbb{F}_p)$ is partitioned into a disjoint union of double cosets parametrized by $S_n$. The difficult part is to realize the phenomenon. The proof is  relatively easy, it follows from a very slight modification of 
the Gaussian elimination algorithm. 
Equation \eqref{bd} is beautiful although slightly annoying in that an element $b_1\dot{w}b_2$ can usually  be also written as $b'_1\dot{w}b'_2,$ for $b_1 \ne b_1'$ and $b_2 \ne b_2'$, but there are ways to refine this partition to solve that. 
This is generally called the \textit{Bruhat decomposition}, because in \cite{Bruhat} François Bruhat formulated for the first time that Equation \eqref{bd} can be generalized for general semisimple groups over $\mathbb{C}$ and stated that he had verified it for all classical groups (this was later proved  over $\mathbb{C}$ by  Harish-Chandra \cite{HC} and for all algebraically closed fields  by Claude Chevalley \cite{Chevalley}).

Bruhat decomposition allows us to reduce questions about $GL_n$ to questions about $S_n$, and thus it is a vital tool both to understand the structure and the representations of $GL_n$ as we will see in later chapters. We have already used Bruhat decomposition without saying it: the order of a Chevalley group over a finite field appearing in Equation \eqref{cheva} was computed in \cite{Chevalley} using
Bruhat decomposition.

Let me finish this section with a fun historical fact that makes apparent how fundamental the Bruhat decomposition is. Equation \ref{bd} stays true if one changes the field $\mathbb{F}_p$ by the complex numbers $\mathbb{C}$. In this new version of Bruhat decomposition, the part of  Equation \eqref{bd}  corresponding to $w_0$, \textit{the longest element} in $S_n$ (as a function from $\mathbf{n}$ to $\mathbf{n}$, $w_0$ is the map $i\mapsto n-i+1$) essentially
appears in Chapter 8 of the Chinese classic ``The Nine Chapters on the Mathematical
Art'' from the second century, Han dynasty \cite{Lusztig3} (that book  is also the first place where negative numbers appear in the literature). The point is that $\dot{w_0}B\dot{w_0}$ is the set of lower triangular matrices and they knew that ``almost any matrix'' could be written as a product of a lower and an upper triangular matrix.  


\subsection{Soul mates}\label{Duality}

If the reader is not yet convinced that $S_n$ and $GL_n$ are two faces of the same coin, let us see a procedure to obtain one from the other. 

Let me first introduce the \textit{group algebra} $\mathbb{C}[S_n]$. As a vector space it is 
$$\mathbb{C}[S_n]:=\bigoplus_{w\in S_n}\mathbb{C}w,$$
i.e. the elements of the symmetric group form a $\mathbb{C}$-basis of $\mathbb{C}[S_n]$. An algebra is a vector space where you can multiply  the elements, and this is done here in the obvious way. To be more precise, the product of two basis elements of $\mathbb{C}[S_n]$ is just their product in $S_n$, and then we extend by linearity. So, for example, if $x,y,z\in S_n,$ in $\mathbb{C}[S_n]$ we would have 
$$(2x-5y)(x+3z)=2x^2-10yx+6xz-15yz. $$

Let $V=\mathbb{C}^n$. In the vector space $V^{\otimes n}$ there are two natural actions. The group $S_n$ acts by permuting the factors $$(v_1\otimes v_2\otimes\cdots \otimes v_n)\cdot w= v_{w(1)}\otimes v_{w(2)}\otimes\cdots \otimes v_{w(n)}, $$
and on $GL_n(\mathbb{C})\cong GL(V)$ acts by
$$g\cdot (v_1\otimes v_2\otimes\cdots \otimes v_n)= g(v_1)\otimes g(v_2)\otimes\cdots \otimes g(v_n), $$

The latter action gives a map $\rho: GL(V)\rightarrow \mathrm{End}(V^{\otimes n})$ that can be enhanced into a highly non injective algebra homomorphism  $\rho': \mathbb{C}[GL(V)]\rightarrow \mathrm{End}(V^{\otimes n})$. By abuse of notation we will call in the next paragraph $\mathbb{C}[GL(V)]$ the image of $\mathbb{C}[GL(V)]$ under $\rho',$ which is also the algebra generated by the image of $\rho$.

One version of a theorem known as \textit{Schur-Weyl duality} is the following two isomorphisms of algebras: 

$$ \ \, \mathrm{End}_{\mathbb{C}[S_n]}(V^{\otimes n})\ \, \cong \,  \mathbb{C}[GL(V)],$$
$$ \mathrm{End}_{\mathbb{C}[GL(V)]}(V^{\otimes n})\cong\, \mathbb{C}[S_n].$$

\vspace{0.2cm}
Here $\mathrm{End}$ means all linear endomorphisms and the subscript means ``invariant under''. So, for example, $$\mathrm{End}_{\mathbb{C}[GL(V)]}(V^{\otimes n})=\mathrm{End}_{GL(V)}(V^{\otimes n})$$ means all linear functions  $f: V^{\otimes n}\rightarrow V^{\otimes n}$ that satisfy $ f(g\cdot v)= g\cdot f(v)$ for all $g\in GL(V)$  and for all $v\in V^{\otimes n}.$ In other words, this means that $f$  commutes with all  $g$. 

This is one of the most beautiful theorems I have ever seen in my short life, and it will be a guiding principle throughout this paper. As in Sections \ref{phantom} and \ref{Bruhat}, there are also versions of this theorem for other classical groups, such as the orthogonal or the symplectic group. This theorem is an example of a general pattern in mathematics (as Poincar\'{e} duality  or Koszul duality in Section \ref{pcan}), that of a duality: apply some procedure to some mathematical object A and obtain B, then apply the same procedure to B, only to obtain A back. 




\section{Looking outside of groups}\label{repsn}
\epigraph{\textit{An ounce of action is worth a ton of theory.}}{{Ralph Waldo Emerson}}
If we are interested in the ``outer structure'' of a group $G$, the first question we should ask is what are the group morphisms from $G$ to other groups? But what ``other groups'' are interesting to study? 

Let us start by asking what are the group morphisms between our two amazingly related groups. The answer is that the set of homomorphisms from $GL_m$ to $S_n$ is probably one of the dullest sets in the world, while the set of homomorphisms from $S_n$ to $GL_m$ is probably one of the most fascinating sets in the world (related to fractals, dynamical systems, homotopy groups of the spheres, Langlands program, etc.). 

\subsection{Nice problem, dull answer} Consider a morphism of groups $\varphi: GL_m(\mathbb{C})\rightarrow S_n.$ As any group homomorphism, $\varphi$ satisfies that $\varphi(h^r)=\varphi(h)^r$ for all $h\in GL_m(\mathbb{C})$ and $r\in \mathbb{Z}$. We also know that any element $w\in S_n$ satisfies $w^{n!}=e$ (by Lagrange's theorem, as the order of $S_n$ is $n!$). So, we have that $\varphi(h^{n!})=e$ for all $h\in GL_m(\mathbb{C})$. But any $g\in GL_m(\mathbb{C})$ admits, for any natural number $r,$ an $r^{\mathrm{th}}$-root $g^{\frac{1}{r}}\in GL_m(\mathbb{C})$ (this is not hard to prove\footnote{For a heuristic argument, it is quite obvious that diagonal matrices admit an $r^{\mathrm{th}}$-root,  so the same is valid for every diagonalizable matrix, a set of matrices that is dense in $GL_m(\mathbb{C})$. }) so $$\varphi(g)=\varphi((g^{\frac{1}{n!}})^{n!})=e.$$ We conclude that the only element in $\mathrm{Hom}(GL_m(\mathbb{C}), S_n)$ is the trivial homomorphism. Boring\footnote{If instead of looking these objects as groups, one looks them as topological spaces, $GL_n(\mathbb{C})$ is connected and $S_n$ is discrete, so any continuous homomorphism from $GL_n(\mathbb{C})$ to $S_n$ is also constant. Boring${}^2$.}. 

\subsection{Nice problem, fascinating answer}\label{Sn} 
Now we come to the question of understanding the set $\mathrm{Hom}(S_n, GL_m(\mathbb{C}))$. A map in this set (for varying $m$) is called a \textit{representation over} $\mathbb{C}$ of the group $S_n$ and it is also called a \textit{a linear action of $S_n$}. One could imagine that this problem is much more fun and interesting that the other way around, given that most groups emerge not as abstract groups but as linear groups (i.e. as groups of matrices). We will now consider the easiest field, the complex numbers $\mathbb{C},$ and in that case we will be able to solve the problem. This case is extremely important for the history of mathematics and related to many things, from statistical mechanics to probability theory. But in my opinion
what is  more mysterious, far from being understood, and related to distant and interesting mathematics, is when you replace $\mathbb{C}$ by a field of positive characteristic. But let us not rush into it. 

\subsubsection{Some ``abstract nonsense''}
When working over $\mathbb{C}$ (the complex numbers make everything easy) we have that there are really building blocks (irreducible representations, as explained in two more paragraphs), and gluing them together (direct sum of  representations, as explained in three more paragraphs) is really easy. It is essentially like Figure \ref{bricks} (not at all like Figure \ref{lego}).  We will think of $GL_m(\mathbb{C})$ as $GL(V)$ with $V$ an $m$-dimensional vector space over the complex numbers.

Consider a homomorphism $\rho:S_n\rightarrow GL(V)$. Then for each $w\in S_n$ we have that $\rho(w)\in GL(V)$, or in other words, each element $w$ of the symmetric group ``is'' a linear map $\rho(w): V\rightarrow V.$

The homomorphism $\rho:S_n\rightarrow GL(V)$
 is called an \textit{irreducible representation} if there is no linear subspace $U\subseteq V$ different from $\{0\}$ and $V$ stable under $\rho$, i.e. such that $\rho(w)(u)\in U$ for all $w\in S_n$ and $u\in U.$

 If $\rho:S_n\rightarrow GL(V)$ and $\rho':S_n\rightarrow GL(V')$ are two representations of $S_n$, the \textit{direct sum} $\rho+\rho':S_n\rightarrow GL(V\oplus V')$ is defined by the only reasonable formula  $$(\rho+\rho')(w)(v,v'):=\rho(w)(v)+\rho(w)(v')\quad \mathrm{for\ any}\ w\in S_n, v\in V, v'\in V'. $$

A fundamental result, known as \textit{Maschke's theorem}\footnote{In categorical terms, this theorem says that the category of representations of $S_n$ over $\mathbb{C}$ is semi-simple.} (due to Heinrich Maschke \cite{Maschke}) says  that every representation of $S_n$ over $\mathbb{C}$  is a direct sum of irreducible representations. So we just need to understand those. 

\subsubsection{A very concrete geometric action }\label{concrete}
\epigraph{\textit{Groups, as people, will be known by their actions.}}{{Guillermo Moreno}}

Every time that $S_n$ acts on a set $X=\{a,b,\ldots, c\}$\footnote{I learnt this notation from Wolfgang Soergel, and I love it, although it is impossible not to ask oneself, what are the names of the elements between $b$ and $c$?},  it also acts on the vector space
\begin{equation}\label{cx} \mathbb{C}X:=\{\lambda_aa+\lambda_bb+\cdots+\lambda_cc\}, 
\end{equation}
via the action 
$$ w(\lambda_aa+\lambda_bb+\cdots+\lambda_cc)=\lambda_aw(a)+\lambda_bw(b)+\cdots+\lambda_cw(c).$$
In other words, we obtain a representation $S_n\rightarrow GL(\mathbb{C}X)$, that is called the \textit{permutation representation} associated to the action. So we must ask ourselves, on what sets does $S_n$ naturally act?

Of course, by definition, it acts on the set $\mathbf{n}:=\{1,2,\ldots, n\}.$ But it also acts on the set of unordered pairs $\{i,j\}$ (this set has cardinality 
${n\choose 2}$). 
We can generalize this easily, and see that $S_n$ acts on the set of  subsets of $\mathbf{n}$ of a fixed cardinality $k$, i.e. in $\mathbf{n}^k$ in the language of Section \ref{phantom}.
We can further generalize this action and make $S_n$ act on \textit{flag varieties}\footnote{This is non-standard notation. I call them like that in analogy with flag varieties in $GL_n(\mathbb{F}_q)$ in the sense of Section \ref{phantom}.}. Instead of giving a general definition of flag varieties, let me show an example. The group $S_8$  (recall the notation $\mathbf{8}:=\{1,2,3,4,5,6,7,8\}$ defined five lines ago) acts on the flag variety
$$ \mathcal{FL}_{1,3,6,8}:=\{A\subset B\subset C \subset \mathbf{8}\ : \ \mathrm{card}(A)=1, \mathrm{card}(B)=3, \mathrm{card}(C)=6  \}, $$
(we add the subscript $8$ in $\mathcal{FL}$ because we are in the $S_8$ case). In the following image we see how the permutation $(12345678)$ (the one that sends $1\mapsto 2, 2\mapsto 3, \ldots 8\mapsto 1$) acts on the flag $\{3\}\subset \{2,3,8\}\subset \{1,2,3,5,7,8\}.$
\begin{figure}[htp]\label{fig}
   \centering
  \includegraphics[width=9cm]{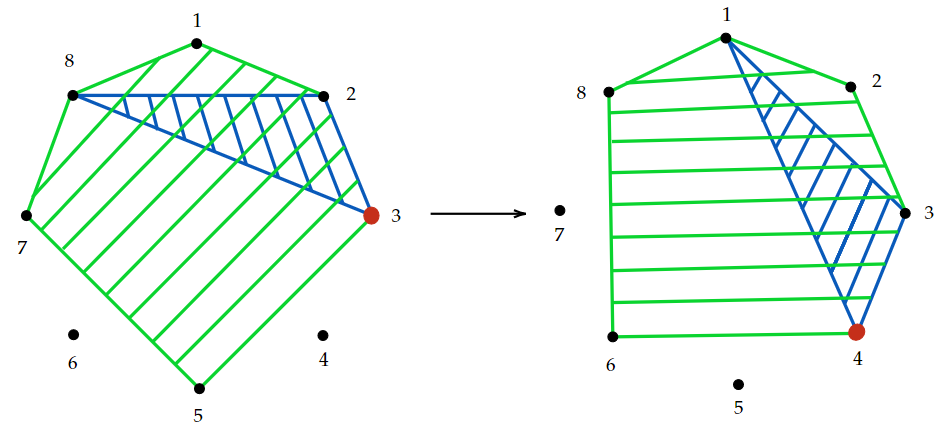}
  \caption{Action of  $(12345678)$ on the flag $\mathcal{FL}_{1,3,6,8}$}
\end{figure}

In each permutation representation $\mathbb{C}\mathcal{FL}_{a,b,\ldots, c}$  (see Equation \eqref{cx}) with $c=n$, there is a special subspace giving an irreducible representation of $S_n$ called a \textit{Specht module}, after Wilhelm Specht, that we will soon define (we will define a Specht module, not a Wilhelm Specht). We will see that Specht modules are the only irreducible representations of $S_n$. 

\subsubsection{Avoiding repetitions}\label{repetitions}
One small business that we have to deal with, is that some permutation representations associated to different flag varieties are obviously isomorphic and we don't want to count Specht modules more than once. For example, for $n>3,$ the map sending the flag $\{a\}\subset \{a,b,c\}\subset \mathbf{n}$ to the flag $\{b,c\}\subset \{a,b,c\}\subset \mathbf{n}$ gives a bijection $\mathcal{FL}_{1,3,n} \cong \mathcal{FL}_{2,3,n}$ that is also a morphism of representations (i.e. if $F\in \mathcal{FL}_{1,3,n}$ is sent to $F'\in \mathcal{FL}_{2,3,n}$ then $w(F)$ is sent to $w(F')$ for all $w\in S_n$). This is called an \textit{isomorphism of $S_n$-actions,} and clearly produce isomorphic permutation representations.  

Let us use the following new convention to explain in general which flag varieties are isomorphic. $$\mathcal{FL}(a,b-a,c-b,\ldots, e-d):=\mathcal{FL}_{a,b,c\ldots d,e}.$$
In this new notation, the fact that we are working with $S_n$ (i.e. $e=n$) translates to $\lambda_1+\lambda_2+\cdots +\lambda_r=n$ in $\mathcal{FL}(\lambda_1,\lambda_2,\ldots ,\lambda_r).$ To this new notation we add another compatible notation: while denoting flags, we will avoid repetition of elements, so $\{a\}\subset\{a,b,c\}\subset\{a,b,c,d,e,f\}$ will be denoted by $\{a\}\subset\{b,c\}\subset\{d,e,f\}.$

So, for example, the flags in Figure 3 belong to $\mathcal{FL}(1,2,3,2).$ As before, the map that sends the flag $$\{a\}\subset\{b,c\}\subset\{d,e,f\}\subset\{g,h\}$$
to the flag $$\{d,e,f\}\subset\{b,c\}\subset\{g,h\}\subset\{a\}$$ gives  an isomorphism of $S_n$-sets between the flag varieties $$\mathcal{FL}(1,2,3,2)\cong \mathcal{FL}(3,2,2,1).$$
It is easy to see that any rearrangement of the sequence $(1,2,2,3)$ will give an isomorphic flag variety, so in order to avoid repetitions of flag varieties, we have to choose one such rearrangement. We will choose the only sequence that is non increasing (in this case $(3,2,2,1)$).

So, as a summary, for any tuple $\lambda=(\lambda_1, \lambda_2, \ldots,\lambda_r)$ such that $\sum_i\lambda_i=n$ and $\lambda_1\geq \lambda_2\geq \cdots \geq \lambda_r$ (these are called \textit{partitions of} $n$) we will define a unique Specht module inside the representation\footnote{In the literature, an element of $\mathcal{FL}(\lambda)$  is called a $\lambda$-\textit{tabloid}. } $$\mathbb{C}\mathcal{FL}(\lambda).$$

\subsubsection{Specht modules, finally}
Imagine that you want a subgroup of $S_n$ to act on a flag variety $\mathcal{FL}(\lambda)$.  Take, for example, the flag 
\begin{equation}\label{F}
    F=\{1,4,5,2\}\subset \{6,7,3\}\subset \{8\}\in \mathcal{FL}(4,3,1).
\end{equation} We can permute it by the group $C_{(F)}=S_{\{1,6,8\}}\times S_{\{4,7\}}\times S_{\{5,3\}}\times S_{\{2\}}$. The notation $S_{\{i,j,\ldots, k\}}$ means the subgroup of $S_n$ permuting only the integers in brackets. This choice of group is natural since $\{1,6,8\}$ are the first elements of each subset of $F$, $\{4,7\}$ the second elements, and so on. Implicit in the definition of this subgroup is the fact that we took the sequence $(1,4,5,2)$ instead of the set $\{1,4,5,2\}$, the sequence $(6,7,3)$ instead of the set $\{6,7,3\}, $ etc. Let us call a flag where each set is replaced by an ordered set an  \textit{ordered flag}, and denote it by $\mathcal{OFL}(\lambda)$\footnote{In the literature, an element of $\mathcal{OFL}(\lambda)$  is called a $\lambda$-\textit{tableau}. }. For example
$$(F)=(1,4,5,2)\subset (6,7,3)\subset (8)\in \mathcal{OFL}(4,3,1).$$

It is clear that for any $\lambda$ partition of $n$ the cardinality of $\mathcal{OFL}(\lambda)$   is $n!$. We defined the ordered flags so that the group $C_{(F)}$ is well defined for $(F) \in \mathcal{OFL}(\lambda)$. Each ordered flag $(F) \in \mathcal{OFL}(\lambda)$ defines uniquely a flag $\{F\} \in \mathcal{FL}(\lambda)$ by forgetting the order. 

Define, for each $(F) \in \mathcal{OFL}(\lambda)$

$$e_{(F)}=\sum_{w\in C_{(F)}}\mathrm{sgn}(w){\{F\}} \in \mathbb{C}\mathcal{FL}(\lambda)  $$
One can prove that $we_{(F)}=e_{w(F)}$ so the subspace of $\mathbb{C}\mathcal{FL}(\lambda)$ spanned by the $e_{(F)}$ is invariant
under $S_n$. It is called the \emph{Specht
module} $S^{\lambda}$. Although the Specht module is generated as a vector space by $n!$ vectors, they are highly linearly dependent (the dimension is much less than $n!$). One can prove that the set $$\{S^{\lambda}\ \mathrm{ for }\ \lambda\ \mathrm{ partition\  of }\ n\}$$ is a complete set (up to isomorphism) of irreducible representations of $S_n$, without repetitions. With this we finish the question of understanding the set $\mathrm{Hom}(S_n, GL_m(\mathbb{C}))$.

\subsubsection{A tiny part of everything}\epigraph{\textit{Everything is representation theory}}{Israel Gelfand}

I will not go into the immense number of applications and relations of this story with other parts of mathematics. I will just mention one very  beautiful relation with probability theory \cite[p. 139]{Diaconis}\footnote{In the book \cite{Diaconis} that I cited there is a mistake in the formula,  the $\prod_i (\lambda_i)!$ is missing (thanks Geordie Williamson for noticing the mistake and Valentin Feray for explaining me the correct version of the theorem).}. Let  $\lambda=(\lambda_1, \lambda_2, \ldots, \lambda_r)$ be a partition of $n.$
If one puts $\lambda_1$-ones, $\lambda_2$-twos, $\ldots$, $\lambda_r$-$r'$s into an urn and draws all the  numbers without replacement,  then the
chance that at each stage of the drawing, the number of ones $\geq $ the number of twos $\geq \cdots \geq$ the number of $r'$s is equal to $$\mathrm{dim}(S^{\lambda})\frac{\prod_i (\lambda_i)!}{n!}.$$

\subsubsection{Educated guess on how Specht came out with this crazy stuff}\label{educated}
There is a big part of representation theory of finite groups that I have ignored until now. It is related to the following concept. If $\rho:S_n\rightarrow GL_m(\mathbb{C})$ is a representation, its \textit{character} is the function $\chi_{\rho}: S_n\rightarrow \mathbb{C}$ defined by the formula $\chi_{\rho}(w)=\mathrm{tr}(\rho(w))$. 

There are lots of theorems about characters, but the most important is that a complex representation  is determined (up to isomorphism) by its character. On the other hand, one can compute characters in many cases without knowing the representations.  So what I think that happened is that someone realized  just using characters Claim \ref{claim1} (that is Formula 2.1.2 in \cite{GK}), claim that I will now explain.

For each $\lambda$ partition of $n$ there is a transpose partition $\lambda'$ that is essentially\footnote{Technically $\lambda'$ is characterized uniquely by the condition $\lambda_i\geq j$ if and only if $\lambda_j'\geq i$.} the partition associated to $C_{(F)}$ above\footnote{Usually in the presentation of this theory one associate to each  $\lambda $ partition of $n$ a Young diagram and the transpose partition is literally the transpose of the Young diagram, thus explaining the name ``traspose''. I will not explain Young diagrams. }. For example, in Equation \eqref{F}, $\lambda=(4,3,1)$ and $\lambda'=(3,2,2,1)$.

On the other hand we have defined, for an action of $S_n$ on $X$ an action on $\mathbb{C}X$, but there is another natural action that we will call $\mathbb{C}X^{\mathrm{sgn}}$, the ``sign action'' defined by 
$$ w(\lambda_aa+\lambda_bb+\cdots+\lambda_cc)=\mathrm{sgn}(w)(\lambda_aw(a)+\lambda_bw(b)+\cdots+\lambda_cw(c)).$$

 \begin{claim}\label{claim1}
 If one writes $\mathbb{C}\mathcal{FL}(\lambda)$ and $\mathbb{C}\mathcal{FL}(\lambda')^{\mathrm{sgn}}$ as a sum of irreducible representations, they have exactly one in common: $S^{\lambda}$.
 \end{claim}
Remark that in this claim one takes the permutation representation associated to $\lambda$ and then the permutation representation associated to the transpose of $\lambda$, twisted by the sign. 
  In this situation  the last fundamental theorem about representations of finite groups that I haven't explained yet, called \textit{Schur's lemma} (due to Issai Schur)  comes in aid. When two representations have exactly one irreducible representation in common, there is exactly one morphism of representations (modulo multiplying by a complex number) between them. So, you need to find this unique map $$f: \mathbb{C}\mathcal{FL}(\lambda')^{\mathrm{sgn}}\rightarrow \mathbb{C}\mathcal{FL}(\lambda),$$
and the image of this map will be your $S^{\lambda}$.

\subsection{Another nice problem, similar/different solution}\label{another}

Now we can wonder about the group homomorphisms \begin{equation}\label{Hom}
    \mathrm{Hom}(GL_n(\mathbb{C}), GL_m(\mathbb{C}))
\end{equation} Let us be more humble and ask that question for $n=1$ and $ m=1,$ i.e. group endomomorphisms of the multiplicative group of complex numbers $\mathbb{C}^{\times}$. Let us be even more humble and ask for some very special group endomomorphisms of $\mathbb{C}^{\times}$, the field automorphisms of $\mathbb{C}$ restricted to $\mathbb{C}^{\times}$. Even those  \cite{Yale} are uncountable, quite crazy  and their construction require the axiom of choice. This proves that this is not the right question. We need to add some structure to  $GL_n(\mathbb{C})$ in order for the question to have a nice answer. This is the oldest mathematical technique with questions:  if you can't beat them, change the questions. In section \ref{fancy} I will explain what an algebraic representation is. The right question (because it has a nice answer) is ``what are the algebraic representations of $GL_n(\mathbb{C})$?'' (so we are asking for a very nice and small subset of the set \ref{Hom}).

The answer is relatively simple, but we will not explain the details. For any $\lambda=(\lambda_1, \lambda_2, \ldots, \lambda_n)$  with $\lambda_1 \geq \lambda_2 \geq \cdots \geq \lambda_n$ (notice that the $n$ in $\lambda_n$ is the same as in $GL_n(\mathbb{C})$) there is   an irreducible representation $V_{\lambda}$ of $GL_n(\mathbb{C})$ called a \textit{Weyl module}. It is  also an irreducible representation of $SL_n(\mathbb{C})$. 

\begin{remark}\label{domi}
For the reader mostly used to associating an irreducible representation to any dominant weight  $\mu_1\omega_1+\cdots +\mu_n\omega_n$ (see Section \ref{gage} for the definition), to this weight you associate the partition $\lambda=(\lambda_1, \lambda_2, \ldots, \lambda_n)$ where  $\lambda_k=\sum_{i=1}^k\mu_i.$ We already saw these two notations when speaking of flag varieties in Section \ref{repetitions}.
\end{remark}

Every  irreducible algebraic representation of $GL_n(\mathbb{C})$ is isomorphic to $V_{\lambda}\otimes \mathrm{Det}^i$ for some $i\in \mathbb{Z}$, where $\mathrm{Det}:GL_n(\mathbb{C}) \rightarrow GL_1(\mathbb{C})\cong \mathbb{C}^{*}$ is the determinant, and $V_{\lambda}\otimes \mathrm{Det}^i \ncong V_{\lambda'}\otimes \mathrm{Det}^{j}$ if $\lambda \neq \lambda'$ or $i\neq j.$ For $SL_n(\mathbb{C})$, the set $\{V_{\lambda}\}$ with $\lambda$ as above is a complete set of irreducible algebraic representations with no repetitions. 

We conclude that the maps from our two favorite groups to $GL_m$ give us a rich theory. But now, as the cherry on the cake, let us see another version of Schur-Weyl duality. Let us use the same terminology as in \ref{Duality}. It says that under the action of the group $S_n\times GL(V),$ the tensor space decomposes into a direct sum of tensor products of irreducible modules for the two groups.

$$ V^{\otimes n}= \bigoplus_{\lambda \vdash n} S^{\lambda}\otimes V_{\lambda}, $$

 \noindent where $\lambda \vdash n$ means that  $\lambda$ is a partition of $n$. In the rest of this paper we will attempt to understand the quantum and the categorical versions of Schur-Weyl duality, thus introducing some of the most exciting objects in modern representation theory, as they appear in the following diagram. 

\begin{figure}[H]
    \centering
    \includegraphics[width=5.5cm]{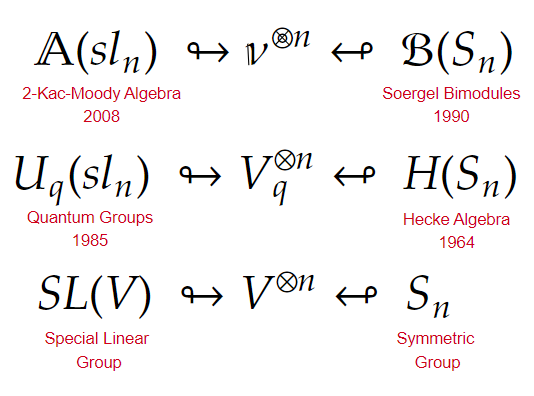}
    \caption{Schur-Weyl dualities}
    \label{swd}
\end{figure}

\section{The affine symmetric group or Aristotle's mistake that lasted two millennia}
\epigraph{\textit{Algebra is the offer made by the devil to the mathematician... All you need to do, is give me your soul: give up geometry.}}{{Michael Atiyah}}

\subsection{You know nothing, Aristotle}\footnote{Game of Thrones reference.} Plato thought that all matter is the disjoint union of a few basic polyhedral units, i.e. that one can fill the space with them without any gaps.  He assumed that there are some  polyhedra (apart from parallelepipeds) that could do that on their own. Later, Aristotle claimed that, from the five regular polyhedra, not only with the cube, but also with regular tetrahedra one can also fill the space  without gaps. That is the mistake in the title of this section, and although I don't like to repeat myself: it  lasted two millennia. Two millennia \cite{Senechal}. 
 
 It is quite funny that many of the later Aristotelian scholars realized that something was wrong, but they assumed that somehow they must be mistaken because Aristotle was oh so very wise and they were oh so very unwise \cite{Senechal}.  While trying to justify Aristotle's
 erroneous assertion, they raised the question of \textit{which tetrahedra actually do fill space}. To this day that problem is open\footnote{Fast question: do open problems give the set of all mathematical problems the structure of a topological space?}. What is clearly solved is that Aristotle was wrong\footnote{See \cite[p. 230]{Senechal} where an amusing (and false) theory of angles is developed by Averroes to explain why Aristotle was right.}  in this (and in some other things).
 
 A sketch of the proof of why Aristotle is wrong goes like this. One can prove that any space-filling tessellation of polyhedra involves packing the polyhedra around their edges, their vertices, or both. So we can't do this with regular tetrahedra because both things are impossible: 
 \begin{itemize}
     \item The dihedral angles (the angles between adjacent faces) of a regular tetrahedron is about $70.32^{\circ}$, which does not divide $360^{\circ}$. Therefore no space-filling tessellation involves regular tetrahedra packed around their edges (see the figure).
     \item The solid angle subtended by a tetrahedron's vertex is about $0.55$ steradians (a \textit{steradian} is a three dimensional analogue of the radians, see wikipedia for details), which does not divide $4\pi$. Therefore no space-filling tessellation involves regular tetrahedra packed around a vertex. 
 \end{itemize}

\begin{figure}[H]\label{Aristoteles}
    \centering
    \includegraphics[width=6cm]{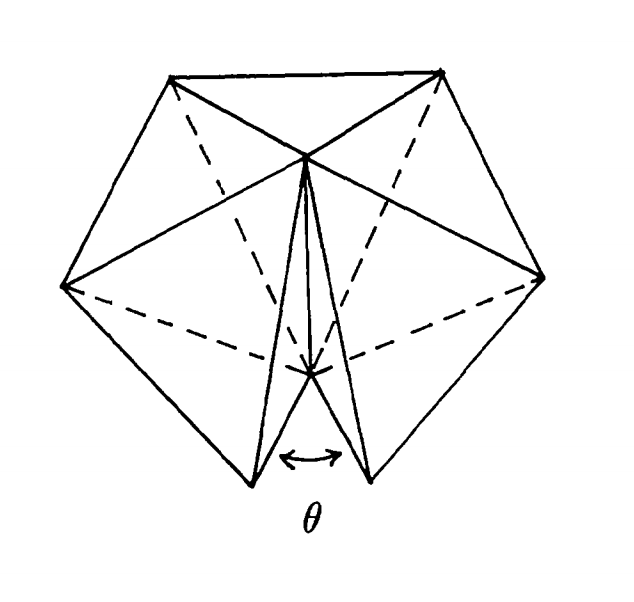}
    \caption{Five regular tetrahedra fitted around an edge}
\end{figure}

Currently  (see \cite{Chen}), the highest percentage of the space achieved with a packing of regular  tetrahedra is $85.63\%$.

\subsection{Why this mistake?}\label{why}
Probably Aristotle's mistake emerges from the fact that there is a tiling of $\mathbb{R}^2$ by ``2-dimensional regular tetrahedra'', a.k.a. equilateral triangles (see the image on the left). 

\begin{figure}[H]
    \centering
    \includegraphics[width=14cm]{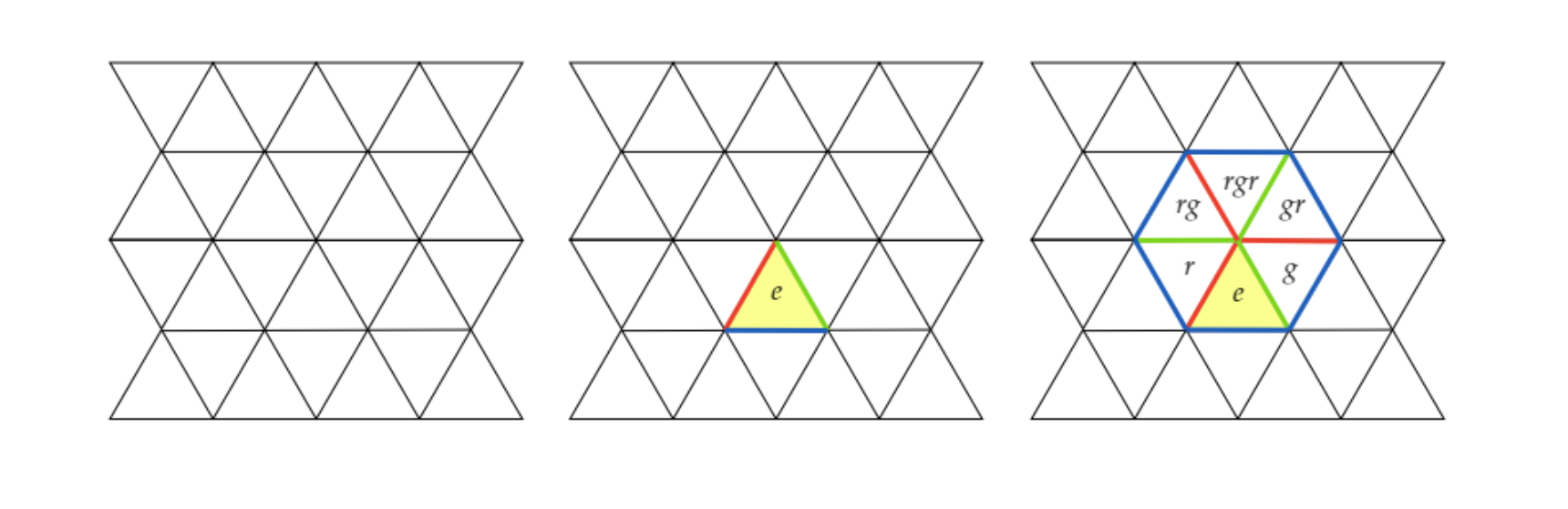}
\end{figure}

One can see the symmetric group $S_3=\{e, g, r, gr, rg, grg\}$ living inside this two dimensional tiling (we call it an $S_3$-polygon). Define arbitrarily  one triangle to be the identity $e\in S_3$ and one of its vertices to be the zero of the vector space (in the middle drawing, the zero is chosen to be the intersection of the green and red edges). Draw each segment of that triangle of a color: red, blue and green (see the image in the middle). 

Call $r$ the reflection of the plane through the line defined by the red edge. Call $g$ and $b$ the same with green and blue. 
By abuse of language I will call $g$ the triangle $g(e)$ and $r$ the triangle $r(e)$ (see the image in the right, where  the edges of the new triangles are painted accordingly). Call $gr$ (resp. $rg, rgr$) the triangle $r(g(e))$ (resp. $g(r(e)), r(g(r(e)))$). Of course we can do the same adding all the blue reflections and obtain the following: 

\begin{figure}[H]
    \centering
    \includegraphics[width=12cm]{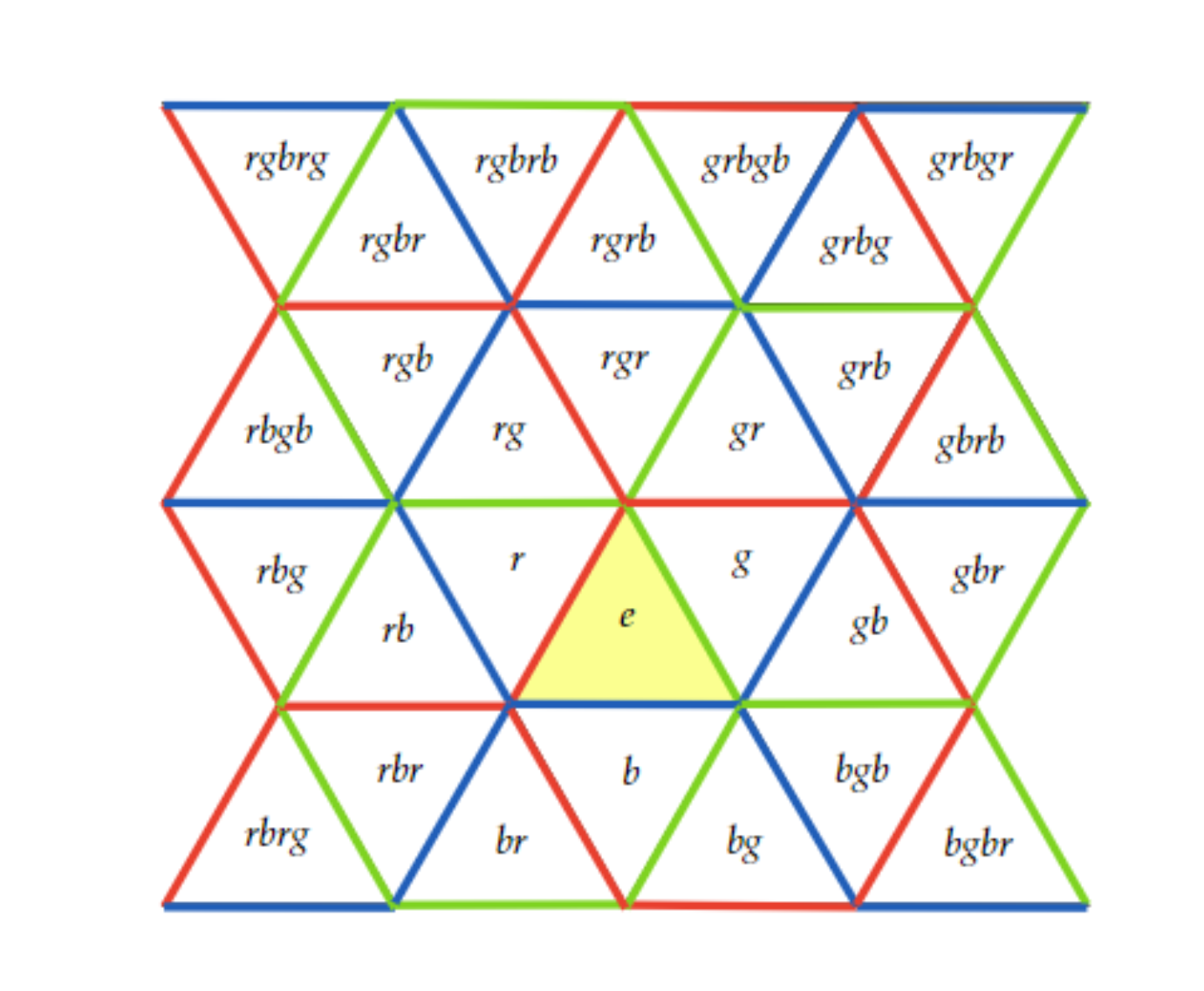}
    \caption{Tiling with equilateral triangles}
\end{figure}

Voilà the affine symmetric group  $\widetilde{S}_3$. The elements of $\widetilde{S}_3$ are exactly the equilateral triangles in the tiling, and multiplication in this group is concatenation, i.e. $$rgb\cdot grbbb:=rgbgrbbb,$$
Of course we are supposing (and it is true) that there will be a triangle named $rgbgrbbb.$ 
Most triangles can be named in several ways. For example $rgr=grg$, but this is not a problem, because $rgr=grg$ as elements of $\mathrm{Aff}(\mathbb{R}^2)$ (the set of affine transformations of $\mathbb{R}^2$)

There are  three particularly nice ways in which $S_3$ sits as a subgroup of $\widetilde{S}_3$, the three  (blue, green and red) hexagons that contain the identity. 

\subsection{Correcting the mistake.}\label{correcting}

In dimension three there is one very nice tiling of the space by a tetrahedron (that we will call \textit{alcove} from now on), it is just not the regular tetrahedron. Let us describe it. As any tetrahedron, it has four faces, so it has six dihedral angles (i.e. the angles that the four faces make pairwise at their lines of intersection). Four of them are $60^{\circ}$ and two of them $90^{\circ}$. All the faces are congruent isosceles triangles, with angles $\mathrm{arccos}(\frac{1}{3})$ and $\mathrm{arccos}(\frac{1}{\sqrt{3}})$, which are approximately  $54.74$ and $70.53$ degrees. This tetrahedron is a special case of a \textit{Disphenoid}, and it is congruent to the tetrahedron with vertices 
$(-1, 0, 0), (1, 0, 0), (0, 1, 1),$ and $(0, 1, -1)$ in $\mathbb{R}^3$.
\begin{figure}[H]
    \centering
    \includegraphics[width=7cm]{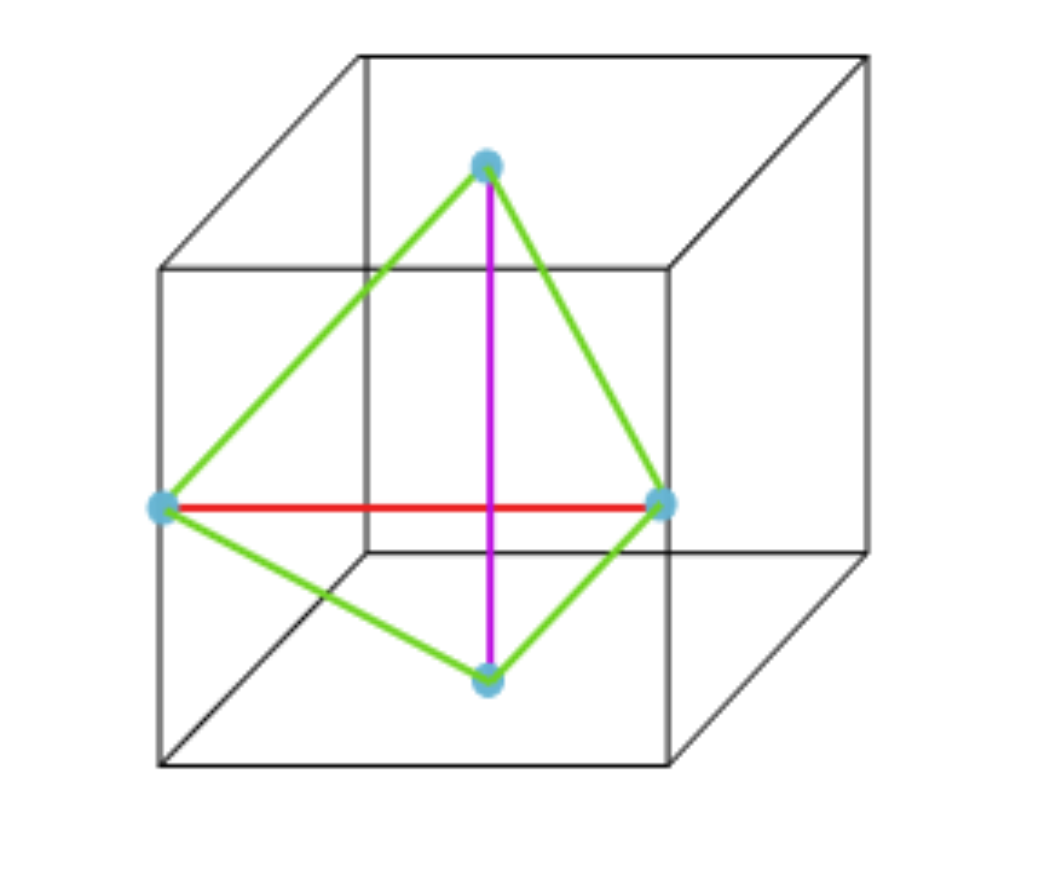}
\end{figure}

But it was really hard for me to understand this tetrahedron until I had it in my hands, that is why I produced an origami for you to print this page and do it yourself. 
\begin{figure}[htp]\label{origami}
    \centering
    \includegraphics[width=12cm]{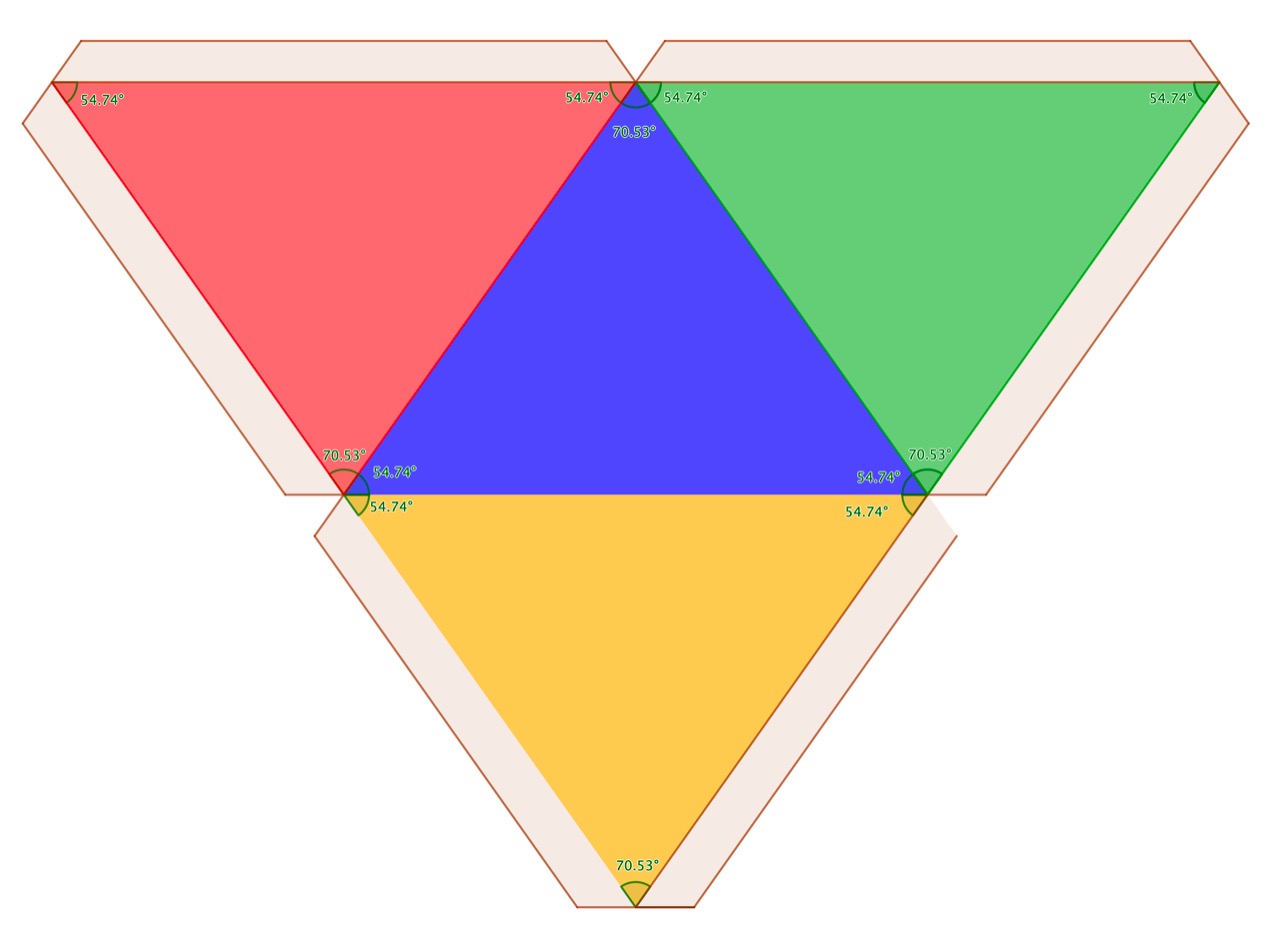}
    \caption{Origami to produce an alcove of $\widetilde{S}_4$}
\end{figure}

Here again one can decide that some element is the identity $e$, color each of the faces of the alcove with one of four colors, and then do all the reflections with respect to three of them, and obtain an $S_4$-polyhedron as $24$ alcoves intersecting in one point. One can think of such an $S_4$-polyhedron as a cube (I drew in pink the vertices of such a cube) such that you put a pyramid (with square basis) in each face. 

\begin{figure}[H]\label{lala}
    \centering
    \includegraphics[width=10cm]{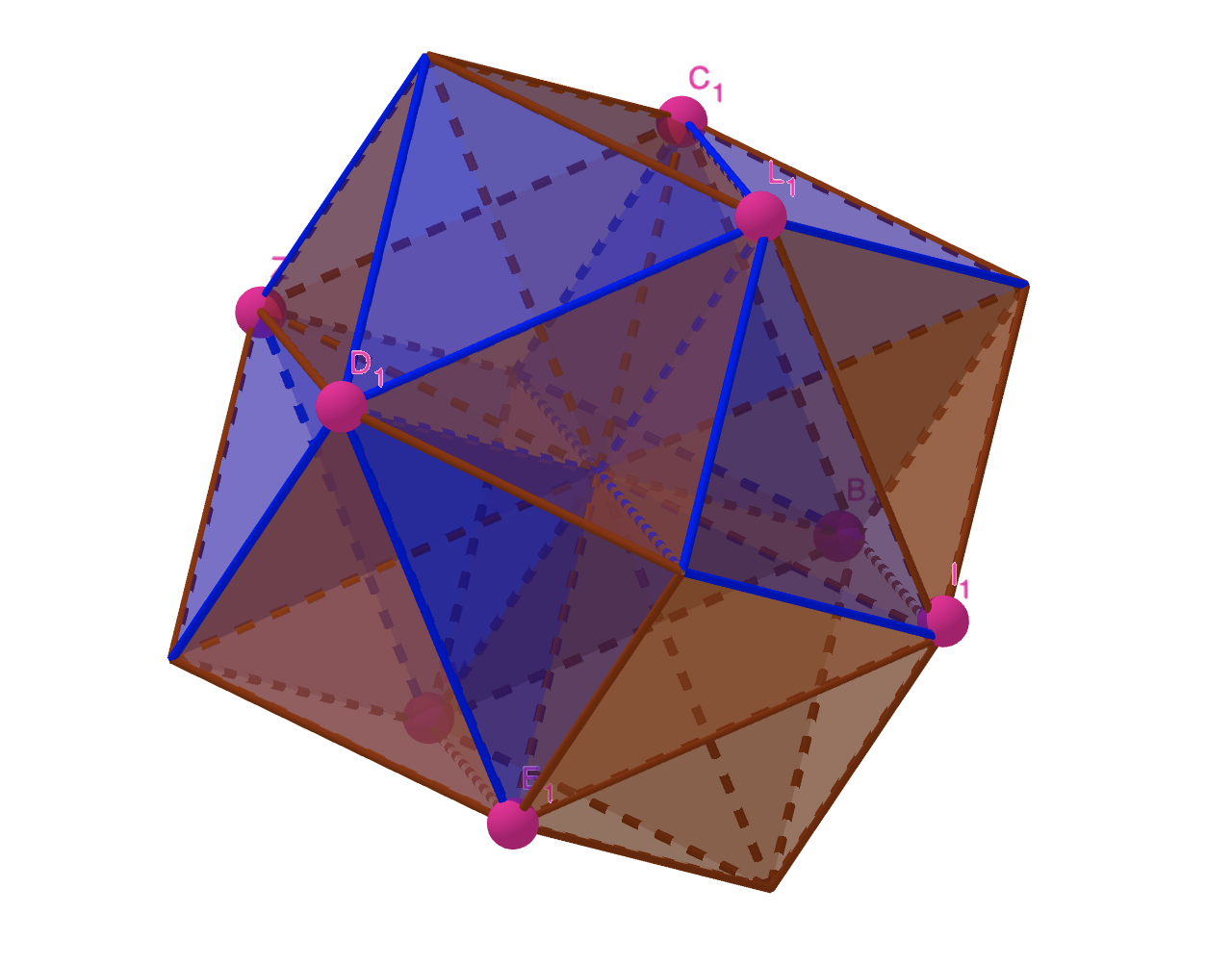}
    \caption{$S_4$ polyhedron}
\end{figure}

\begin{figure}[H]
    \centering
    \includegraphics[width=6cm]{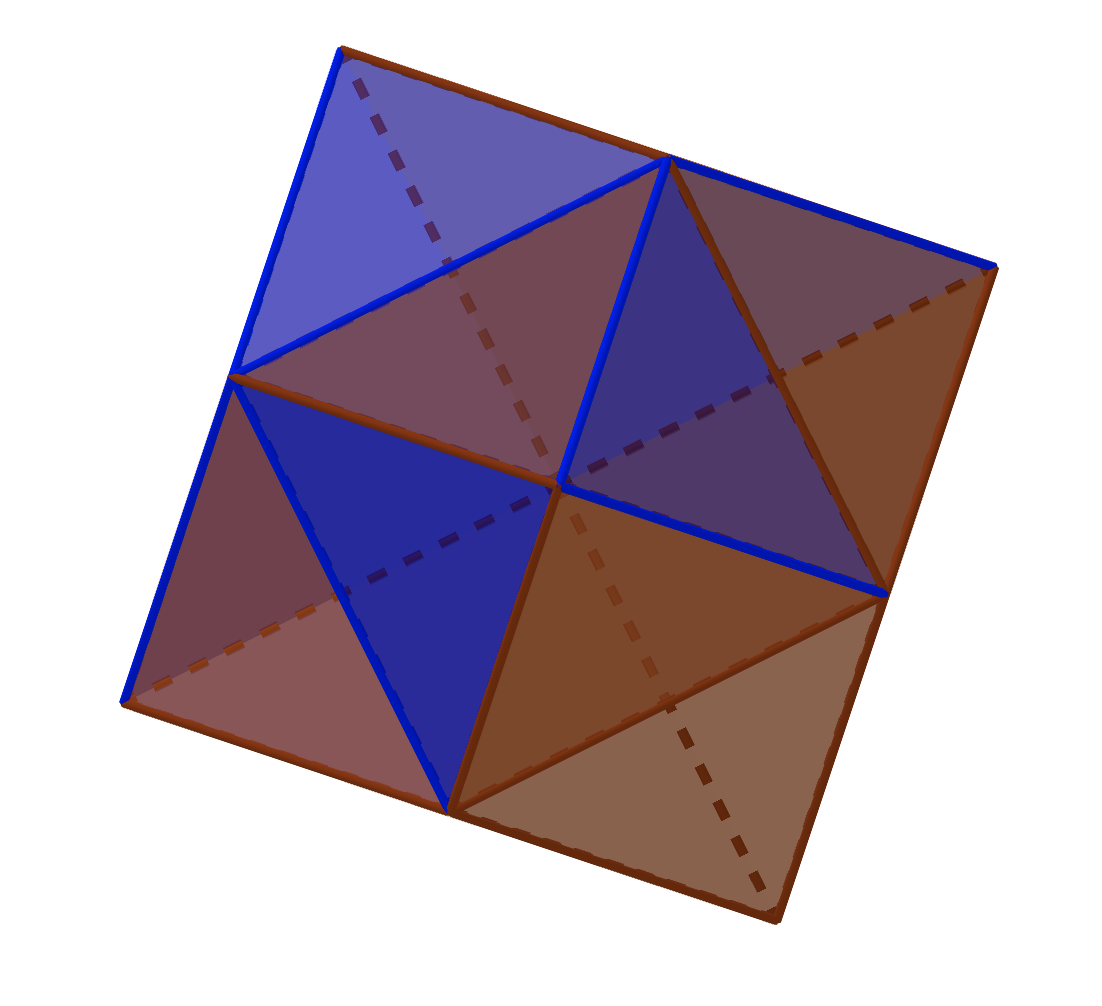} \includegraphics[width=6cm]{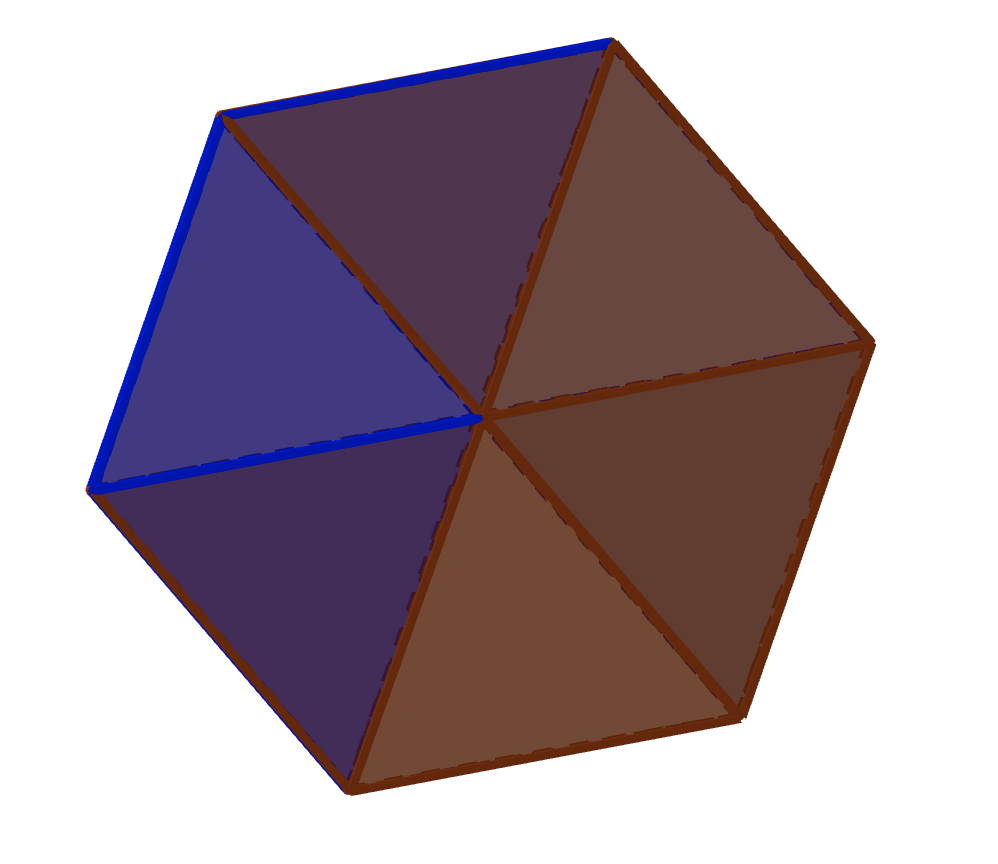}
    \caption{Different views of the $S_4$ polyhedron }
    \label{angles}
\end{figure}

But if one looks at this $S_4$-polyhedron from different angles, one sees these very interesting images. The one on the right of Figure \ref{angles} is part of the two dimensional tessellation that we saw before and the one on the left is the product of 
two tessellations of the same kind in dimension 1 (we will see the general definition of these tessellations in Section \ref{bigger}). 

\begin{figure}[H]
    \centering
    \includegraphics[width=4cm]{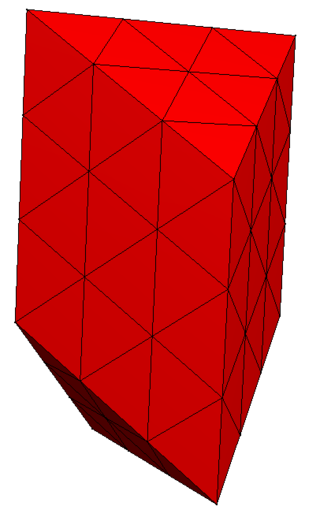}
    \caption{A bigger part of our tessellation}
\end{figure}

I find it beautiful that the dual\footnote{Two polyhedra are called \textit{dual} if the vertices of one correspond to the faces of the other, and the edges of one correspond to the edges of the other, in an incidence-preserving way.} of this $S_4$-polyhedron is the permutohedron (see wikipedia for the definition), that is exactly the polyhedron obtained as the convex hull of the set $\{w\cdot \rho\, \vert \, w\in S_4\}$ (the definition of $\rho$ will be given in Section \ref{bigger}).  If the reader is not familiar with Weyl groups, please skip the next paragraph.

What we said in the last paragraph is beautiful, but easy to prove, also for other Weyl groups. The reason is that a scalar multiple of $\rho$ hits exactly the  centroid  of the ``external'' face (the face not containing $0$) of the alcove corresponding to $w_0$ in the $W$-polyhedron (for $W$ a Weyl group) and the $W$-orbit hits exactly  the centroid of all the  external faces of all the alcoves of the $W$-polyhedron, thus giving the dual by definition.

\subsection{Fancy intermezzo}\label{fancy}\epigraph{\textit{But in my opinion, the most significant reason for the usefulness of perverse sheaves is the following secret known to experts: perverse sheaves are \textit{easy}, in the sense that most arguments come down to a rather short list of tools}}{Pramod Achar}
The reader might ask at this point: ``why should I care about $\widetilde{S}_n$? What is the point of all of this? Is it because we hate Aristotle? Please end my misery and speak up.'' 
I would answer to the reader: wow, you are intense, but don't worry... Explanations are coming.

We have understood in Section \ref{Sn} the representation theory of the symmetric group over the field of complex numbers (it would have been similar over any field of characteristic zero).

\vspace{0.1cm}

\textbf{But what about the representation theory of the symmetric group over a field of positive characteristic? }

\vspace{0.2cm}

This is usually known as the \textit{modular representation theory} of the symmetric group, i.e.  group homomorphisms $S_n\rightarrow GL(V)$ for $V$ a $k$-vector space, where $k$ is a field of positive characteristic. Then things start to get spicy... We know how to classify the irreducible representations but  imagine how poor our understanding is, if we don't even know their dimensions!

We know by now that  $SL_n(\mathbb{F}_q)$ is (a deformation, a partition, a dual) strongly related to the symmetric group, so it is not surprising that one studies its representation theory to understand that of $S_n$. But in 1963 Robert Steinberg changed the problem  into an easier one (although still extremely hard). He proved  \cite{Steinberg} that all the irreducible representations of the finite group $SL_n(\mathbb{F}_q)$ could be obtained from the irreducible algebraic (I will explain  this concept in a minute) representations of the algebraic group $SL_n(\overline{\mathbb{F}}_p)$ by ``restriction''. Here $\overline{\mathbb{F}}_p$ is the algebraic closure of the finite field with $p$ elements $\mathbb{F}_p$: a countably infinite field that contains a copy of the field of order $p^n$ for each positive integer $n$ (and is in fact the union of these copies). This is how we see $SL_n(\mathbb{F}_q)$ as a subset of $SL_n(\overline{\mathbb{F}}_p)$, and thus it makes sense to speak about ``restriction''.

For a field $k$, I will explain what is a (finite dimensional\footnote{By a basic result of the theory, all representations are direct
limits of finite-dimensional representations. That is why I assume our representations to be
finite-dimensional.}) algebraic representation of $SL_2(k)$, and the reader will correctly guess what  it is for $SL_n(k)$. Recall that elements of $SL_2(k)$ are matrices $\left ( \begin{matrix} a & b \\ c & d \end{matrix} \right )$ with determinant equal to $1.$ An \textit{algebraic representation} of $SL_2(k)$ is a group homomorphism 
$$SL_2(k)\rightarrow GL_n(k),$$
where each entry in $GL_n(k)$ is a polynomial in $a, b, c$ and $d$ (as before, we call $n$ the \textit{dimension} of the representation). 

Here an example of a representation of dimension $3$:
\begin{align*}
\left ( \begin{matrix} a & b \\ c & d \end{matrix} \right ) & \mapsto
\left ( \begin{matrix} a^2 & ab & b^2  \\ 2ac & ad + bc & 2bd \\ c^2 & cd & d^2 \end{matrix} \right )
\end{align*}
Here an example of a representation of dimension $2$, with $\mathrm{char}(k)=3:$
\[
\left ( \begin{matrix} a & b \\ c & d \end{matrix} \right )   \mapsto 
\left ( \begin{matrix} a^3 & b^3 \\ c^3 & d^3 \end{matrix} \right )  
\]

Why did I say that the new problem (understanding algebraic representations of  $SL_n(\overline{\mathbb{F}}_p)$) is easier? Because the extra structure
that comes from considering $SL_n(\overline{\mathbb{F}}_p)$ as an algebraic group
is  immensely useful.

Now let us go to business, and answer the question. Why do we care about $\widetilde{S}_n$? Because it is the group behind the representation theory of $SL_n(\overline{\mathbb{F}}_p)$. We will see this in two different ways, the first one relating it to a diagrammatic category, and the second one relating it with a geometric category. 

The following characterizations of the algebraic representations of $SL_n(\overline{\mathbb{F}}_p)$ are probably among the top five most important contributions to representation theory in the last 20 years (and both are for general reductive groups, not only for $SL_n(\overline{\mathbb{F}}_p)$). 
For application (1) I will assume that the reader knows what a category  and what a Coxeter system are (if you don't, just read those concepts in wikipedia or skip the rest of this section without  guilt). For part (2) the reader needs to know what an algebraic variety is (and for the footnote much more). Let me say that this is quite fancy, so don't get discouraged if you don't understand every word or every idea.

\begin{enumerate}
    \item We will see in Section \ref{bigger} that there is a subset $S\cup \{s_0\}\subseteq \widetilde{S}_n$ that gives $\widetilde{S}_n$ the structure of a Coxeter system. In Section \ref{sbim} we will see that to the Coxeter system $(\widetilde{S}_n,S\cup \{s_0\})$ (in fact, to any Coxeter system $(W,S)$) one can associate  a category $\mathcal{H}_{\mathrm{BS}}(\widetilde{S}_n)$ called \textit{the Bott-Samelson Hecke category} where objects      are arbitrary sequences $(r_1,r_2,\ldots,r_m)$, with $r_i \in S\cup \{s_0\}$. The good thing about this category is that it is quite easy to compute stuff, as we will see in Section \ref{bigger}.
    One defines the \textit{Bott-Samelson anti-spherical category} $\mathcal{N}_{\mathrm{BS}}(\widetilde{S}_n)$ as the following quotient category of the Bott-Samelson Hecke category: the objects of $\mathcal{N}_{\mathrm{BS}}(\widetilde{S}_n)$ are the same as those of $\mathcal{H}_{\mathrm{BS}}(\widetilde{S}_n)$. A morphism is zero in $\mathcal{N}_{\mathrm{BS}}(\widetilde{S}_n)$ if and only if it factors through an object $(r_1,r_2,\ldots,r_m)$ with $r_1\neq s_0.$
    
In 2018, Simon Riche and Geordie Williamson proved \cite{GeordieSimon}, very  roughly speaking,  that the category $\mathcal{N}_{\mathrm{BS}}(\widetilde{S}_n)$ is equivalent to the category of algebraic representations of $GL_n(\overline{\mathbb{F}}_p)$. (For more general groups it has been proved this year, independently by Simon Riche and Roman Bezrukavnikov \cite{RicheB} and by Joshua Ciappara \cite{Ciappara}.) This is not only an equivalence of categories, it preserves much more structure, but one can only say it using fancy words: it is an equivalence of additive right $\mathcal{H}_{\mathrm{BS}}(\widetilde{S}_n)$-module categories. For the precise statement of the theorem, see \cite[Th. 1.2]{GeordieSimon}.
    \item The second theorem involving $\widetilde{S}_n$ is called the \textit{geometric Satake equivalence} proved by Ivan Mirkovic and Kari Vilonen \cite{MV} in 2007. It gives an equivalence between the category of algebraic representations of $GL_n(\overline{\mathbb{F}}_p)$ and a geometric category that I will not explain in detail, but will  try to briefly describe. 
    \begin{itemize}
        \item  The first ingredient for this equivalence is something called the ``affine Grassmannian'' associated to $GL_n$. Let $\mathbb{O}=\mathbb{C}[[x]]$ be the ring of formal powers in $x$. This is the ring of ``infinite polynomials'', with elements of the form   $\sum_{i=0}^{\infty}\lambda_ix^i$ (with $\lambda_i\in \mathbb{C}$). Let $\mathbb{K}=\mathbb{C}((x))$ be its fraction field, the field of formal Laurent series, with elements of the form   $\sum_{i=-\infty}^{\infty}\lambda_ix^i.$
    As a set, the \textit{affine Grassmannian}  is $$\mathcal{G}r=GL_n(\mathbb{K})/GL_n(\mathbb{O}),$$
    It is not an algebraic variety itself but one can prove that it is the increasing union of some very concrete algebraic varieties (this induces a topology on $\mathcal{G}r$). There is a group homomorphism $$GL_n(\mathbb{O})\rightarrow GL_n(\mathbb{C})$$ given by $x\mapsto 0.$ Define  an \textit{Iwahori subgroup} $I$ of $GL_n(\mathbb{K})$ as the preimage of the upper triangular matrices under this map. There is an easy map (that I will not describe) $(\mathbb{Z}/n\mathbb{Z})\ltimes \widetilde{S}_n\rightarrow GL_n(\mathbb{K})$ that we denote $w\mapsto \dot{w}$ and a generalization of the Bruhat decomposition (see Equation \eqref{bd}) is given by
    $$GL_n(\mathbb{K})=\coprod_{w\in (\mathbb{Z}/n\mathbb{Z})\ltimes \widetilde{S}_n}I\dot{w}I$$
    We see that the symmetric group $S_n$ in Equation \eqref{bd} is replaced here by the affine symmetric group $\widetilde{S}_n$ (allow me to ignore the cyclic group $\mathbb{Z}/n\mathbb{Z})$). 
    \item The second main ingredient are ``perverse sheaves''. This is hard to define but not so hard to work with, as Pramod Achar told us in the epigraph. The category of perverse sheaves on a variety (which is almost our case here) is a very good category to compute stuff in, because it has nearly everything that one would ask\footnote{It is an abelian category, Noetherian, Artinian, it is the heart of a t-structure, the simple objects are very easy to compute, etc.} of a well behaved category.
        \end{itemize}
Now we can say the exact formulation of the geometric Satake equivalence. 
It says that the category of spherical\footnote{Let us ignore this word for now.} perverse $\overline{\mathbb{F}}_p$-sheaves on the complex affine Grassmannian associated to $GL_n$  is equivalent\footnote{Something very important about this equivalence and that is being ignored here, is that the general equivalence is between perverse sheaves on the affine Grassmannian associated to a group $G$ and algebraic representations of ${}^LG,$ the Langlands dual group. We don't see this phenomenon here, because  $GL_n$ is self-dual.} to the category of algebraic representations of $GL_n(\overline{\mathbb{F}}_p)$. This equivalence is also enriched, in the sense that one can define a convolution on perverse sheaves, and it corresponds to tensor product on representations. So the point is that this gives a fresh, geometric look on algebraic representations of $GL_n(\overline{\mathbb{F}}_p)$.  
\end{enumerate}

\subsection{The bigger picture}\label{bigger}
The two tessellations that we have studied, are special cases of a family of tessellations parametrized by a positive integer $n$. Consider $V$ the hyperplane of $\mathbb{R}^{n+1}$ consisting of points with the sum of the coordinates equal to zero. The group $S_{n+1}$ acts on $\mathbb{R}^{n+1}$ by permuting the variables, and so it acts on $V.$ Let $\varepsilon_i\in \mathbb{R}^{n+1}$ be the $i^{\mathrm{th}}$ coordinate vector (one in position $i$ and zero elsewhere). 

Define the \textit{root system} as the set of vectors (or \textit{roots})
$$\Phi:=\{\varepsilon_i-\varepsilon_j\, :\, i\neq j, 1\leq i,j\leq n+1\},$$
and for $\alpha\in \Phi$ and $i\in \mathbb{Z}$ define  $$H_{\alpha, i}=\{v\in V\,:\, \langle v,\alpha \rangle=i \}.$$

\subsubsection{First definition: reflections} The \textit{affine Weyl group} $\widetilde{S}_n$ is the group of affine transformations of $V$ generated by the orthogonal reflections with respect to the $H_{\alpha, i}$. To be more precise, this reflection is defined by the formula $$ s_{\alpha, i}(v)=v-(\langle\alpha,v \rangle-i)\alpha.$$ This definition is elegant but not so useful in practice. 

\subsubsection{Second definition: semi-direct product}  Consider $Q\subset \mathbb{R}^{n+1} $ to be the set of vectors with integer coordinates in $V$ (i.e. with sum of the coordinates equal to zero). We see elements of $Q$ as affine transformations of $V$, given that each vector $v\in V$ determines a translation by $v$. Recall that $S_{n+1}$ acts on $\mathbb{R}^{n+1}$ by permuting the variables, so it acts on $V$ and on $Q.$ Then $\widetilde{S}_{n+1}=S_{n+1}\ltimes Q$.

\subsubsection{Alcoves}
Let $H$ be the union of all $H_{\alpha, i}$.
Each connected component of $H^c$ (the complement of $H$ in $V$) is called an \textit{alcove}. All alcoves are  congruent to each other, i.e. they are the same open $n$-simplex (the $n$-analogue of the interior of  a tetrahedron)  translated and rotated. In the case $n=2$ this gives our first tessellation (Section \ref{why}) and if $n=3$ this gives our second tessellation (Section \ref{correcting}).

Look at Equation \eqref{omega} for the definition of $\omega_i$. Let us define the alcove $A_{\mathrm{id}}$ whose $n+1$ vertices are $$\{0\}\cup \{-\omega_i\}_{1\leq i\leq n}$$ (one can check that this is indeed an alcove, i.e. a connected component of $H^c$).

One can prove that $\widetilde{S}_{n+1}$ acts simply transitively on the set of alcoves, and what this means in English is that  there is a bijection between $\widetilde{S}_{n+1}$ and the set of alcoves via $w\mapsto w(A_{\mathrm{id}})$ (in particular, $A_{\mathrm{id}}$ corresponds to the identity element of the group). 

\subsubsection{Third definition: gen and rel.}
One can prove that the orthogonal reflections through the hyperplanes supporting the faces of $A_{\mathrm{id}}$ generate $\widetilde{S}_{n+1}$. More precisely, for $1\leq i\leq n$, let $s_i$ be the orthogonal reflection through  the only face not containing $-\omega_i$ and let $s_0$ be the reflection through the remaining face of $A_{\mathrm{id}}$.

We use the convention that $s_{n+1}:=s_0$
Then  $\widetilde{S}_{n+1}$ has a  presentation given by generators $\{s_0, s_1, s_2, \ldots, s_{n}\}$ and relations:
\begin{enumerate}
    \item $s_i^2=1,$
    \item $s_is_j=s_js_i$ if $\vert i-j\vert>1,$
    \item $(s_is_{i+1})^3=1.$
\end{enumerate}
Let me remark that if one eliminates $s_0$ from the generators, the same relations define the group $S_{n+1}$. The group $\widetilde{S}_{n+1}$ is infinite while $S_{n+1}$ is finite. These two groups are best understood in their natural habitat: they are Coxeter groups.

\subsubsection{Coxeter groups}\label{Coxeter}\epigraph{\textit{In our times, geometers are still exploring those new Wonderlands, partly for the sake of their applications to cosmology and other branches of science, but much more for the sheer joy of passing through the looking glass into a land where the familiar lines, planes, triangles, circles and spheres are seen to behave in strange but precisely determined ways.}}{Harold Scott MacDonald Coxeter }  

Some days ago I stumbled into this majestic \Laughey[1.4]  introduction  to Coxeter groups \cite{Gentle} (and some things on their relation to non-euclidean geometries). Here I will be very brief and give some basic definitions. 

A \textit{Coxeter system} $(W,S)$ is a group $W$ together with a generating subset $S\subseteq W$ such that $W$ admits a presentation by generators $s\in S$ and relations $(sr)^{m(s,r)}=1$ for $s,r\in S$, with $m(s,s)=1$ and $m(s,r)\geq 2$ (and potentially $m(s,r)=\infty$) for $s\neq r.$ The group $W$ is called a \textit{Coxeter group} and the set $S$ are called the \textit{simple reflections}.

For $x\in W$, if $n$ is the minimal integer such that $x$ can be written as a product of $n$ simple reflections, then $n$ is the \textit{length} of $x$ and is denoted by $l(x)=n.$ In this case, any expression $x=s_1s_2\cdots s_n$ with $s_i\in S$ is called a \textit{reduced expression}.

Fix a reduced expression of some element $x=s_1s_2\cdots s_n$. An element $y$ is called \textit{lesser in the Bruhat order} and is denoted by $y\leq x$ if $y=s_{a}s_b\cdots s_c$, where $1\leq a\leq b\leq \cdots \leq c\leq n.$

\subsection{Géométrie alcovique et géométrie euclidienne}\label{gage}\footnote{This is a humble reference to the famous GAGA paper (géométrie algebrique et géométrie analytique) \cite{Serre} by Jean-Pierre Serre.} 

Something that has had little significance in the history of humanity, but  I believe  is a powerful and extremely beautiful source of insights, is the relation between Euclidean geometry and alcovic geometry. I will give a (conjectural) example that I would like to prove in the near future with my collaborators Damian de la Fuente and David Plaza. 

There are two important bases of $V$. The first one consists of the \textit{simple roots} $\alpha_i:=\varepsilon_i-\varepsilon_{i+1}$ for $1\leq i\leq n,$ and the other one, consists  of the  \textit{fundamental weights}   \begin{equation}\label{omega}
    \omega_i:= (\varepsilon_1+\cdots+\varepsilon_i)-\frac{i}{n+1}\sum_{j=1}^{n+1}\varepsilon_j,
\end{equation} for $1\leq i\leq n.$
They are dual to each other with respect to the dot product in $\mathbb{R}^{n+1}$, i.e. $$\alpha_i\cdot \omega_j= \delta_{i,j}$$ ($\delta_{i,j}$ is $1$ if $i=j$ and zero otherwise). A vector of the form $\lambda=\lambda_1\omega_1+\cdots+\lambda_n\omega_n$ with all $\lambda_i\in \mathbb{Z}_{\geq 0}$ is called a \textit{dominant weight}. An important vector is $$\rho=\omega_1+\omega_2+\cdots +\omega_n.$$

For a dominant weight $\lambda=\lambda_1\omega_1+\cdots+\lambda_n\omega_n$, define $\theta(\lambda)\in \widetilde{S}_n$ as the unique alcove containing $\lambda +\epsilon\rho$, where $\epsilon,$ as usual, is a very small real number (the  $\theta(\lambda)$'s are very important in representation theory). One would like to understand the set $$\leq \theta(\lambda):=\{w\in \widetilde{S}_n : w\leq \theta(\lambda)\}.$$
Waldeck Sch\"{u}tzer \cite{Schutzer} proved that the cardinality of that set, $\mathrm{card}(\leq \theta(\lambda))$ is a polynomial of degree $n$ in the $n$ variables $\lambda_1,\ldots, \lambda_n.$ For example, when $n=2$ and $\lambda=\lambda_1\omega_1+\lambda_2\omega_2,$ then
$$\mathrm{card}(\leq \theta(\lambda))= (3\lambda_1^2 + 3\lambda_2^
2 + 12\lambda_1\lambda_2) + (9\lambda_1 + 9\lambda_2) + (6).$$ I put the parentheses only to distinguish the graded parts of this polynomial of degrees $2, 1,$ and $0.$ For $n=3$ it is already quite a big polynomial\footnote{If the reader is tempted to believe that the sets $\leq \theta(w)$ with $w\in W_a$ are easy, she can check the paper \cite{Bjorner}  by Anders Bj\"{o}rner  and Torsten Ekedahl that appeared in Annals of Mathematics in 2009, where the main theorem (proved with fancy mathematics) is that if $f_i$ is the number of elements in $\leq \theta(w)$ of length $i$ then $f_i\leq f_j$ if $0\leq i<j\leq l(w)-i.$ }. Notice that $\mathrm{card}(\leq \theta(\lambda))$ counts the number of alcoves in the set $\leq \theta(\lambda)$,  so
$$\mathrm{card}(\leq \theta(\lambda))=\frac{\mathrm{Vol}(\leq \theta(\lambda))}{\mathrm{Vol}(A_{\mathrm{id}})} $$
where Vol stands for the volume with respect to  the standard metric on $V$ (induced form that of $R^{n+1}$. By the way, $\mathrm{Vol}(A_{\mathrm{id}})$ is just $\sqrt{n+1}$).

We have a conjectural geometric (in the sense of Euclidean geometry) interpretation of this polynomial. We produce sets $P_i(\lambda)\subseteq V$ such that  $$\leq \theta(\lambda)=\coprod_{0\leq i\leq n}P_i(\lambda),$$ 
i.e. $\leq \theta(\lambda)$ is partitioned into a disjoint union of sets satisfying that 
$$\frac{\mathrm{Vol}(P_i(\lambda))}{\mathrm{Vol}(A_{\mathrm{id}})} $$
is the degree $i$ part of the polynomial $\mathrm{card}(\leq \theta(\lambda))$. I don't want to give the general recipe to produce the $P_i(\lambda)$, but let me show an illustrative example. 

\begin{figure}[H]\label{Pi}
    \centering
    \includegraphics[width=15cm]{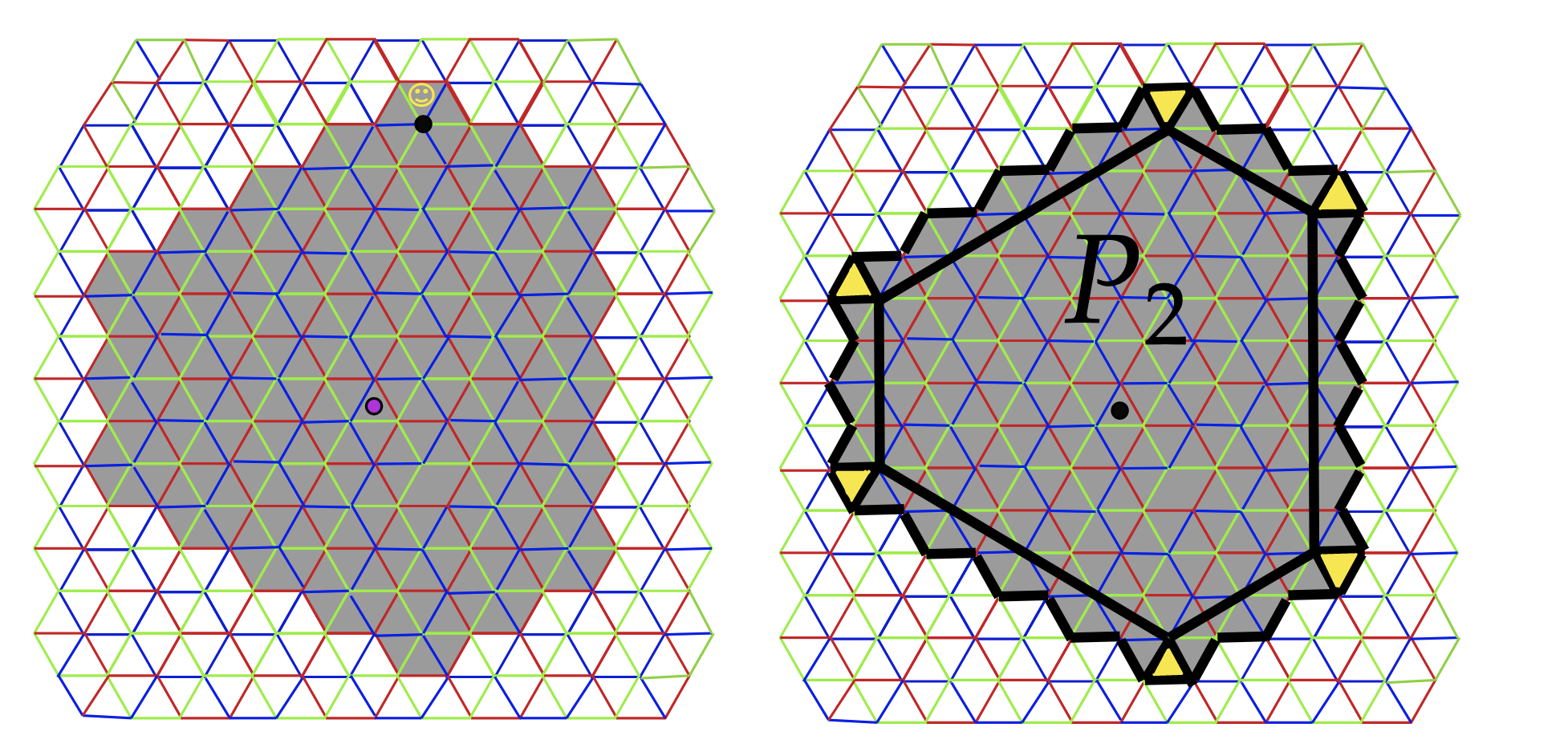}
    \caption{The partition of $\theta(\lambda)$}
\end{figure}
In the left drawing, the purple dot is in the identity alcove, the dark dot is $\lambda \in V$, the happy face is in the alcove $\theta(\lambda)$ and the grey region is $\leq \theta(\lambda)$. In the right picture, the big hexagon is $P_2,$ the six yellow alcoves together form $P_0$ and the rest of $\leq \theta(\lambda)$ is $P_1.$ It is a general thing that $P_n$ is the convex hull of the set $S_n\cdot \lambda.$

\part{Quantum level}

\section{Iwahori-Hecke algebras}\label{IHalg}
\subsection{Deforming the symmetric group}
\epigraph{\textit{One can learn a lot about a mathematical object by studying how it behaves
under small perturbations.}}{Barry Mazur}
In this part of the paper we aim to understand what  ``quantum Schur-Weyl duality'' should mean (recall that this is the guiding principle throughout this paper). Having this in mind, our first ingredient will be the famous Iwahori-Hecke algebra. But let us not rush into it.

 In Section \ref{phantom} we saw that $PGL_n(\mathbb{F}_q)$ is a deformation of $W=S_n$. This means that there is a parameter $q$ and  a family $W_q$ of groups such that 
 \begin{equation}\label{lim}
     \lim_{q \to 1}W_q = W.
 \end{equation} But this is not fully satisfying for two reasons. Firstly $W$ is a Coxeter group and $PGL_n(\mathbb{F}_q)$ is not one (there is a classification of finite Coxeter groups and $PGL_n(\mathbb{F}_q)$ does not belong to it\footnote{This argument is like killing a fly with a bazooka, but I don't know an easier argument.}. Oddly enough, $PGL_2(\mathbb{Z})$ is an infinite Coxeter group). Secondly, $q$ is some prime number to the power of some natural number, and powers of primes are quite rare  (they don't account for ``small perturbations''). We would like that $q$ could be any  complex number.
 
 But these complaints are not easy to handle by the discrete nature of Coxeter groups and the continuous nature of $\mathbb{C}$. In plain terms, the set of Coxeter systems is countable and $\mathbb{C}$ is uncountable, so we would be forced to have  $W_q=W_{q'}$ for a lot of $q,q'$. Moreover a Coxeter group is determined by its $m$-matrix, which is a set of natural numbers. So for Equation \eqref{lim} to be satisfied, we would be forced to have a ball $B$ containing $1$ such that if $q,q'$ are in $B$ then  $W_q=W_{q'}$. This is far from the ``small perturbations''
 philosophy.
 
 But representation theory comes to our aid. Recall from Section \ref{Duality} the group algebra $\mathbb{C}[S_n]$. For a representation theorist this algebra and $S_n$ are ``almost the same''. There is a canonical bijection between representations of $S_n$ over $\mathbb{C}$ and $\mathbb{C}[S_n]$-modules. And this bijection identifies the following concepts (that I have not necessarily explained, they are just to convince the reader  about the ``almost the same'' claim):
 
 
\begin{center}
 \begin{tabular}{l|l}
\hline
Finite dimensional reps. of $S_n$ over $\mathbb{C}$&  Finitely generated non-zero $\mathbb{C}[S_n]$-modules\\
Subrepresentations   & Submodules \\ 
 Irreducible reps. & Simple modules \\ 
Tensor product of representations & Tensor product of modules \\ 

etc.    & etc. \\ 
\hline
\end{tabular}
 \end{center}
 
  \vspace{0.5cm}
 
 So the genius idea is to stop trying to deform the undeformable  Coxeter group $W$ and try instead to deform the algebra $\mathbb{C}[S_n]$, which is much easier. We know that this algebra can be described as the $\mathbb{C}$-algebra with generators $\mathbf{s}$ for all $s\in S$ and relations: 
 \begin{itemize}
     \item \textbf{Quadratic relations:} $\mathbf{s}^2=1$
     \item \textbf{Braid relations:} $\mathbf{s}\mathbf{r}\mathbf{s}\cdots= \mathbf{r}\mathbf{s}\mathbf{r}\cdots $ where each side has $m_{sr}$ factors and $s,r$ belong to $S$.
 \end{itemize}
To deform this algebra is easy! One can, for example, just multiply by $q$ any part of any of the relations, or add $(q-1)$ times anything, and we would have that the resulting set of algebras would converge\footnote{Although with this procedure the family might not be ``flat'' i.e. writing down an arbitrary definition would produce something that has smaller dimension in general. It is a bit of a miracle that this doesn't happen for Coxeter groups.} to $\mathbb{C}[S_n]$ when $q$ tends to $1.$ For example, one could deform the braid relation like this: $$\mathbf{s}\mathbf{r}\mathbf{s}\cdots= q\mathbf{r}\mathbf{s}\mathbf{r}\cdots +(q-1)\mathbf{s}. $$
We will make an interlude to give heuristics on why we shouldn't  deform the braid relations but instead we should deform the quadratic relations. 
\subsection{Interlude: braid groups}\label{braidgroups}\epigraph{\textit{In these days the angel of topology and the devil of abstract algebra fight for the soul of every individual discipline of mathematics.}}{Hermann Klaus Hugo Weyl}

\vspace{0.2cm}

Sometimes you are using a mental image to understand some algebraic object and... {\color{red} paf! the mental image comes to life:  it becomes a \textbf{{\color{violet}topological object}}.}

\vspace{0.2cm}

What is your favorite way to imagine an element of the symmetric group? Maybe as a product of simple transpositions, or as a product of cycles? Maybe as a signed rotation of the $n$-simplex or with the one-line notation $w(1)w(2)\cdots w(n)$. The way in which you imagine it is extremely relevant, as we will see. There is one way of depicting an element of the symmetric group called the \textit{strand diagram notation}. For example this diagram 
\begin{figure}[H]
    \centering
    \includegraphics[width=4cm]{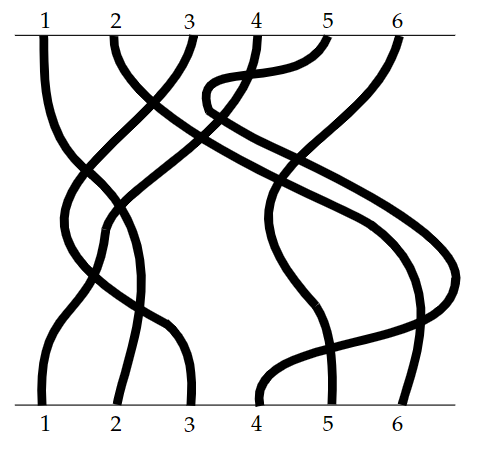}
\end{figure}
 \noindent represents the map sending $1\mapsto 4, 2\mapsto 1, 3\mapsto 3, 4\mapsto 5, 5\mapsto 6, 6\mapsto 2 $ (we read from  bottom to  top). And  we will now give life to this mental image and transform it into a topological object. This will be the same process of thought that we will need to do in order to go from  Soergel bimodules to the diagrammatic Hecke category, which is probably the single most important idea in the waterfall of discoveries in the last ten years mentioned in page 2 of \cite{Gentle}.

We imagine that the drawing before lives in three dimensional space and that each strand that we draw is a strand in $\mathbb{R}^3$ not intersecting the other ones. 

\begin{figure}[H]
    \centering
    \includegraphics[width=5cm]{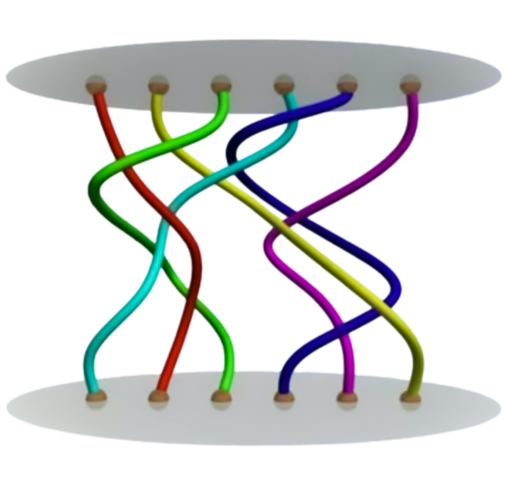}
\end{figure}

(Strands are  one dimensional and have no colors. The last picture is only for illustrative purposes). 
Imagine that the bottom of the diagram belongs to the set $\{(x,y,z): z=0, y=0\}$ and the top  belongs to $\{(x,y,z): z=1, y=0\}$.
We also suppose that these strands are monotone in their $z$-coordinate (one does not allow a strand to go up and  down and  up again). These objects are called \textit{geometric braids in $n$ strands.}

We also need to suppose that the intuition of what a ``real life'' braid is, takes place. 
Two geometric braids $b$ and $b'$ on $n$ strands are isotopic if $b$ can be continuously deformed into $b'$ in the class of geometric braids. A \textit{braid} is a class of geometric braids modulo isotopy. 
As an example, all of these geometric braids on 4 strands are  the same braid: 

\begin{figure}[H]
    \centering
    \includegraphics[width=11cm]{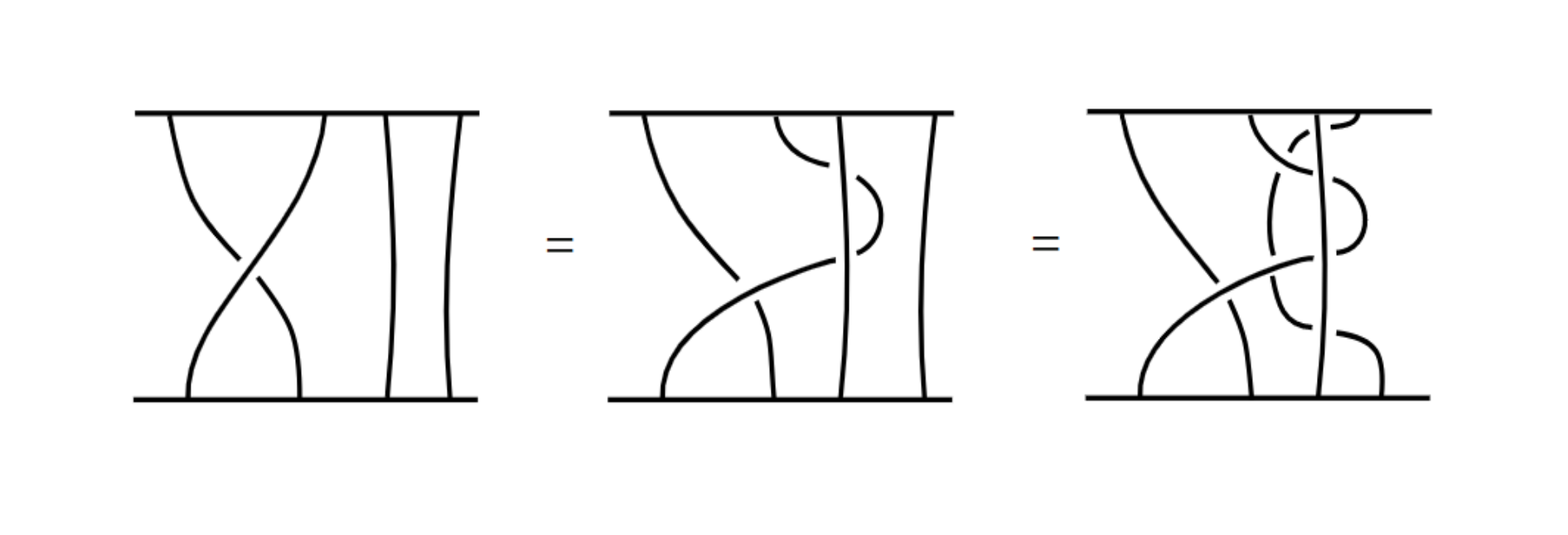}
\end{figure}
Finally,  the composition of braids is quite simple: just put one diagram on top of the other one (the top of the diagram will now live in $\{(x,y,z): z=2, y=0\}$) and then ``compress it''  in order for the composition to live again between $\{(x,y,z): z=0, y=0\}$ and $\{(x,y,z): z=1, y=0\}.$ For a precise definition of how to ``compress'' the diagrams,  more details on the construction of the braid group and everything that I will mention in this section, I strongly recommend the beautiful book by Christian Kassel and Vladimir Turaev \cite{Kassel}. Another  nice book  is \cite{Kamada}.

The set of braids on $n$ strands with the operation described before is called the \textit{braid group} $B_n$. It is the third most beautiful group in the Universe, in my opinion. It appears everywhere!
There are several other equivalent (and quite diverse) descriptions of this group: as braid automorphisms of the free group, as mapping class groups (this interpretation is so incredibly beautiful! And it connects this theory with William Thurston's classification theorem of homeomorphisms of  compact orientable surfaces), as fundamental groups of configuration spaces, etc. These  interpretations might sound difficult, but they are rather easy to understand, very  visual and well explained in the book mentioned above. 

One version of this group that is particularly easy is the following one. Let $(W,S)$ be a Coxeter system. Then the \textit{Artin braid group} is 

$$B(W)=\langle \sigma_s, s\in S \ \vert \  \underbrace{\sigma_s\sigma_r\sigma_s\cdots}_{m_{sr}}= \underbrace{\sigma_r\sigma_s\sigma_r\ldots \rangle}_{m_{sr}}.$$
One can prove that  $$B_n \cong B(S_n).$$

So one can see clearly that $S_n$ is a quotient of $B_n$ and not a subgroup as one could have thought. Geometrically the quotient is quite clear, given that a geometric braid permutes the sets in the bottom and in the top, both sets being in canonical bijection with $\{1,2,3,\ldots,n\}$. 

One last amazing theorem  before this interlude ends, is that of James Waddell Alexander II. In 1923 he proved \cite{Alexander} that every knot (or even more generally, every \textit{link} which is a collection of non-intersecting knots) can be obtained ``closing'' a braid as in the following picture. 

\begin{figure}[H]
    \centering
    \includegraphics[width=7cm]{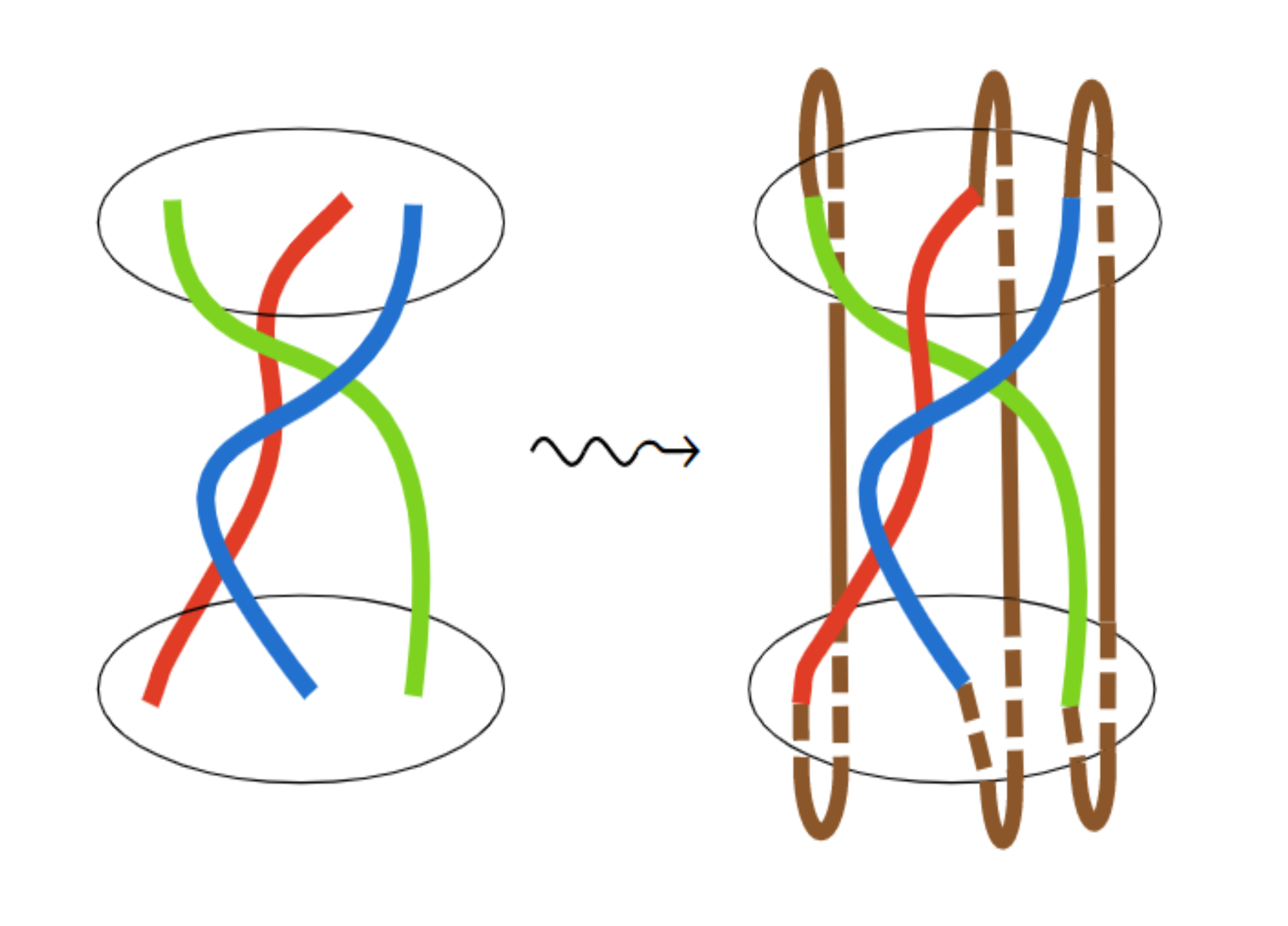}
\end{figure}

In 1985 \cite{Jones} Vaughan Jones, while working on operator algebras, discovered new representations of the braid groups, from which he derived his celebrated polynomial of knots and links (for the discovery of the Jones polynomial, he won the Fields Medal). In the same paper he realised that passing through the Hecke algebra was the best idea... So, kind reader, please ask yourself: do you still want to deform the braid relation?

\subsection{Weil, Shimura, Hecke, Iwahori...}\label{Weil}
Before we decide what is ``the correct deformation of a Coxeter group'' let me tell you a story. 

This story starts in  a land far, far away, called France. At some point André Weil\footnote{There is an old controversy in the mathematical community on whether André Weil and Andrew Wiles are the same person or not. Some argue it is just the pronunciation of the same name in French and English. Some go even further and suggest that André Weil's sister, famous philosopher Simone Weil is the same person as Andrew Wiles's sister,  the most decorated gymnast of all time: Simone Wiles. 
On such a delicate matter I prefer not to pronounce myself (just to be clear, this footnote is a joke). }  had an idea. He told it to Goro Shimura  \cite[p304]{Shimura}. Because I read the lovely 
six-pages article \cite{Shimura2}, I  imagine that it was a rainy day and that they were in Weil's favorite
restaurant  \textit{Au Vieux Paris}. They were probably having radish with buttered rabbit.    Weil's idea was the following.

Let $G$ be a group and $H$ a subgroup of $G$ such that 
$$[H:(H\cap xHx^{-1})]<\infty$$
for all $x\in G.$ Let $\widetilde{H}(G,H)$ be the free $\mathbb{Z}$-module generated by the double $H\times H$-classes of $G$ (of the type $HxH$ for $x\in G$). The set $\widetilde{H}(G,H)$ can be enriched with  the structure of  an associative algebra (for details, see  \cite[Section 7]{Shimura} or the next paragraph). After Weil proposed this idea, the rain probably stopped and Shimura smiled. Weeks later, while writing this down, Shimura had the horrible idea of not giving a name to this construction. Iwahori \cite{Iwahori} called  this the \textit{Hecke ring} (it should have been called ``Weil ring'' or ``Shimura ring'') because if $G=SL_2(\mathbb{Q})$ and $H=SL_2(\mathbb{Z})$ we get the abstract ring behind Hecke operators in the theory of modular forms.

It was a sunny but cold day of spring in 1964 when Nagayoshi Iwahori  \cite[Section 1]{Iwahori}  gave a simpler description of this ring in general, using measures and integration and that kind of stuff. But for the particular case when $G$ is finite he proved that $\widetilde{H}(G,H)$ is
$^H\mathbb{Z}G^H$, i.e. the set of functions $f:G\rightarrow \mathbb{Z}$ such that $f(g)=f(h_1gh_2)$ for all $g\in G$ and $h_1,h_2\in H,$ with the obvious addition and multiplication  defined by the convolution 
$$(f*f')(x)=\sum_{y\in G}f(y)f'(y^{-1}x).$$

Furthermore, Iwahori  realised that if $B$ is the set of upper triangular matrices, then $\widetilde{H}(SL_n(\mathbb{F}_q),B)$ is the ring generated over $\mathbb{Z}$ by $T_{s_i}$, the characteristic function on $B\dot{s_i}B$ divided by the order of $B$ (we recall that $\dot{s_i}$ is the permutation matrix defined by the transposition $s_i=(i,i+1)$, see Section \ref{Bruhat}), subject to the braid relations on $T_{s_i}$ and to the quadratic relation 
\begin{equation}\label{quadratic}
T^2_{s_i}=(q-1)T_{s_i}+q, \quad \text{where } q \text{ is the order of the field.}   
\end{equation}

(His theorem, as most of the time in this paper, was for more general Chevalley groups).

\subsection{Iwahori-Hecke algebras}\label{Hecke}
Enough of talking. Let $\mathbb{Z}[v,v^{-1}]$ be the Laurent polynomial ring. The \textit{Hecke algebra} $H(W,S)$ of a Coxeter system $(W,S)$ is the $\mathbb{Z}[v,v^{-1}]$-algebra with generators $h_s$ for $s\in S$ and relations 
\begin{itemize}
\item  $h_s^2=(v^{-1}-v)h_s + 1$
\item $\underbrace{h_sh_rh_s\cdots}_{m(s,r)} =\underbrace{h_rh_sh_r\cdots}_{m(s,r)}$ for all $s,r \in S.$
\end{itemize}

You might wonder where the $v$ and the $h_s$ came from? It is just a slight normalization that makes the next section much prettier. We define $v^{-2}=q$ and $h_s:=vT_s$ and obtain Equation \eqref{quadratic}.

To be more precise, a ring morphism $\varphi:\mathbb{Z}[v,v^{-1}] \rightarrow \mathbb{C} $ is determined by an invertible complex number $\varphi(v)$. This map gives $\mathbb{C}$ the structure of a $\mathbb{Z}[v,v^{-1}]$-module and we may form the \textit{specialization} of the Hecke algebra
$$H_{\varphi(v)}:=\mathbb{C}\otimes_{\mathbb{Z}[v,v^{-1}]}H.$$
The specialization of $H$ at $v\mapsto 1$ (i.e. $H_1$) is the group algebra $\mathbb{C}W$, and by Iwahori theorem above, when $q$ is the power of a prime number,
$$H_{q^{-\frac{1}{2}}}\cong \mathbb{C}\otimes_{\mathbb{Z}}\widetilde{H}(SL_n(\mathbb{F}_q),B).$$


\vspace{0.2cm}

\textbf{Note to the reader:} If you invent a concept in mathematics, try to find a really bad name for it. It will then be named after you. For example, Soergel called ``Special bimodules'' the objects that now are called ``Soergel bimodules'' and  Vaughan Jones originally called the Jones polynomial  ``trace invariant''. Never do what Shimura did, to not give a name to your invention, nor what Iwahori did, to call your invention by someone else's name. That name will stick. Don't be humble. You want stuff named after you!

\vspace{0.2cm}


\begin{remark}

A natural question that might arise is the following. How is this ``Hecke algebra''  related to $GL_n(\mathbb{F}_q)$ (the other $q$-deformation of the symmetric group, as explained  in Section \ref{phantom})? They are related in several ways, one of them is  explained in Section \ref{Weil}, but another one that I find particularly fascinating and deep was discovered \cite{Joyal} by André Joyal and  Ross Street  (they invented the concept of ``braided monoidal category'' \cite{Joyal2}) and it  says something like this. Fix $q$ and ``assemble together for all $n$'' the categories  of finite dimensional complex representations of $GL_n(\mathbb{F}_q)$. By assembling  all these categories one obtains a big category (which is braided monoidal, as the representations of quantum groups also are (see Section \ref{qg}))  equivalent to what they call the \textit{Hecke algebroid} that is roughly a category obtained by assembling together all the Hecke algebras $H(S_n)$. This Hecke algebroid can also be described in terms of representations of generalised Hecke algebras in the sense of Robert Howlett and Gustav Lehrer \cite{HL}, so this last interpretation could be summarized as this: ``representations of  $GL_n(\mathbb{F}_q)$ for all $n$ at the same time are equivalent to representations of  generalized Hecke algebras $H'(S_n)$ for all $n$ at the same time''.
\end{remark}




\section{Kazhdan-Lusztig theory}
Hecke algebras are important for several reasons. We have mentioned the Jones polynomials, but let us briefly mention some more. For example, $H(\widetilde{S}_n)$ is closely related to the representation theory of $SL_n(\mathbb{Q}_p)$ (see \cite{IM}) as well as to the modular representation theory of $SL_n(\overline{\mathbb{F}}_p),$ in the vein of $(1)$ in Section \ref{fancy}.  Hecke algebras are closely related to the Temperley-Lieb algebras
which arise in both statistical physics and quantum physics (the related examples were key in the discovery of quantum groups). They appear in papers by Richard Dipper and Gordon James on modular representations of finite groups of Lie type (see \cite{Dipper} for example). But by far the most important reason of why they are important, in my opinion, is  that Kazhdan-Lusztig polynomials  naturally live inside the Hecke algebra, as I will now explain.

\subsection{Definition of Kazhdan-Lusztig polynomials}\label{KLP}
We start with a Coxeter system $(W,S)$.
Let $x=sr\cdots t$ be a reduced expression of $x$ (recall that this means that it can not be written using less simple reflections). A theorem by Hideya Matsumoto \cite{Matsumoto} says that if one defines $h_x=h_sh_r\cdots h_t$, this definition only depends on $x$ and not on the chosen reduced expression. One can prove that $$H(W,S)=\bigoplus_{x\in W}\mathbb{Z}[v,v^{-1}]h_x.$$

The set $\{h_x\}_{x\in W}$ is called the \textit{standard basis}. There are two involutions of $H(W,S)$ (we will just call it $H$) that are impressively deep, but they seem quite innocent at first view. Kazhdan-Lusztig theory emerges from the first one, and  Koszul duality from the second one. 

\begin{itemize}
\item Define $d:H\rightarrow H,$ the $\mathbb{Z}[v,v^{-1}]$-module homomorphism (that can be proved to be a ring homomorphism)  by $d(v)=v^{-1}$ and $d(h_x)=(h_{x^{-1}})^{-1}$. To prove that this last equation makes sense (i.e. that there is an element $(h_y)^{-1}$ in $H$) it is enough to see that $(h_s)^{-1}=h_s+v-v^{-1}\in H$.
\item Define $\iota : H \rightarrow H,$ the $\mathbb{Z}[v,v^{-1}]$-module homomorphism (that can be proved to be a ring homomorphism)
by $\iota(v)=-v^{-1}$ and $\iota(h_x)=h_x$.
\end{itemize}

The main theorem of the revolutionary paper \cite{KL} is that for every element $x\in W$ there is a unique element $b_x\in H$ such that $d(b_x)=b_x$ and such that 
$$b_x\in h_x+\sum_{y\in W}v\mathbb{Z}[v]h_y.$$

In this formulation you see that we need $v$ instead of $q$. The set $\{b_x\}_{x\in W}$ is a $\mathbb{Z}[v,v^{-1}]$-basis of $H$ called the \textit{Kazhdan-Lusztig basis}. If $b_x= \sum_{y\in W}h_{y,x}h_y$ then the \textit{Kazhdan-Lusztig polynomials} are defined by the normalization $p_{y,x}=v^{l(x)-l(y)}h_{y,x}$.

These polynomials are quite ubiquitous in representation theory, they appear in tens (maybe hundreds) of formulas. Usually they appear as the multiplicity of some kind of representation theoretical object in some other kind of representation theoretical object (examples of representation theoretical objects: simple, tilting, projective, injective, standard, costandard, etc.). But it also appears in  geometry: perverse sheaves, Lagrangian subvarieties, Springer resolution, etc. The website  \cite{Vaz05} by Monica Vazirani is very nice as a source of references for Kazhdan-Lusztig theory.

\subsection{Pre-canonical bases}\label{precan}
Probably the most important Kazhdan-Lusztig polynomials appear when $W$ is an affine Weyl group. To fix ideas, as we have done in the rest of the paper, we will suppose $W=\widetilde{S}_n$. Among the most important elements for representation theory are the $y$ (in $p_{y,x}$) of the form $\theta(\lambda)$ for $\lambda$ a dominant weight.

One fundamental appearance of  these polynomials (after evaluating  at $1$) is that they give the formal characters\footnote{I have not explained what a formal character is because I don't want to introduce Lie algebras, but they are analogues of characters for simple groups as explained in \ref{educated} and are probably the most important piece of data one can extract from a representation. } of  the irreducible representations $V_{\lambda}$ (see Remark \ref{domi}) of the  group $SL_n(\mathbb{C})$. 

There are several famous formulas for $p_{\theta(\lambda),x}(1)$. For example, Weyl's character formula by Hermann Weyl \cite{Weyl} is an  extremely useful formula, Hans Freudenthal's formula  \cite{Freudenthal} is good for actual calculations. Bertram Kostant's formula  \cite{Kostant} is just beautiful. By far my favorite formula is that of Peter Littelmann \cite{Littelmann} with level of beautifulness close to infinity (and it is also a sum of positive terms, unlike the previous formulas, that have signs). Finally there is  Siddhartha Sahi's formula  \cite{Sahi}. There is also a formula for the full Kazhdan-Lusztig polynomial $p_{\theta(\lambda),x}$ by Alain Lascoux and Marcel-Paul Sch\"{u}tzenberger,  \cite{Lascoux1} and another one by Alain Lascoux, Bernard Leclerc,  and Jean-Yves Thibon  \cite{Lascoux2} (for a new interpretation of these formulas using geometric Satake, see Leonardo Patimo's paper \cite{Patimo}).

All of these formulas have pros and cons, but they all have the problem that they are not easy to compute. One could reply that they are computing some very complicated objects so the answer has to  be complicated, but I believe that we could produce a more simple understanding of these polynomials.  This was the starting point for Leonardo Patimo, David Plaza and I in the paper \cite{precan}, where for $\widetilde{S}_n$ we divide the problem into $n$ steps. I will explain this approach to do some publicity of this new concept that I absolutely love, and also because it will give you a taste for the complexity of the result in its easiest form, to my knowledge. You don't need to read in detail the formulas, they are just written in detail to show you
that they are very good-looking.

Let me be more precise. Computing $p_{\theta(\lambda),x}$ for all $x$, is equivalent to find $b_{\theta(\lambda)}$ in terms of the standard basis .
Define 
\begin{equation}\label{N}
    N_x=\sum_{y\leq x}v^{l(x)-l(y)}h_y,
\end{equation}
where $\leq$ is the Bruhat order and $l(\cdot)$ is the length function, both defined in Section \ref{Coxeter}.
I will give formulas for $b_{\theta(\lambda)}$ in terms of the $N_x$. For $\widetilde{S}_3$  this formula was probably known  and for $\widetilde{S}_4$ we found it this year with Leonardo Patimo and David Plaza \cite{precan}.
\subsubsection{``Easy'' case $\widetilde{S}_3$}  We denote $\theta(m,n):=\theta(m\omega_1+n\omega_2)$. Here the answer is quite simple, although the proof is not simple at all. 
$$b_{\theta(m,n)}=\sum_{i=0}^{\mathrm{min}(m,n)}v^{2i}N_{\theta(m-i,n-i)}.$$
In the paper \cite{LP} with Leonardo Patimo (following ideas of Geordie Williamson) we not only proved this formula, but furthermore found explicitly the  Kazhdan-Lusztig polynomials for all pairs of elements in this group. It is the only infinite group for which explicit formulas are available. The formula above was probably known to the experts. 
\subsubsection{Fun case $\widetilde{S}_4$}\label{funcase}
Recall the definitions in the beginning of Section \ref{gage}.
Let  $$\lambda = a\omega_1 + b\omega_2 + c\omega_3 \in X^+ $$ 
where $X^+$ is the set of  dominant weights (i.e. with all $\lambda_i\geq 0$). We use the notation $\alpha_{13}:=\alpha_1+\alpha_2+\alpha_3$, $\alpha_{12}:=\alpha_1+\alpha_2$ and $\alpha_{23}:=\alpha_2+\alpha_3$. For each dominant weight $\lambda$ there are elements $N^3_{\lambda}\in H$ and $N^2_{\lambda}\in H$ that are particularly interesting, they are called the ``pre-canonical bases'' (for the definition and some properties see \cite{precan}). We will not need their definition in this section, just  that they are some elements of the Hecke algebra. We will give three formulas that together with Equation \eqref{N} describe $b_{\theta(\lambda)}$ in terms of the standard basis. 

\begin{equation*} \label{eq n4 en n4 intro}
    b_{\theta(\lambda)} = \displaystyle\sum_{k=0}^{\min(a, c)} v^{2k} \mathbf{N}_{\lambda -k\alpha_{13}}^3.
\end{equation*}

For the second formula we need to introduce the set $I_{\lambda}$.
 It is the set of $\mu\in X^+$ such that there exist $n, m, l\in \mathbb{N}$  with  \[\mu=\lambda-n\alpha_{12}-m\alpha_{23} \ \ \ (\mathrm{in\ this\ case\ we\ consider\ }l\ \mathrm{to\ be\ }0)\] or \[\lambda-n\alpha_{12}-m\alpha_{23}\in \mathbb{N}\omega_1+\mathbb{N}\omega_3\ \mathrm{ and}\  \mu=\lambda-n\alpha_{12}-m\alpha_{23}-l\alpha_{13}\in X^+.\]   
For $\mu\in I_{\lambda}$ we define $d(\mu):=n+m+2l$. 
\begin{equation*}\label{eq intro N3 in N2}
   \mathbf{N}_{\lambda}^3 = \displaystyle\sum_{\mu \in I_{\lambda}} v^{2d(\mu)} \mathbf{N}_{\mu}^2.   
\end{equation*}
 
The last formula is 
\begin{equation*}\label{1}
    \mathbf{N}_{\lambda}^2 = \displaystyle \sum_{\substack{\mu \in X^+ \\ \mu \leq \lambda \ \  } }  v^{2\mathrm{ht}(\lambda-\mu)} \mathbf{N}_{\theta(\mu)},
\end{equation*}
where  $\mu\leq \lambda$ means that $\lambda-\mu\in \mathbb{N}\alpha_1+\mathbb{N}\alpha_2+\mathbb{N}\alpha_3$. The notation $\mathrm{ht}$ denotes the \textit{height of a weight}: if $\lambda=\sum_im_i\alpha_i$ then $\mathrm{ht}(\lambda)=\sum_im_i.$

\subsection{Representations of the Hecke algebra}
Kazhdan and Lusztig in their foundational paper \cite{KL} introduced the concept of (left, right, two-sided) cells, certain partitions of the Coxeter group $W$.

For $x,y\in W$ we write  $x\preceq_{\mathrm{L}}y$, resp. $x\preceq_{\mathrm{R}}y$
if there exist $s\in S$ such that $b_x$ appears with nonzero coefficient when $b_sb_y$ (resp. $b_yb_s$) is written in the Kazhdan-Lusztig basis $\{b_w\}.$ We write $x\preceq_{\mathrm{LR}}y$ if either $x\preceq_{\mathrm{L}}y$ or $x\preceq_{\mathrm{R}}y$. These relations are not necessarily transitive. 
 The preorders $\leq_{\mathrm{L}},\   \leq_{\mathrm{R}},\   \leq_{\mathrm{LR}},$ defined as the transitive closures of the relations $\preceq_{\mathrm{L}},\   \preceq_{\mathrm{R}},\   \preceq_{\mathrm{LR}},$ generate
equivalence relations
$\sim_L, \sim_R, \sim_{LR}$ ($x\sim_L y$ means that $x\leq_{\mathrm{L}}y\leq_{\mathrm{L}}x$). Its equivalence classes are  the \textit{left, right} and \textit{two-sided Kazhdan–Lusztig cells of} $W$, respectively.

Given a left cell $C$, consider the left ideals
$$H_{\leq_{\mathrm{L}}C}:=\bigoplus_{w\leq_{\mathrm{L}}C}\mathbb{Z}[v,v^{-1}]b_w, \hspace{1cm} H_{<_{\mathrm{L}}C}:=\bigoplus_{w<_{\mathrm{L}}C}\mathbb{Z}[v,v^{-1}]b_w.$$
of $H$, and define the  \textit{left cell module}
$$H_C:=H_{\leq_{\mathrm{L}}C}/H_{<_{\mathrm{L}}C}.$$
Kazhdan and Lusztig \cite[Theorem 1.4]{KL} prove that when $W$ is the symmetric group $S_n$, these left cell modules give a full set of irreducible representations of $H(S_n)$. Moreover, there is a beautiful bijection (using the Robinson-Schensted correspondence) between left cells and partitions of $n$ (see \cite[Section 22.2]{GBM} for details), and under this bijection one can see that when $v\rightarrow 1$, the cell module goes to the corresponding Specht module. 

\section{Fast \& Furious intro to Quantum groups}\label{qg}
\epigraph{\textit{The rigid cause themselves to be broken; the pliable cause themselves to be bound.}}{Xun Kuang}
\subsection{A rigid situation}
The fact that most interesting Lie groups and Lie algebras are rigid (i.e. can not be deformed) is known since the sixties. For compact Lie groups it was proved in \cite{Pal} by Richard Palais and Thomas Stewart. For semisimple Lie algebras (such as $\mathfrak{sl}_n(\mathbb{C})$), at least for formal deformations, it
was proved in \cite{Gers} by  Murray Gerstenhaber (adding Whitehead's second lemma). There is an even stronger version of this theorem that says that semisimple Lie algebras are ``geometrically rigid''. Roughly speaking, a Lie algebra $L$ is geometrically rigid if every Lie algebra near $L$ is isomorphic to $L$. 

This situation is similar to that of Coxeter groups, but we learnt in Section \ref{IHalg} that if we can not deform a group, we should deform a natural algebra attached to it, one that encodes all of its information. In that case it was the group algebra, that admitted a presentation by generators and relations and so we deform one of the relations and pam! we are done.

We can and we will follow the same approach, but not with the group algebra because it is monstrous, or, to be more precise, it is Brobdingnagian\footnote{I was searching in the   dictionary for a synonym of ``monstrous'' and this beautiful, self-explanatory word appeared.}. Think of  $SL_n(\mathbb{C})$ and its group algebra. It has a basis consisting of every element in  $SL_n(\mathbb{C})$ which is quite a lot. But you could reply that we did the same for infinite Coxeter groups and there was no problem there, so the problem is not to be an infinite group. The problem is to be uncountably infinite, because it is impossible to describe such a group algebra with a finite (or even countable) set of generators. So our approach collapses using the group algebra (there are other, more technical reasons for this not to work). But there are other nice (i.e. finitely generated, admitting a presentation by generators and relations) algebras naturally attached to $SL_n(\mathbb{C})$. One of these is the algebra of ``polynomial functions''.

\subsection{Deforming polynomial functions}
Some parts of this section follow  $III.75$ (Section written by Shahn Majid) of the fantastic book \cite{Comp} where all of Mathematics is explained by amazing mathematicians! Although in my opinion, an even better book  that explains all of Mathematics is the $300$ pages book  \cite{Dieudonne} by Jean Dieudonné, a masterpiece. Other parts of this section are inspired in \cite{Chari} that is a  very influential, very big book. For a nice history of Hopf algebras see \cite{Nicolas1}. For the reader that wants to go deeper into quantum groups, I strongly recommend the incredibly nice paper by one of the main creators of quantum groups Vladimir Drinfel'd \cite{Drinf}, that is also maybe the most cited paper in representation theory  that I have ever seen. For political reasons Drinfel'd was not able to give his talk in international congress of mathematicians in Berkeley, so all we have is the proceedings. (If you do get a hold of the proceedings, readings Manin's Laudatio is a real pleasure!)

\subsubsection{A ten lines introduction to algebraic geometry (including this title)} The starting point of algebraic geometry is this following link between geometry and algebra. Every subset $X\subseteq \mathbb{C}^n$ defined by polynomials gives rise to an algebra $\mathbb{C}[X]$, called its \textit{polynomial functions} which are restrictions to $X$ of polynomial maps $\mathbb{C}^n\rightarrow \mathbb{C}$ (see Hilbert's Nullstellensatz to deepen this relation). This map $X\rightarrow \mathbb{C}[X]$ sending a variety to an algebra is ``functorial'' in the sense that if $X\subseteq \mathbb{C}^n$ and $Y\subseteq \mathbb{C}^m$, a \textit{polynomial map} $\phi:X\rightarrow Y$ (i.e. a restriction of a polynomial function in every coordinate $\mathbb{C}^n\rightarrow \mathbb{C}^m$)  induces an algebra homomorphism $\phi^*:\mathbb{C}[Y]\rightarrow \mathbb{C}[X]$ (note that the order is reversed) defined by the formula $$\phi^*(f)(x)=f(\phi(x)), \quad \mathrm{for}\ x\in X, f\in \mathbb{C}[Y].$$

\subsubsection{Baby example}
Consider the group $SL_2(\mathbb{C})$.We think of this group  as a subvariety of $\mathbb{C}^4,$  as it is the set of matrices $\begin{psmallmatrix}\alpha & \beta \\\gamma  & \delta \end{psmallmatrix}$ such that $\alpha\delta-\beta\gamma=1$.  
Here is an example of polynomial function
$$\begin{psmallmatrix}\alpha & \beta \\\gamma  & \delta \end{psmallmatrix}\mapsto 3\alpha^3\delta-i\alpha \beta^{10}\gamma^2\delta^{57}.$$As we are restricting to $SL_2(\mathbb{C})$,   two such maps are equal if they differ by a multiple of $0= \alpha\delta-\beta\gamma-1.$

In other words, if  $a: SL_2(\mathbb{C})\rightarrow \mathbb{C}$ is the function $$\begin{psmallmatrix}\alpha & \beta \\\gamma  & \delta \end{psmallmatrix}\mapsto \alpha,$$ and $b,c,d$ are defined in the obvious way, the ring of polynomial functions is  
\begin{equation}\label{C}
    \mathbb{C}[SL_2]=\mathbb{C}[a,b,c,d]/\langle ad-bc-1\rangle,
\end{equation}
where we quotient by the ideal generated by $ad-bc-1.$

The problem here is that $\mathbb{C}[SL_2]$ is only an algebra (where you add polynomial functions, multiply polynomial functions and multiply  a complex number by a polynomial function point-wise). Why is this a big problem? Because if $A$ is a $\mathbb{C}$-algebra and $M, N$ are representations of $A$ (i.e. $M, N$ are $A$-modules) there is no natural way\footnote{Beware that the tensor product is over $\mathbb{C}$. If it was over $A$ one could see an $A$-module as a coherent sheaf and take  the standard tensor product of coherent sheaves.} to give $M\otimes_{\mathbb{C}}N$ the structure of an $A$-module \footnote{The naive idea of defining $h\cdot(m\otimes n):=(h\cdot m)\otimes (h\cdot n)$ is not well defined as \newline  $((h+h')\cdot m)\otimes ((h+h')\cdot n)$ is usually different from $(h\cdot m)\otimes (h\cdot n)+(h'\cdot m)\otimes (h'\cdot n)$.} Other things that you can not do: define the trivial module $\mathbb{C}$ or the dual module $M^*.$ We want to deform $SL_2(\mathbb{C})$ but we also want to deform its representation theory, so to see $\mathbb{C}[SL_2]$ just as an algebra will give us a very poor theory. But we have a solution for this. 

The functoriality that we mentioned before has immense applications. For example, there is a map $m:SL_2(\mathbb{C})\times SL_2(\mathbb{C})\rightarrow SL_2(\mathbb{C})$  that defines the multiplication in the group (it is just matrix multiplication). One can prove that there is an isomorphism of algebras $\mathbb{C}[SL_2\times SL_2]\cong \mathbb{C}[SL_2]\otimes \mathbb{C}[SL_2]$, so the map $m$ gives rise to an algebra homomorphism $m^*$ that we denote by $$\Delta: \mathbb{C}[SL_2]\rightarrow \mathbb{C}[SL_2]\otimes \mathbb{C}[SL_2].$$
This map is known as the \textit{coproduct} and is given explicitly in the generators of $\mathbb{C}[SL_2]$ (see Equation \eqref{C}) by the formula
$$ \Delta \begin{pmatrix}
a & b \\
c & d 
\end{pmatrix}= \begin{pmatrix}
a & b \\
c & d 
\end{pmatrix}\otimes \begin{pmatrix}
a & b \\
c & d 
\end{pmatrix}.$$
which means that $\Delta a=a\otimes a+b\otimes c,$ etc. 
The associativity of the group $SL_2(\mathbb{C})$ is the expressed in the formula $m(m\times \mathrm{id})=m( \mathrm{id} \times m)$. One can prove that it is equivalent to the equation $(\Delta\otimes \mathrm{id})\Delta=( \mathrm{id}\otimes \Delta)\Delta$.

The fact that the group has an identity can be expressed by the \textit{unit map} $$u:\{1\}\rightarrow SL_2(\mathbb{C}),$$ defined by $1\mapsto e$. In this language, the equations $eg=g=ge$ for all $g\in SL_2(\mathbb{C})$ translate into $m(u\times \mathrm{id})=\mathrm{id}=m(\mathrm{id}\times u)$. 
This produces a counit map $\epsilon:\mathbb{C}[SL_2]\rightarrow \mathbb{C}$ that satisfies some equivalent axiom. Finally there is an antipode map $S:\mathbb{C}[SL_2]\rightarrow \mathbb{C}[SL_2],$ which  corresponds to the group inversion, and that also satisfies an axiom equivalent to that of the inverse. Adding up, all of this structure makes $\mathbb{C}[SL_2]$ into a Hopf algebra. For more details on this example see \cite[Section 7.1]{Chari}.

\begin{definition}
A \textit{Hopf algebra} over a field $k$ is a quadruple $(H,\Delta, \epsilon, S)$, where
\begin{enumerate}
    \item $H$ is a unital algebra over $k,$
    \item $\Delta:H\rightarrow H\otimes H$ and $\epsilon:H\rightarrow k$ are algebra homomorphisms such that $(\Delta\otimes \mathrm{id})\Delta=( \mathrm{id}\otimes \Delta)\Delta$ and $(\epsilon\otimes \mathrm{id})\Delta=( \mathrm{id}\otimes \epsilon)\Delta=\mathrm{id}$,
    \item $S:H\rightarrow H$ is a $k-$linear map such that $m(\mathrm{id}\otimes S)\Delta=m(S\otimes \mathrm{id})\Delta=1\epsilon$ where $m:H\otimes H\rightarrow H$ is the product operation on $H$.
\end{enumerate}
\end{definition}

Mostly based in South America, there is a big community of mathematicians led by Nicolas Andruskiewitsch trying to classify Hopf algebras \cite{NicolasS}. It is a huge collective effort and it makes me proud that it is happening here, in South America. 

Let us come back to our example. With a representation theory point of view it is natural to see $\mathbb{C}[SL_2]$ as a Hopf algebra: if $H$ is a Hopf algebra over $k$ and $M,N$ are $H$-modules, then $M\otimes_k N$ is also an $H$-module with 
$$h(m\otimes n):=\Delta(h)(m\otimes n)=\sum_ih_im\otimes h'_in,$$
where $\Delta(h)=\sum_ih_i\otimes h'_i.$ The trivial representation is $k$ with the action $$h(z):=\epsilon(h)z,$$ for $h\in H, z\in k$. Finally, if $M$ is an $H$-module and $M^{*}$ the dual $k$-vector space, then $M^{*}$ is an $H$-module with the formula
$$ (hf)(m):=f(S(h)m),$$
for $h\in H, f\in M^{*}, m\in M.$ So we see that $\Delta, \epsilon$ and $S$ are fundamental to define fundamental representation theory objects. 

There is no universal definition of quantum groups, but they are all Hopf algebras. For example, for $q\in \mathbb{C},$ the quantum group $\mathbb{C}_q[SL_2]$ is the free associative (noncommutative!) algebra on symbols $a,b,c$ and $d$ modulo the relations
$$ ba=qab,\ \ \ bc=cb,\ \ \ ca=qac,\ \ \ dc=qcd,$$
$$ db=qbd,\ \ \ da=ad+(q-q^{-1})bc,\ \ \ ad-q^{-1}bc=1. $$
There are suitable maps $\Delta, \epsilon$ and $S$ that give this algebra the structure of a Hopf algebra (in fact, $\Delta$ is given by the same formula as it is for $\mathbb{C}[SL_2]$) and when $q\rightarrow 1$  one obtains $\mathbb{C}[SL_2]$ back -- as a Hopf algebra. This example generalizes to give quantum groups $\mathbb{C}_q(G)$ for all complex simple Lie groups $G.$

\subsection{Drinfeld-Jimbo quantum groups}
\epigraph{\textit{It is important that a quantum group is not a group, nor even a group object
in the category of quantum spaces.}}{Vladimir Drinfeld}

This epigraph reminds me of another famous quote by Beilinson, Bernstein and Deligne when they introduced perverse sheaves: \textit{``Les faisceaux pervers n'etant ni des faisceaux, ni pervers, la terminologie requiert une explication''}.

The kind of quantum groups explained before seem to me as a direct approach towards a deformation of $SL_n$. But there is another kind of quantum groups,  probably more interesting for representation theory and for knot theory. This kind of quantum groups and the ones explained in the last section are ``dually paired'' in a precise sense. These quantum groups go through the following path. 

$$ SL_n(\mathbb{C})\xrightarrow[algebra]{Lie} \mathfrak{sl}_n(\mathbb{C}) \xrightarrow[algebra]{Enveloping}U_{\mathbb{C}}(\mathfrak{sl}_n)\xrightarrow{Quantization} U_q(\mathfrak{sl}_n)$$

All of the objects appearing in this scheme are fundamental objects in Lie theory, and from a representation theory point of view, to study representations of $SL_n(\mathbb{C})$ or representations of its Lie algebra $\mathfrak{sl}_n(\mathbb{C})$ or  representations of its universal enveloping algebra $U_{\mathbb{C}}(\mathfrak{sl}_n)$ are all very similar. I will not enter into each of these concepts and will just jump to the algebra $U_q(\mathfrak{sl}_n),$ which is a Hopf algebra that converges to $U_{\mathbb{C}}(\mathfrak{sl}_n)$ when $q\rightarrow 1.$

I will give a definition of $U_q(\mathfrak{sl}_n)$ for the reader to have a precise idea of ``how it looks like'', the exact relations  not being so relevant. Of course, when $q=1$ these relations are exactly the relations giving $U_{\mathbb{C}}(\mathfrak{sl}_n)$ that were known before quantum groups appeared in this world.

Let $q\in \mathbb{C}$. Define $a_{i,i}=2, a_{i,j}=-1$ if $\vert i-j\vert =1$ and $a_{i,j}=0$ in the rest of the cases. The \textit{quantum group} $U_q(\mathfrak{sl}_n)$ is the unital associative algebra with generators $\{E_i,F_i,K_i,K_i^{-1}: 1\leq i\leq n-1\}$ subject to the following relations. 
     $$K_iK_i^{-1}=K_i^{-1}K_i=1, \quad K_iK_j=K_jK_i,$$
     $$K_iE_jK_i^{-1}=q^{a_{i,j}}E_j,\quad K_iF_jK_i^{-1}=q^{-a_{i,j}}F_j,$$
     $$E_iE_j=E_jE_i\ \ \mathrm{ and }\ \ F_iF_j=F_jF_i \ \  \mathrm{if}\ \ \vert i-j \vert\geq 2,$$
     $$E_iF_j-F_jE_i=\delta_{i,j}\frac{K_i-K_i^{-1}}{q-q^{-1}},$$
     $$E_i^2E_j-(q+q^{-1})E_iE_jE_i+E_jE_i^2=0\ \ \mathrm{ if }\ \ \vert i-j \vert=1,$$
     $$F_i^2F_j-(q+q^{-1})F_iF_jF_i+F_jF_i^2=0 \ \ \mathrm{ if }\ \  \vert i-j \vert=1.$$

\vspace{0.2cm}

We now give $U_q(\mathfrak{sl}_n)$ a  Hopf algebra structure. The coproduct is given by 
$$\Delta(E_i)=E_i\otimes K_i+1\otimes E_i,\ \ \ \Delta(F_i)=F_i\otimes 1+K_i^{-1}\otimes F_i,\ \ \  \Delta(K_i)=K_i\otimes K_i,$$
the counit by $\epsilon(E_i)=\epsilon(F_i)=0,$ and $\epsilon(K_i)=1$, and the antipode is given by 
$$S(E_i)=-E_iK_i^{-1}, \ \ \ S(F_i)=-K_iF_i,\ \ \ S(K_i)=K_i^{-1}. $$

An important point is that a representation of $U_q(\mathfrak{sl}_n)$ is just an $U_q(\mathfrak{sl}_n)$-module (we ignore $\Delta, \epsilon$ and $S$). There is one particular representation $V_q$ of $U_q(\mathfrak{sl}_n)$, called the \textit{fundamental representation} that has a basis $\{v_1,\ldots,v_n\},$ and the action of the generators of $U_q(\mathfrak{sl}_n)$ in this basis is given by 
$$K_i\mapsto qE_{i,i}+q^{-1}E_{i+1,i+1}+\sum_{j\neq i,i+1}E_{j,j}$$
$$E_i\mapsto E_{i,i+1},\ \ \ \ \ \ \ \  F_i\mapsto E_{i+1,i},$$
where $E_{i,j}$ is the $n\times n$-matrix with a $1$ in position $(i,j)$ and zero elsewhere.

Finally we can enunciate ``Quantum Schur-Weyl duality'', a theorem due to Michio Jimbo \cite{Jimbo}. Define ${\check{R}}_{V,V}:V_q^{\otimes 2}\rightarrow V_q^{\otimes 2}$ by the equations
\begin{equation*}
{\check{R}}_{V,V}(v_i\otimes v_j) =
\begin{cases}
 qv_i\otimes v_j & \text{if }  i=j,\\
v_j\otimes v_i & \text{if } i<j,\\
v_j\otimes v_i+(q-q^{-1})v_i\otimes v_j & \text{if } i>j.
\end{cases}
\end{equation*}
Then the map $h_{s_i}\mapsto \mathrm{id}_{V_q^{\otimes(i-1)}}\otimes {\check{R}}_{V,V}\otimes \mathrm{id}_{V_q^{\otimes(n-i-1)}}$ extends to an action of the Hecke algebra $H(S_n)$ on the space $V_q^{\otimes n}$. With this action one has that 
$$\mathrm{End}_{U_q(\mathfrak{sl}_n)}(V_q^{\otimes n})\cong H(S_n), \ \ \mathrm{and}$$
$$\mathrm{End}_{H(S_n)}(V_q^{\otimes n})\cong U_q(\mathfrak{sl}_n).$$
Of course, in the classical Schur-Weyl duality as explained in Section \ref{Duality} one can replace $GL(V)$ by $U_{\mathbb{C}}(\mathfrak{sl}_n)$ and it works as well.

\subsection{Quantum groups at roots of unity}\epigraph{\textit{A mathematician is a person who can find analogies between theorems; a better mathematician is one who can see analogies between proofs and the best mathematician can notice analogies between theories. One can imagine that the ultimate mathematician is one who can see analogies between analogies.}}{Stefan Banach}

It is interesting that the quantum group $U_q(\mathfrak{sl}_n),$ where $q$ is specialized at a $p^{\mathrm{th}}$-root of unity (with $p$ a prime), has many similarities to $SL_n(\overline{\mathbb{F}_p})$. One avatar of this analogy is that if $\epsilon$ is a primitive $p^{\mathrm{th}}$-root of unity there is a ring homomorphism from $U_{\epsilon}(\mathfrak{sl}_n)$ to the enveloping algebra  $U_{\mathbb{F}_p}(\mathfrak{sl}_n)$ of $\mathfrak{sl}_n$ over $\mathbb{F}_p$, and under certain conditions, representations of $U_{\epsilon}(\mathfrak{sl}_n)$ can be specialized to give representations of $U_{\mathbb{F}_p}(\mathfrak{sl}_n)$.  
For a deep relation between these theories (mostly between their representation theories), see \cite{Lu5}  and for the baby example $SL_2$ see \cite{Tubb} where both categories of representations are generalized at once. 

I find interesting the analogy between (the analogy between $S_n$ and $GL_n$) and (the analogy between the quantum group $U_q(\mathfrak{sl}_n)$ at root of unity and the algebraic group $SL_n(\overline{\mathbb{F}_p})$). They both seem to exist in a deeper level in their representation theories. And yes, I am the ultimate mathematician because I just found an analogy between analogies, thank you very much Banach for your kind words.


\part{Categorical level}

\section{More recent objects of study.}
\subsection{Categorifications}\epigraph{\textit{If one proves the equality of two numbers a and b by showing first that $a \leq b$ and then that $a \geq b$, it is unfair; one should instead show that they are really equal by disclosing the inner ground for their equality.}}{{Emmy Noether}}

When saying this, Emmy Noether was foreseeing a whole field of mathematics that would develop almost a century later, that of ``categorifications''. 

A first instance of categorification is that of the natural numbers $\mathbb{N}$ (including zero). One can think of the category $\mathcal{F}$ of finite sets or the category $\mathcal{V}$ of finite dimensional vector spaces. Both categories categorify $\mathbb{N}$, in the sense that, if $\mathcal{F}'$ and $\mathcal{V}'$ are the categories obtained from $\mathcal{F}$ and $\mathcal{V}$
by identifying isomorphism classes of objects, one obtain isomorphisms  $\mathcal{F}'\cong \mathbb{N}\cong \mathcal{V}'$. In other words, in $\mathcal{F}'$ there is one object for each cardinality (because two sets with equal cardinality are isomorphic in $\mathcal{F}$) and in $\mathcal{V}'$ there is one object for each dimension. I remark that decategorification is a well-defined process while categorification is not.

What is the great thing about categorifications? That there are no maps between the number $3$ and the number $5$ but there are tons of maps between a set with $3$ elements and a set with $5$ elements, as well as tons of linear maps between a $3$-dimensional and a $5$-dimensional vector space. But these categorifications are much more impressive than that. In the natural numbers you can add $3+5$ and this corresponds to disjoint union in $\mathcal{F}$ and direct sum in $\mathcal{V}$. In $\mathbb{N}$ you can multiply $3\cdot 5$ and this corresponds to the cartesian product in $\mathcal{F}$ and to the tensor product in $\mathcal{V}$... impressed?

In the mid nineties Louis Crane and Igor Frenkel \cite{CrFr} suggested that there should exist a ``categorified Lie group'' whose representations would give rise to $4$-dimensional quantum field theories. Louis Crane went even further and suggested \cite{Cr} that  this might give the right setting for  Quantum Gravity\footnote{In Crane's paper one can find the very surprising fact that even Albert Einstein was at the end of his career thinking these kind of ideas. In a letter to Paul Langevin Einstein said ``The other possibility leads in my opinion to a renunciation of the space-time
continuum, and to a purely algebraic physics.''}, some kind of ultimate desire for physicists, one that would unite quantum physics and general relativity. The categorification of the Lie groups, or more precisely, of quantum groups, was achieved independently by Raphaël Rouquier and by Mikhail Khovanov and Aaron Lauda. That is the object denoted by $\mathbb{A}(sl_n)$ in Equation \eqref{swd} and is called \textit{2-Kac-Moody algebra}. But we will speak more about the other part of the categorical action, namely about Soergel bimodules $\mathcal{B}_n$.

\subsection{Some historical background}
The first part of the story is the geometric part. Recall from Section \ref{Weil} that the Hecke algebra of $S_n$ (specialized at $q^{-1/2}$) is the algebra $^B\mathbb{Z}SL_n(\mathbb{F}_q)^B$. Grothendieck has a philosophy to categorify stuff: one replaces functions by sheaves. So this algebra was replaced by Soergel with the algebra $\mathrm{Der}^{\mathrm{ss}}_{B\times B}(SL_n(\mathbb{F}_q))$. I will not define this ``$B\times B$-equivariant semisimple derived category'' (see \cite{Lunt}), but just say that it is not as difficult as it appears, as we will see in the next paragraphs, and that it has a convolution defined between its objects, categorifying the convolution defined in Section  \ref{Weil}. 

The second part is a Lie algebra tale. The category of finite dimensional representations of the Lie algebra $\mathfrak{sl}_n$ (traceless $n\times n$ matrices) is similar to that of $SL_n(\mathbb{C})$ and it is well known, irreducible modules are in bijection with the dominant weights $X^+$ and their characters are given by the formulas discussed in Section \ref{precan}. On the other hand, understanding all infinite dimensional representations of $\mathfrak{sl}_n$ is  suicide\footnote{The only Lie algebra for which a classification of the irreducible representations is known is $\mathfrak{sl}_2(\mathbb{C})$, see \cite{Block}.}.  One needs to study something that is still studiable but fun. That is exactly the BGG category $\mathcal{O}$. It is essentially the category that contains the finite dimensional representations of $\mathfrak{sl}_n$ but also some important modules called ``Verma modules'' and their simple quotients.   

Both stories collide in that of Soergel bimodules. The category of Soergel bimodules is equivalent to $\mathrm{Der}^{\mathrm{ss}}_{B\times B}(SL_n(\mathbb{F}_q))$ and tensor product of bimodules corresponds to the convolution mentioned above. The category of Soergel bimodules is extremely similar to the category $\mathcal{O}$, and modelled upon it. One can make this very precise via Soergel's functor $\mathbb{V}$ (see, for example \cite[Section 7]{LL}).

\subsection{Soergel bimodules}\label{sbim}
In this section we study a categorification of the Hecke algebra (we will explain what this means in a minute). 
For $k$ a field, the symmetric group $S_n$ acts on $R=k[x_1,\ldots, x_n]$ by permutation of the variables. This ring is graded in an obvious way, by adding the exponents of a monomial, so $x_2, 2x_3^2$ and $4x_1^7x_2^5$ have degrees $1,2$ and $12$ respectively. But here we will multiply this natural grading by two\footnote{Recall the business of needing a square root of $q$ in the Hecke algebra.}, so for us $x_2, 2x_3^2$ and $4x_1^7x_2^5$ have degrees $2,4$ and $24$. For a graded $R$-bimodule $\oplus_{i\in \mathbb{Z}}M_i$, we define the \textit{shifted bimodule} $M(j)$ by the formula $$M(j)_i:= M_{j+i}.$$

Recall that the set of simple reflections $S$ is composed by the simple transpositions $s_i=(i,i+1)$. For $s\in S$, call $R^s$ the subring of $R$ fixed by $s$. 
For a sequence $s,r,\ldots, t\in S$, one defines a \textit{Bott-Samelson bimodule}
$$BS(s,r,\ldots,t):=R\otimes_{R^s}R\otimes_{R^r}\otimes R\cdots \otimes_{R^t}R(l),$$
where $l$ is the number of elements in the sequence $(s,r,\ldots,t)$. The category of \textit{Bott-Samelson bimodules} $\mathcal{B}_{\mathrm{BS}}(S_n)$ is defined as the category  whose objects are direct sums of shifts of Bott-Samelson bimodules and the morphisms are  $(R,R)$-bimodule morphisms. To obtain the category of \textit{Soergel bimodules} $\mathcal{B}(S_n),$ you need to enlarge this category and consider also direct summands of the objects in $\mathcal{B}_{\mathrm{BS}}(S_n)$.  Let me explain this a bit. 

If $x\in W,$ and $x=sr\cdots t$ is a reduced expression, inside the Bott-Samelson bimodule $BS(s,r,\ldots,t)$ there is a well-defined indecomposable summand $B_x$, i.e. $BS(s,r,\ldots,t)=B_x\oplus M$ and $B_x$ can not be written as the direct sum of two non-trivial bimodules. Small remark: in principle $B_x$ sits inside $BS(s,r,\ldots,t)$ in several different ways but there is one favorite way, called in \cite{LLL} the ``favorite projector''. Another small remark: for $x=s$ a simple reflection, the indecomposable $B_s$ is the same as $BS(s)=R\otimes_{R^s}R(1).$

Then Soergel bimodules are just bimodules of the following sort: $$\sum_{i\in \mathbb{Z}}\bigoplus_{x\in W}(B_x(i))^{\oplus n_i},$$
where both sums have finite support. Just to be clear, $$M^{\oplus n}:=\underbrace{M\oplus M\oplus \cdots \oplus M}_{n\ \mathrm{terms}}.$$ 

The category $\mathcal{B}(S_n)$ categorifies the Hecke algebra $H(S_n)$. This means that  multiplication by $v$, addition and  product in $H(S_n)$ is translated to shifting by $(1)$, direct sum and tensor product over $R$ in $\mathcal{B}(S_n)$. Also, equalities are translated to isomorphisms. Even more, Soergel's conjecture, a massive theorem proved by Ben Elias and Geordie Williamson \cite{EWhodge}, says that when $k$ is the field of real numbers $\mathbb{R}$ then $b_x$ is translated to $B_x$. This paragraph will become clear with the following set of examples when $m(s,r)=3$ and $k=\mathbb{R}.$ 

\begin{center}
 \begin{tabular}{|l|l|}
\hline
$b_sb_s=vb_s+v^{-1}b_s$ & $B_sB_s\cong B_s(1)\oplus B_s(-1)$\\
\hline
 $b_sb_{srs}=vb_{srs}+v^{-1}b_{srs}$  & $B_sB_{srs}\cong B_{srs}(1)\oplus B_{srs}(-1)$ \\ 
\hline
$b_{srs}b_{srs}=(v^3+2v^1+2v^{-1}+v^{-3})b_{srs}$ & $B_{srs}B_{srs}\cong B_{srs}(3)\oplus B_{srs}(1)^{\oplus 2}\oplus$\\ 
& $\qquad \qquad \ \ \ \  \oplus B_{srs}(-1)^{\oplus 2}\oplus B_{srs}(-3) $ \\ 
\hline
$b_sb_rb_s=b_{srs}+b_s$ & $B_sB_rB_s\cong B_{srs}\oplus B_s$ \\ 
\hline
\end{tabular}
 \end{center}
 
 \vspace{0.4cm}
 
That this will happen with all possible formulas written in the Kazhdan-Lusztig basis is implied by the theorem mentioned in the last paragraph \cite{EWhodge}. In particular, that theorem proves that all Kazhdan-Lusztig polynomials for all Coxeter systems have positive coefficients, because by results of Soergel \cite{Soergel1}, \cite{Soergel2}  this implies that the coefficients of the Kazhdan-Lusztig polynomials are graded dimensions of some Hom spaces between Soergel bimodules (and that is why categorification is cool!). The positivity of Kazhdan-Lusztig polynomials was a very important conjecture that first appeared in the  foundational paper \cite{KL}.

Everything that I have said in this section could have been done with any Coxeter system over the real numbers. One starts with a Coxeter system $(W,S)$ and defines the \textit{geometric representation} as the vector space $V:=\oplus_{t\in S}\mathbb{R}e_t$, with an action of $W$  given  by the formulas 
$$s(e_t)=e_t +2\mathrm{cos}\bigg(\frac{\pi}{m(s,t)}\bigg) e_s.$$
Then one defines $R$ to be the symmetric algebra on $V^*$ or the polynomial functions on $V$.
In \cite{Libeq} I prove that with the geometric representation  the theory of Soergel bimodules works as well. On the other hand, not so much is known when the field we are working in has positive characteristic. In the Appendix of \cite{LLL} I prove that for Weyl groups and $\mathrm{char}(k)\neq 2,3$ all the theory works just as well. It is not known whether the theory works for other groups in positive characteristic, in particular, it is not known for the most important case, $\widetilde{S}_n$. 

There is an interesting variation \cite{Abe} of the category of Soergel bimodules, that in few words one could describe as ``Soergel bimodules+a specified localization''. This category (provided that a few  technical conditions coming from dihedral groups are satisfied \cite{Abe3}) is better suited for positive characteristic (although Soergel calculus, not needing that technical condition  seems like the winner in the competition of Hecke categories in positive characteristic). For example, for Weyl groups it works in all characteristics.

\subsection{Rouquier complexes}\label{RC}
There is a very nice and deep relation between  the Artin braid group and the category of Soergel bimodules. For $s\in S,$ define the complex of Soergel bimodules
$$F_{s}:=0\rightarrow R\otimes_{R^s}R(1)\rightarrow R\rightarrow 0,$$
where $R\otimes_{R^s}R(1)$ lives in cohomological degree $0$ and the only non trivial map in the complex is the multiplication map $a\otimes b\mapsto ab$. We see this complex as an element of the bounded homotopy category $K^b(\mathcal{B}(W))$.

Raphaël Rouquier proved that the map $\sigma_s\mapsto F_s$ extends to a group morphism from the Artin braid group $B(W)$ to the group of isomorphism classes of invertible objects of $K^b(\mathcal{B}(W))$. I never fully understood his proof  \cite[Theorem 9.7]{Rouquier}, so I provided another, more concrete proof in \cite[Section 2.5]{thesis}. See Sub-Question \ref{inj} below.

 The objects $F_{\sigma}$
for $\sigma \in B(W)$ are called \textit{Rouquier complexes}. It has been proved by Mikhail Khovanov \cite{Khovanov} that using them one can categorify the Jones polynomial, and even more, the HOMFLYPT polynomial, a two-variable generalization of the Jones polynomial. He does an algebraic procedure (taking Hochschild homology) that is analogous to the topological process of closing a braid (see Section \ref{braidgroups}), and thus produces a triply graded vector  space, such that using its graded dimension one can recover the HOMFLYPT polynomial. For details, see \cite{Khovanov}.

\subsection{Soergel calculus}
Let $(W,S)$ be any Coxeter system and $\mathcal{B}(W)_{\mathrm{BS}}$ its Bott-Samelson category over the real numbers $\mathbb{R}$. Raphaël Rouquier told me in 2007 the genius idea of giving a presentation of this category by generators and relations. Let me spell out what this means. 

If $M$ and $N$ are Soergel bimodules we will write $MN$ to mean $M\otimes_RN$. Let $a,b,c,d,e\in S$ and suppose that you have a map $g:B_bB_c\rightarrow B_s$, then we have a map $$\mathrm{id}\otimes g\otimes \mathrm{id}^2:B_aB_bB_cB_dB_e \rightarrow B_aB_sB_dB_e.$$

To be a \textit{generating set of morphisms} $\mathcal{F}$ means that every morphism between two Bott-Samelson bimodules $BS_1$ and $BS_2$ is an $R$-linear combination of maps of the form 
\begin{equation}\label{form}
 (\mathrm{id}^{n_1}\otimes f_1 \otimes \mathrm{id}^{m_1})\circ \cdots \circ (\mathrm{id}^{n_l}\otimes f_l \otimes \mathrm{id}^{m_l}),   
\end{equation}
 
\noindent with $n_i,m_i\in \mathbb{Z}_{\geq 0}$ and $f_i\in \mathcal{F}.$ Let us say that $(n_1,f_1,m_1,\ldots n_l,f_l,m_l)$ is an \textit{expression} of the corresponding morphism in \eqref{form}.

In the search for a generating set of morphisms, I found   a basis for $\mathrm{Hom}(BS_1,BS_2)$ as a left $R$-module, that I called \textit{light leaves} basis in \cite{LL} and \textit{double leaves basis} in \cite{LLL} (these two bases are slightly different, and now the common usage for  a light leaf is one half of a double leaf). This is a combinatorial basis  that allows one to really work\footnote{While writing this paper, a new preprint \cite{KST} was posted on the arXiv with the quite surprising comment that light leaves might be relevant in cryptography.} with and prove stuff (see \cite{Gentle} for a nice explanation). In particular, it is a ``cellular basis'' (some kind of very well behaved basis), which is key for the proof of the ``monotonicity conjecture'' of Kazhdan-Lusztig polynomials \cite{monotonicity} by David Plaza. It was also key for the proof of most recent results in the theory of Soergel bimodules.

Anyways, Rouquier's idea was to present this category by generators and relations. 
This means that one can pass from any expression of a morphism between Bott-Samelson bimodules to any other expression, using only local relations in the generators (for example $\mathrm{id}\otimes f_1=f_2\circ f_3$ with $f_1,f_2\in F$).
 I thought for a while that I have found these generators and relations but there was a flaw in my proof, so I had to restrict  to some particular Coxeter systems \cite{Libgr}. My strategy was to start from a general expression of a morphism and reduce it to an $R$-linear combination of light leaves. 

But some months later, new magic came to this theory. Mikhail Khovanov \cite{KhovanovElias} and his student at that point, Ben Elias, had the following genius idea: to draw morphisms (they restricted to the case $W=S_n$).
So, the generators that I found were, for every $s,r\in S$ some maps
$$m_s:B_s\rightarrow R\quad , \quad j_s:B_sB_s\rightarrow B_s \quad , \quad f_{sr}:\underbrace{B_sB_rB_s}_{m(s,r)}\rightarrow \underbrace{B_rB_sB_r}_{m(s,r)}$$
$$\hspace{-4.7cm} m'_s:R\rightarrow B_s\quad , \quad j'_s:B_s\rightarrow B_sB_s. $$
There are explicit formulas for these maps (for example $m_s(a\otimes b)=ab$) but a short and elegant way of describing them is to say that, modulo multiplying by a non-zero scalar, they are the only degree $1$ (resp. $-1, 1, -1, 0$) maps having the prescribed source and target.  Anyways, if you are not an elegant person, and you need the formulas, see \cite{KhovanovElias}. Elias and Khovanov draw these generating maps (I recall that this is for $W=S_n$, for more general groups one should know how to draw $f_{s,r}$ for a general $m(s,r)$) in the following way.

\begin{figure}[H]
    \centering
    \includegraphics[width=12cm]{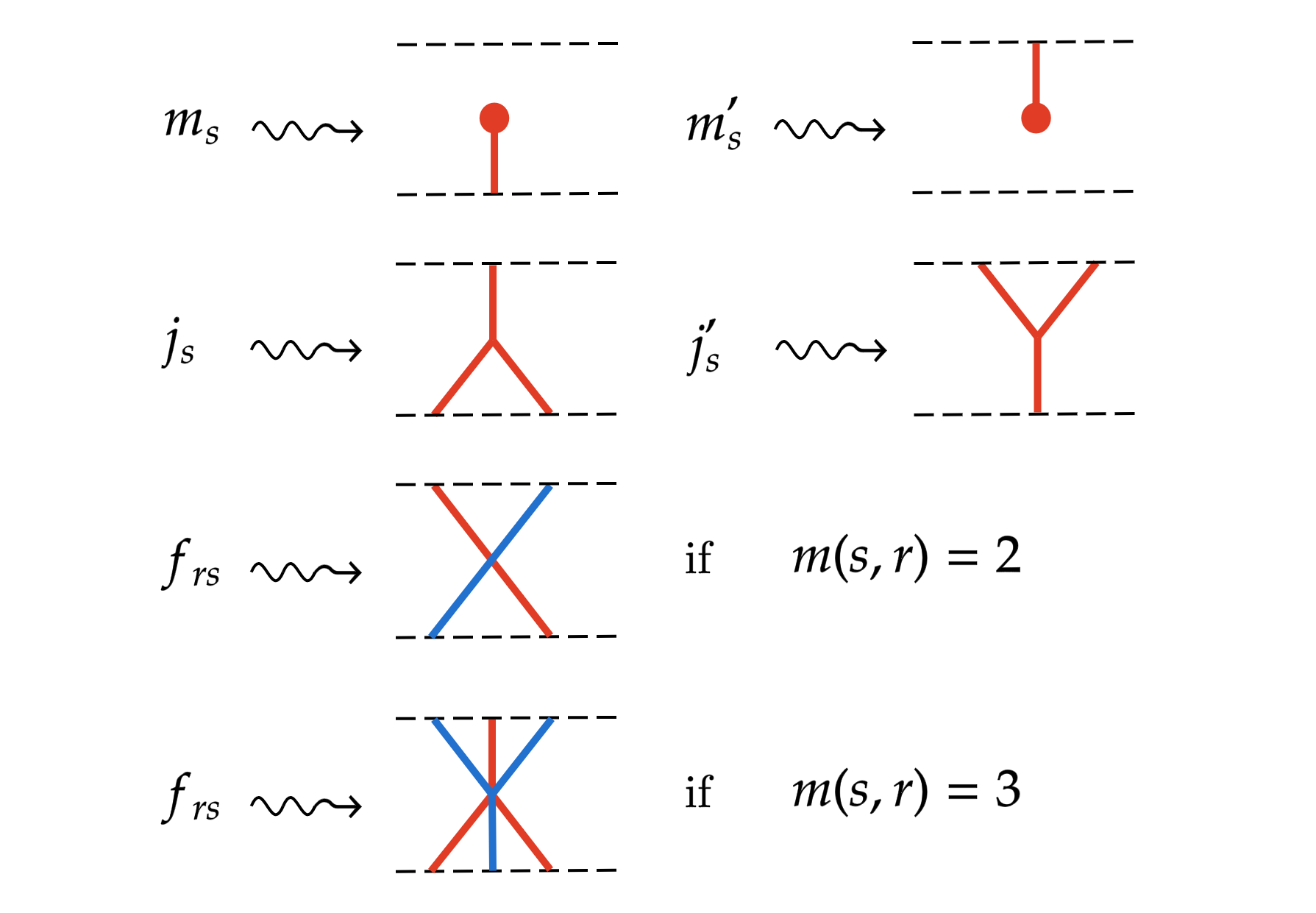}
    \caption{Generating maps}
    \label{gm}
\end{figure}
We draw $s$ red and $r$ blue. The intersection of the bottom dashed line with the colored diagram represents the source, and the intersection of the upper dashed line with the diagram represents the target. For example, in $f_{sr}$ the source is the sequence $(s,r,s)$ and the target $(r,s,r)$ (of course we identify these sequences of simple reflections with the corresponding Bott-Samelson bimodule). The empty sequence is identified with $R.$ Maps go ``from bottom to top''.
Tensor product is drawn putting the diagrams side-by-side. For example $$j_t\otimes f_{sr}:B_tB_tB_sB_rB_s\rightarrow B_tB_rB_sB_r$$
is drawn 
\begin{figure}[H]
    \centering
    \includegraphics[width=5cm]{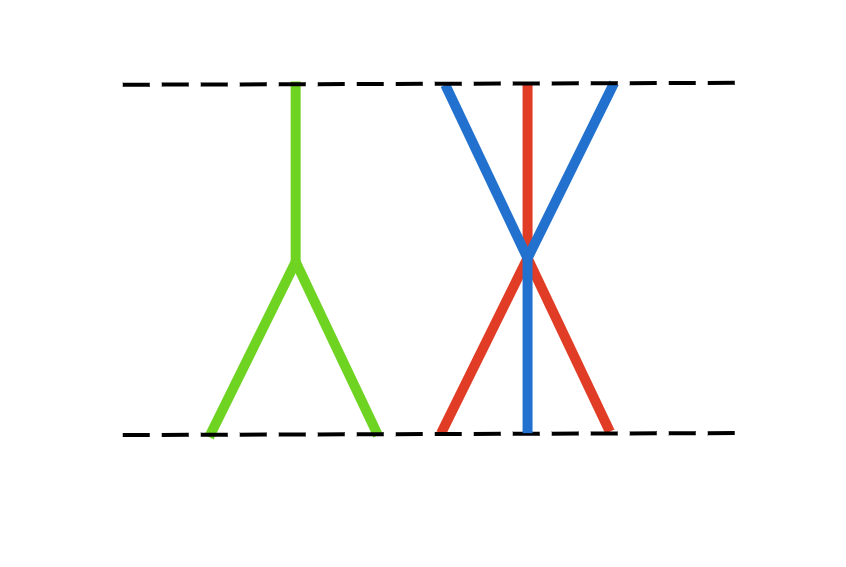}
\end{figure}
\noindent where $t$ is green. Composition is drawn by glueing bottom-up  in the obvious way. For example $j_s\circ (\mathrm{id}\otimes j_s)$ is represented by

\begin{figure}[H]
    \centering
    \includegraphics[width=5cm]{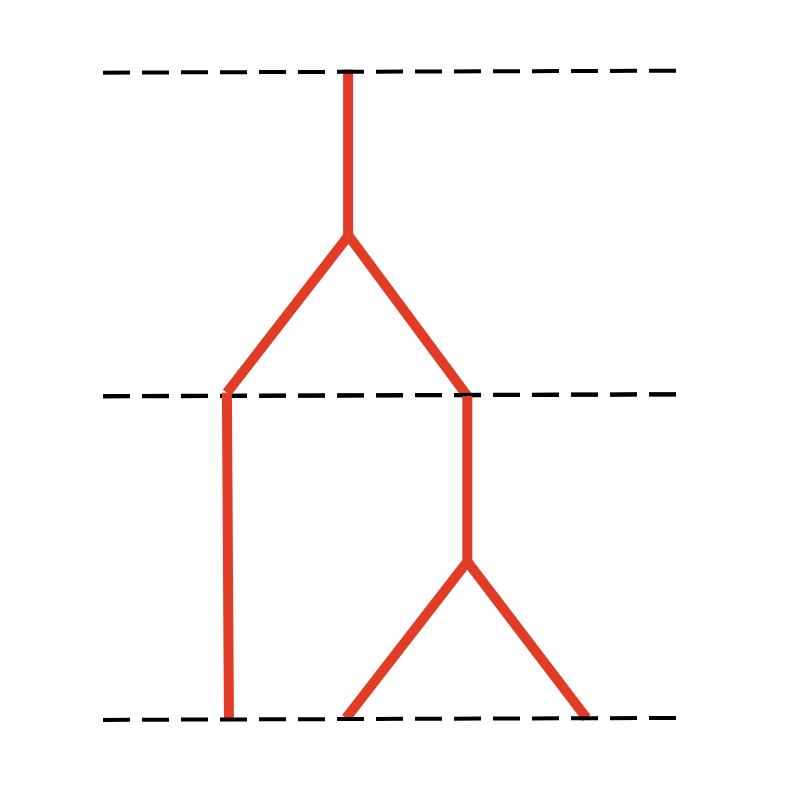}
\end{figure}
\noindent and finally, we use the following notation 
\begin{figure}[H]
    \centering
    \includegraphics[width=7cm]{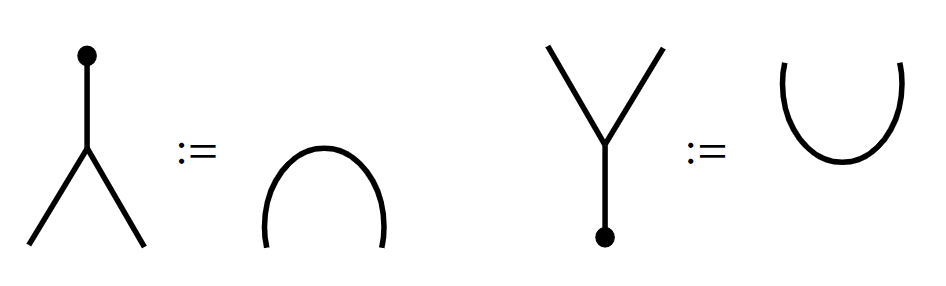}
\end{figure}
But why did I say ``magic'' before? It is because if one draws the maps so, the drawings are isotopy invariants:

\begin{figure}[H]
    \centering
    \includegraphics[width=6cm]{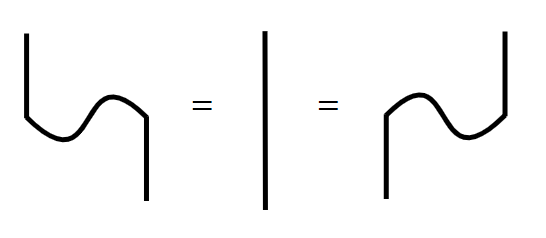}
\end{figure}
\begin{figure}[H]
    \centering
    \includegraphics[width=6cm]{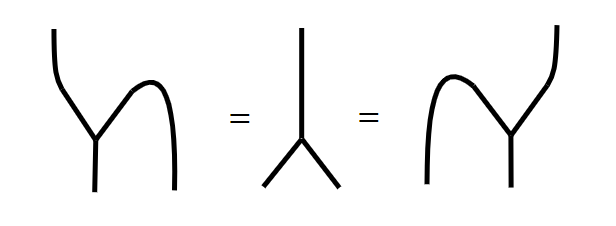}
\end{figure}
\begin{figure}[H]
    \centering
    \includegraphics[width=6cm]{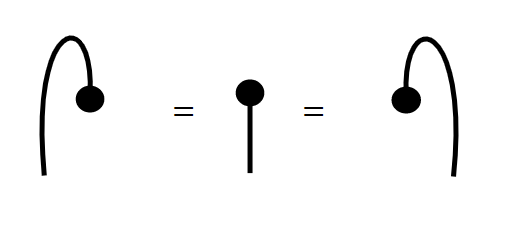}
    \caption{Isotopies}
    \label{isot}
\end{figure}

So Ben Elias and Mikhail Khovanov gave (always in the  $S_n$-case)  with the generators in Figure \ref{gm} a precise set of relations (some twenty of them). We have seen some (one color relations) in Figure \ref{isot}. Here we can see some more (two-color) relations where the two colors $s,r$ involved are such that $m(s,r)=3$: 

\begin{figure}[H]
    \centering
    \includegraphics[width=6cm]{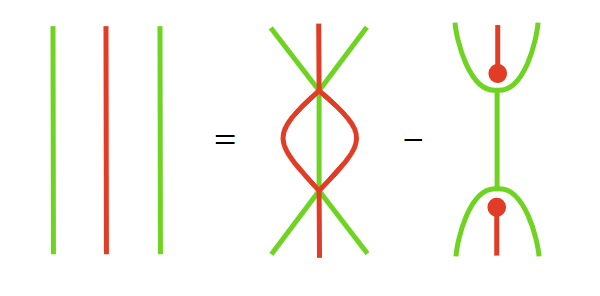}
\end{figure}
\begin{figure}[H]
    \centering
    \includegraphics[width=6cm]{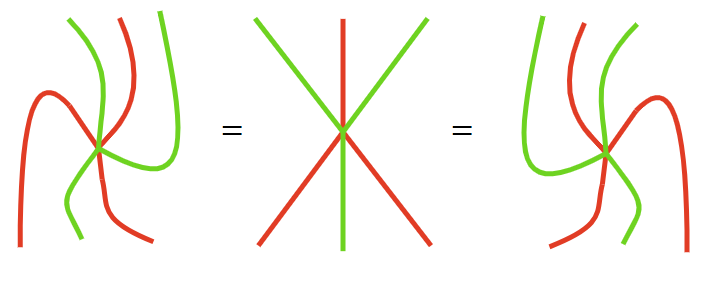}
\end{figure}
\begin{figure}[H]
    \centering
    \includegraphics[width=6cm]{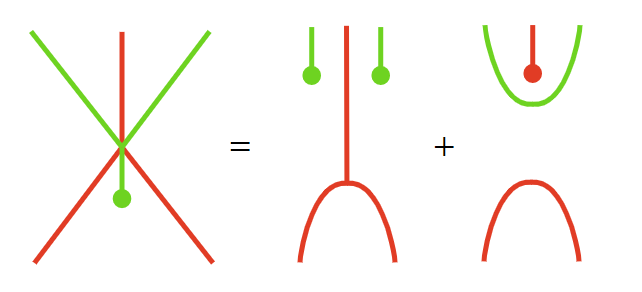}
\end{figure}

This   was later generalized by Ben Elias for every dihedral group \cite{Bendihedral} and finally it was done in \cite{Soergelcalculus} for every Coxeter group. This ``presenting Soergel bimodules by generators and relations'' or simply put \textit{Soergel calculus}, is the single most important tool in the explosion of proofs and counter-examples to conjectures that appeared in the last ten years (for a non-exhaustive list of them see the second page of \cite{Gentle}). We will call this new category presented by generators and relations  the ``diagrammatic Hecke category'' and denote it by $\mathcal{H}(W)$ or $\mathcal{H}^k(W)$ if we want to specify the field $k$ of definition.

Soergel calculus or the diagrammatic Hecke category has three fundamentally new features. First, it allows you to give topological arguments to prove algebraic statements that would otherwise would be too hard to prove algebraically. Second, one can see that the relations defining the Hecke algebra involve one simple reflection (the quadratic relation) or two simple relations (the braid relations) and the new feature here is that $\mathcal{H}(W)$  has relations involving $1, 2,$ and $3$ simple reflections\footnote{One could think that going up one step in the categorical ladder gives rise naturally to one more rank involved, but this is not always the case as for ``singular Soergel bimodules'', the relations have no bound in the number of simple reflections involved.}. The last new feature, and probably the most important, is that this construction ``works'' over all fields and even over more general rings. With this I mean that the diagrammatic Hecke category  over a general field $k$ still categorifies the Hecke algebra  and still has indecomposable objects  $\{B_w\}_{w\in W}$, modulo grading shift. This categorification comes with a map $\mathrm{ch}:\mathcal{H}^k(W)\rightarrow H(W)$ that behaves well\footnote{Technically, to ``behave well'' means that the map ch  descends to an isomorphism of $\mathbb{Z}[v,v^{-1}]$-algebras $\mathrm{ch}:[\mathcal{H}^k(W)]\cong H(W)$, where $[\cdot]$ means the split Grothendieck group of an additive category.} as explained in Section \ref{sbim}.

Usually the category $\mathcal{H}(W)$ depends heavily on the representation of $W$ chosen to define $R$ (see Section \ref{sbim}). If $W$ is the symmetric group $S_n$, or the affine symmetric group $\widetilde{S}_n,$ one does not have to worry about that because with any such representation (satisfying some minimal things), there will be, for $x\in W,$ a well defined element of the Hecke algebra $\mathrm{ch}(B_x)= {}^pb_x,$
called the $p$-\textit{canonical basis}, that only depends on $x$ and on the characteristic $p$ of the field $k.$ When I was doing my Ph.D. I understood that ``Soergel conjecture in positive characteristic'' (although that was never actually conjectured by Soergel) said that ${}^pb_x=b_x$, i.e. that the $p$-canonical basis is the Kazhdan-Lusztig basis for the symmetric group $S_n$.

\subsection{ The $p$-canonical basis. Fractals? Billiards?}\epigraph{\textit{This problem is not made for humans.}}{Wolfgang Soergel}
\label{pcan}
The year was $2008$. I had spent the whole year of $2006$ trying to understand three and a half pages of \cite{Soergel1}, the proof of the ``categorification theorem'' that I have explained before. Those three and a half pages were so incredibly dense that I wrote a  dense text of 30 pages explaining them. It was the first chapter of my Ph.D. thesis. In my mind Wolfgang (Amadeus) Soergel was an extremely serious giant. 

There was  a conference called something like ``Categorifications in mathematics and philosophy'' in Sweden, where real philosophers interacted with us, imagine that. I knew that Soergel was in the conference and was extremely anxious to meet him. Raphaël Rouquier, my Ph.D. advisor, asked me to cross with him a huge garden and meet Soergel. It was raining. We ran and finally got to a little house where some of the talks were given. 
And there he was, a barefooted man, with a huge smile and a flowery shirt. I was trembling. Rouquier said, ``Hi, Wolfgang, this is Nicolas and he wants to prove your conjecture'' (he really said that, while introducing me, can you believe it?). I will never forget what happened next. The answer that Wolfgang Soergel gave to this poor soul was the most impressive and yet the most simple answer imaginable: he burst into an unstoppable laughter. In the middle of that surreal scene, with tears of joy in his eyes, he shouted: ``that is impossible, my conjecture is false'' and then continued to laugh. There was my giant,  my barefooted idol,  mocking me for trying to solve his own \textit{false} conjecture.  Mocking me for having ideals. Mocking me for believing in him.

Anyways, there were several things that I didn't know at that time. First of all, that Wolfgang Soergel is a really funny person. But mathematically, there were already several hints that Soergel's conjecture in positive characteristic, as stated before, was false. What no one knew at that point was to what extent it was 
false. It was really really false. 

Why was Soergel's conjecture so important? Because, in characteristic zero it implied the famous Kazhdan-Lusztig conjectures, some fundamental conjectures concerning multiplicities of simple modules in standard modules in the BGG category $\mathcal{O}$ and it also implied the positivity of the coefficients of Kazhdan-Lusztig polynomials (a massive conjecture in algebraic combinatorics). This part of the conjecture was later proved using Hodge theoretic methods, in \cite{EWhodge}. But even more importantly, Soergel's conjecture for  $S_n$ in characteristic $p>n$ was equivalent to a part of Lusztig conjecture. This last conjecture was the Holy Grail of representation theory of (Lie type) objects. It predicted the characters of simple $SL_n(\overline{\mathbb{F}}_p)$-modules in terms of Kazhdan-Lusztig polynomials for $\widetilde{S}_n$ (for more on the history of this conjecture, in particular the huge odyssey to prove it for $p\gg n,$ see the introduction to \cite{LLL}). 

This reminds me of the quests of Einstein to prove that quantum mechanics was wrong, that was fundamental to prove it right. In the quest to prove Lusztig's conjecture (that Soergel thought was correct due to its inherent beauty) he paved the way to find counter-examples. In \cite{Soergel00} Soergel introduced a subquotient of the category of rational representations of $SL_n(\overline{\mathbb{F}}_p)$,  the ``subquotient around the Steinberg weight” and he realised that it was controlled by $S_n$ and behaved like a modular version of category $\mathcal{O}.$ So it was a key step in this quest that Soergel essentially reduced the problem from Soergel bimodules over $\widetilde{S}_n$ to Soergel bimodules over $S_n.$

 Soergel's laughter was due to the fact that in 2002 Tom Braden (Appendix to \cite{Braden}) had found examples in which ${}^pb_x\neq b_x$ with $p=2$ and $x\in S_8$. Later, in 2011, Patrick Polo (unpublished) proved something even more impressive, that there are $x\in S_{4n}$ such that ${}^pb_x\neq b_x$ if $p$ divides $n-1.$ This was key, since it showed that the $p$-canonical basis is far from equal to the Kazhdan-Lusztig basis. Before that point, one could hope that ``they are equal except for very small primes'', but after Polo's result it was impossible to continue saying so. 

In the quest to prove it or to find counter-examples, using Soergel's approach, in \cite{LLL} I proved that  for $S_n$, if one finds a specific light leaf that satisfies certain properties, then this gives an element $x\in S_n$ such that  ${}^pb_x\neq b_x$ for some $p$ that one can explicitly compute. Using this approach, that he found independently (\cite[Section 4]{Geordiesc}), Geordie Williamson found in \cite{Geordiesc} a set of light leaves satisfying the ``certain properties'' mentioned before, and proved that if $f(n)$ is the lowest prime number $p$ such that there exist  $x\in S_n$ with ${}^pb_x\neq b_x$, then $f(n)$ grows at least exponentially in $n$ (this is what I meant when I said that Soergel's conjecture in positive characteristic was really really false). This also implied that the linear bounds in Lusztig's original paper were too optimistic and it also implied that James' conjecture,  concerning the simple
representations of the symmetric group in characteristic $p$, was  false. 

This result by Williamson was an earthquake in the field. At this point Soergel said what is written in the epigraph of this section (always with a charming smile of acceptance in his face). But slowly light started to emerge from the bottom of this world of despair and a new world started to appear. 

Simon Riche and Geordie Williamson produced ``a new version of Lusztig's conjecture'' in \cite{GeordieSimon} where Kazhdan-Lusztig polynomials are replaced by $p$-Kazhdan-Lusztig polynomials, and using the theory of 2-Kac-Moody algebras, or $\mathbb{A}(sl_n)$ (and we see how all the concepts seen so far start to glue together with Schur-Weyl duality plaining over all of this) and they proved it for $GL_n$. In the paper \cite{AMRW} the authors  prove the general form of this new version of Lusztig's conjecture, i.e. a character formula for irreducible modules for a connected reductive algebraic group in characteristic $p$. They do so by categorifying the $\iota$ involution explained in Section \ref{KLP}, the so called Koszul duality. For nice introductions to this duality see \cite{Maki} and \cite[Section 26]{GBM}.

These and other related papers show that now the ultimate quest is to understand the $p$-canonical basis (or, equivalently the $p$-Kazhdan-Lusztig polynomials). The only case where this is known is for $\widetilde{S}_2$, the infinite dihedral group. In this case one has that 
 the Kazhdan-Lusztig basis is given by the simple formula
$$b_x=\sum_{y\leq x}v^{l(x)-l(y)}h_y. $$
Define $b_n:=b_{\underbrace{srs\cdots}_{n+1}}.$
 We will write the decomposition of the $3$-canonical basis in terms of the canonical basis in the first few examples.
\begin{equation*} 
\begin{split}
{}^3b_0 & = b_0 \\
{}^3b_1 & = \quad \ \  b_1 \\
{}^3b_2 & =\quad \quad \ \ \ \  b_2\\
{}^3b_3& = \quad \ \ b_1\  \  +\quad   b_3\\
{}^3b_4& = b_0\quad \quad \ + \quad \quad \  b_4\\
\end{split}
\end{equation*}
We represent this in the next figure (borrowed from the last page of the paper \cite{JW3}) where each black box represents a one.

\begin{figure}[H]
    \centering
    \includegraphics[width=7cm]{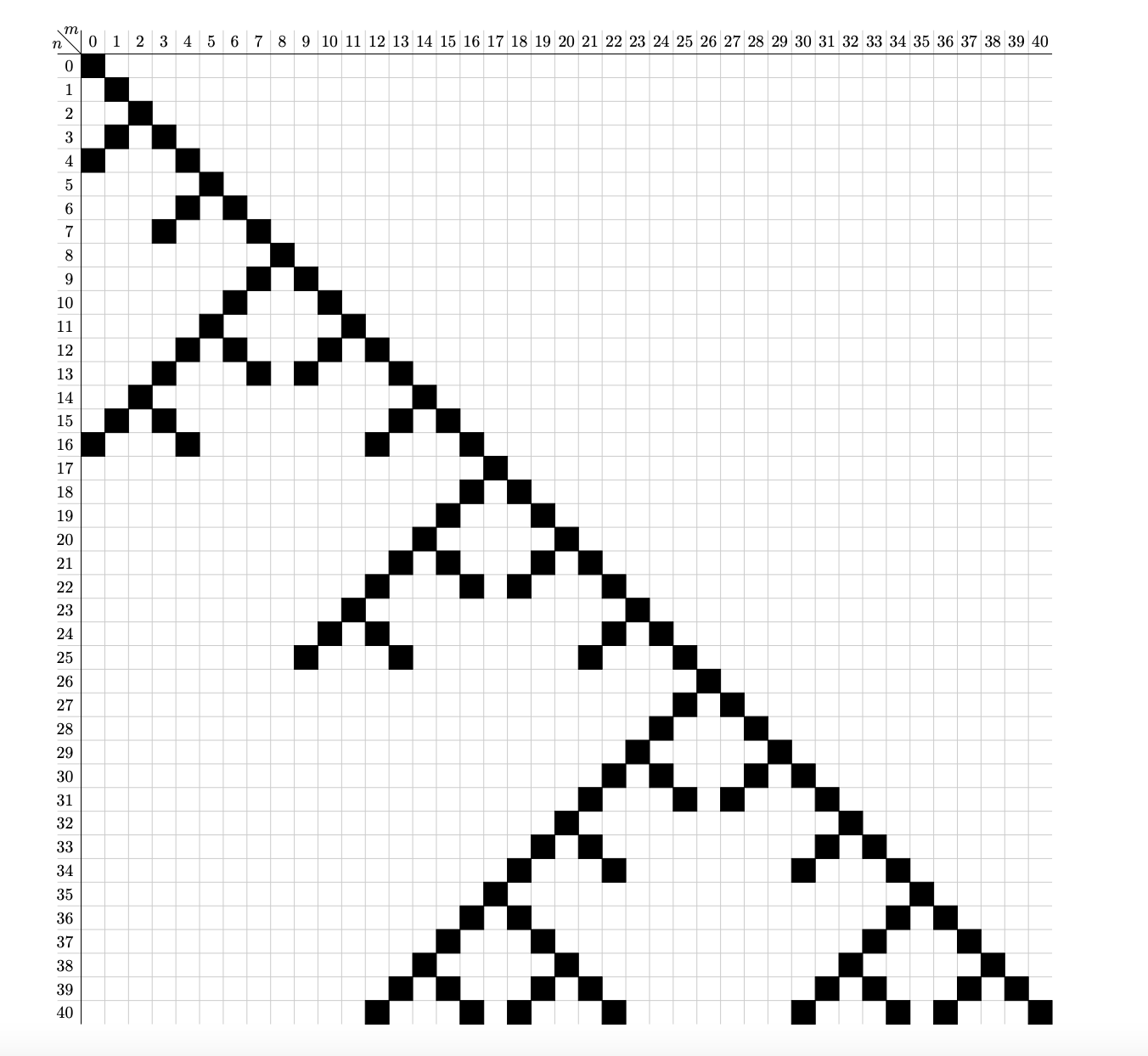}
    \caption{Fractal}
    \label{Frac}
\end{figure}
The reader might see that this is a fractal-like figure. It is a fun exercise to try to guess how this fractal continues, or in other words, to understand the precise sense in which this is a fractal.  Hint: look at the upper-left $(3^n-2)\times (3^{n}-2)$ square, and from there construct the  upper-left $(3^{n+1}-2)\times (3^{n+1}-2)$ square.

As I said before, this is the only case where the $p$-canonical basis is known. But there is a beautiful conjecture by Geordie Williamson and George Lusztig \cite{LusWil} where they predict the decomposition of ${}^pb_x$ in terms of the Kazhdan-Lusztig basis $\{b_y\}$ in a finite part of $\widetilde{S}_3$.  It is not easy to explain the conjecture. I think that the easiest way to understand it is to see the YouTube video by Geordie Williamson called ``Billiards and tilting characters for ${\rm SL}_3$''. In that video there is a black dot that goes from the bottom-left part of the image to the upper-right. That black dot is an element $x_m=s_0s_1s_2s_0s_1s_2s_0\cdots$ (with $m$ terms) (or in the language of Section \ref{why} it is an element of the form $grbgrbgrb\cdots $), where $s_0,s_1$ and $s_2$ are the simple reflections of $\widetilde{S}_3$. The other dots that appear (green, blue and purple) mean the $b_y$'s that appear in the decomposition of ${}^7b_{x_m}$ and the color reflects the multiplicity of $b_y$ in ${}^7b_{x_m}$ (darker color means higher multiplicity). It is a crazy but incredibly beautiful conjecture. I think that the $\widetilde{S}_2$-case could also be seen as a dynamical system if, in Figure \ref{Frac} you started with an horizontal line on top, take that horizontal line down at constant speed and only see what are the black dots that intersect your horizontal line. I would like to think that this is the same phenomenon we see in $\widetilde{S}_3,$ that the billiards or dynamical system behaviour is just a shadow of a fractal behaviour. 

\section{Twenty Love Poems and a Song of Despair}

In this section I will write twenty-three  open questions that I see as love poems, all  related to what I have explained in this paper. I end this section with a song of despair. When I have some vague idea about how to prove these results, I will say it. If I am actively working on it, I will put a sign \danger $\ $ in the beginning of the question.

The problems are ordered from more concrete to more abstract (not from easier to harder). I will also add three numbers for each question. How important I think the problem is, how well-known and how difficult I imagine the problem to be. These three numbers are entirely subjective and essentially  useless, but I do it anyways ``pour la beauté du jeu''. Probably everyone would disagree with me and I am probably very very wrong, so please just take these numbers for what they are: a  very rough approximation.  For example $(I3,WK10,D5)$ means that a question has a $3$ in importance level, it is a $10$ in ``well-knownedness'' and $5$ in difficulty.  The scale is from 1 to $1000$. For example, the Riemann hypothesis would be $(I900,WK1000,D600)$. A paper $(I1,WK1,D1)$  should be published in a top specialist journal like ``Journal of Algebra'' for example. I will just ``throw'' the questions out there and then explain them and maybe give some idea. 

\begin{question} {\color{blue}$(I30,WK1,D30)$}
 For affine Weyl groups, classify all Bruhat intervals modulo poset isomorphism. 
\end{question}
A \textit{Bruhat interval} in an affine Weyl group $W$ is a set  $$[y,z]:=\{x\in W : y\leq x\leq z\}$$ for some $y,z\in W$. The interval inherits a partial order from $\leq$.  This problem seems more like a program, it doesn't seem like something one can do in a single paper. 

A vague idea for the solution is related to Question \ref{egag}. I would like to think that ``most'' isomorphisms of Bruhat intervals come from euclidean isometries (between the connected components) when you consider  Bruhat intervals as subsets of $\mathbb{R}^n$.

\begin{question} {\color{blue}$(I20,WK1,D20)$}
Produce a formula for the number of elements of  every Bruhat  interval in an affine Weyl group. 
\end{question}

There is a formula for the  intervals of this form $[e,\theta(\lambda)]$ for $\lambda$ a dominant weight in \cite{Schutzer} (in a representation theory language), and it is not hard to deduce a formula $[e,x]$ for any $x\in \widetilde{S}_3$ from \cite[Lemma 1.4]{LP}, using the same techniques as in that lemma. 

\vspace{0.2cm}

\hspace{-0.4cm}\textrm{\textbf{Sub-Question 2.}} \danger  $\ $ {\color{blue} $(I6,WK1,D4)$} \textit{Find a formula in the particular case  of intervals of the form $[\theta(\lambda), \theta(\mu)]$ for $\lambda$ a dominant weight. }

\vspace{0.2cm}
Idea: try to find a geometric formula in the sense of Section \ref{gage} for $e,\theta(\lambda)$, we hope (with Damian de la Fuente and David Plaza) that this formulation is more easy to generalize to lower elements different from the identity. 

\begin{question}\label{inj} {\color{blue}$(I50,WK25,D30)$}
For any Coxeter system $(W,S)$, the map $\sigma_s\mapsto h_s$ extends to an injection from the Artin braid group $B(W)$ to the Hecke algebra $H(W)$.  
\end{question}

I have no ideas for this. For a very particular kind of Coxeter groups it was solved recently by Paolo Sentinelli \cite{Paolito}, although it doesn't seem very generalizable. 

\vspace{0.3cm}

\hspace{-0.4cm}\textrm{\textbf{Sub-Question 3.}} {\color{blue}$(I13,WK25,D15)$}
Prove that the map $B(W) \rightarrow K^b(\mathcal{B(W)})$ sending $\sigma$ to $F_{\sigma}$ is injective. 
\vspace{0.3cm}

 Question \ref{inj} implies this question (oddly enough, the categorical conjecture is weaker than the non-categorical one). This Question appears in \cite[Conjecture 9.8]{Rouquier}. It is important because it would imply that Rouquier complexes give a categorification of the Artin braid group (and not just of a quotient of it). This conjecture was proved in \cite{KhovanovSeidel} when $W$ is the symmetric group and in \cite{thorge}  when $W$ is finite using the beautiful ``Garside theory''.  

\begin{question}\label{egag} {\color{blue}$(I15,WK1,D8)$}
Understand the relations between euclidean geometry and alcovic geometry. 
\end{question}
This is a vague question, but I hope that there is a concrete answer. 

\begin{question}\label{Burau} \danger $\ $ {\color{blue} $(I50,WK100,D30)$}
Is the Burau representation of the braid group $B(S_4)$  faithful?   
\end{question}
Go to the entry ``Burau representation'' in wikipedia for an explicit definition of this representation in terms of matrices. I am working on this problem with Geordie Williamson and David Plaza following \cite{KhovanovSeidel}. We have reduced the problem several times to something simpler, but it still seems elusive.  

\begin{question} $\ $ {\color{blue} $(I70,WK200,D45)$}
Is there a nontrivial knot with Jones polynomial equal to that of the unknot?
\end{question}
Historically, the most prominent approach in the search of a positive answer to this question is via Question \ref{Burau}. But I am pretty convinced that this is a bad approach, because I am almost sure that the answer to Question \ref{Burau} is yes (so this would not give the other knot with Jones polynomial equal to the unknot). So probably the best approach for this question is to try and study other representations of the braid group.

\begin{question}
 {\color{blue}$(I10,WK3,D10)$} Give a ``global'' description of Kazhdan-Lusztig double cells of the symmetric group. 
\end{question}
Kazhdan-Lusztig double cells for the symmetric group $S_n$ are parametrized by partitions $\lambda \vdash n$. 
For $w\in S_n$, there is a beautiful combinatorial way to see in what double cell does $w$ live, called ``Schensted's bumping algorithm'' (\cite[Section 22.2]{GBM}). Although a beautiful thing, for me it is too ``local'' in the sense that you have to do a long algorithm to know if some particular element lives in that particular double cell. Once you have calculated this algorithm, you can do left and right Knuth moves to move from $w$ to all other $w'$ living in the same double cell. Still local. 

I will quickly explain some other description that with Geordie Williamson we thought for a while  could work, but it doesn't. In any case, it can give an idea of what I mean by ``global''.

Let us call a partition $\{1, \dots, n\} = \bigsqcup P_m$
\emph{connected} if each $P_m$ is of the form $\{i, i+1, \dots, j \}$
for some $1 \le i \le j \le n$. We have an obvious bijection
\[
\{ \text{connected partitions of $\{ 1, \dots, n \}$} \} \longrightarrow \{
\text{standard parabolic subgroups of $S_n$}\}
\]
given by sending a partition to its stabiliser in $S_n$. Given a
partition $P = \bigsqcup P_m$ of $\{1, \dots, n \}$ we obtain the decreasing reordering Young diagram  $\gamma_P$. We define the \emph{type} of
$P$ to be $\lambda_P := \gamma_P^t$ (transpose diagram).
 Say  that $x \leq_w y$ (here the $w$ stands for ``weak'') if there exist a reduced expression $y = s_{i_1}s_{i_2} \cdots
s_{i_m}$ such that $x = s_{i_a} \cdots s_{i_{b}}$ for some $1 \le a
\le b \le m$. For $I \subset S$ let  $w_I$ be the longest element in the parabolic subgroup $W_I\subseteq W$ generated by $I$.
Given a partition $\lambda$ of $n$, let
\[
E_{\lambda} := \{ x \; | \; w_I \le_{w} x \text{ for some $I \subset S$
  of type $\lambda$} \},
\]
 Set
\[
D_{\lambda} := E_{\lambda} \setminus \bigcup_{\mu < \lambda } E_\mu.
\]
I think that $D_\lambda$ is the two-sided cell in $S_n$ corresponding to $\lambda$, for $n\leq 7$. I have not proved it, although this is certainly false for higher ranks. I am sorry for explaining a ``non-theorem'', but hope that it explains what kind of description could be expected. For more inspiration (on cells for infinite Coxeter groups) see \cite{Belo1}, \cite{Belo2}, \cite{Belo3} and \cite{Belo4}.

\begin{question}
{\color{blue}$(I80,WK70,D40)$} Prove that if $W$ and $W'$ are two Coxeter groups, $[x,y]$ and $[x',y']$ Bruhat intervals in $W$ and $W'$ respectively, then if  $[x,y]$ is isomorphic as a poset to $[x',y']$, then $p_{x,y}=p_{x',y'}$. 
\end{question}

This is called the \textit{combinatorial invariance conjecture} and it is due independently to George Lusztig and Matthew Dyer. I find it fascinating, because in other geometric situations (see the introduction of \cite{BLP}) this is not true, but in this case some strange miracle seems to happen. This has been a hot topic of research in the last years, mostly by Francesco Brenti's school (see the first page of \cite{BLP} to find lots of references). While I was writing this paper the preprint \cite{George} appeared where, with the help of artificial intelligence (Deep Mind by Google), Geordie Williamson and collaborators from Deep Mind elaborate a way to prove this conjecture for the symmetric group. A version of this paper appeared in the famous journal Nature \cite{Nature}, because it explains how Artificial Intelligence can help gain intuition and even converge to the ``correct'' intuition in mathematics.

\begin{question} {\color{blue}$(I13,WK5,D10)$}
Prove that $\widetilde{R}$
polynomials are log-concave. 
\end{question}
$R$ polynomials are strongly related to Kazhdan-Lusztig polynomials, as the latter were defined using the former in Kazhdan and Lusztig's original formulation. They are defined by the formula $$(h_{y^{-1}})^{-1}=\sum_{x\in W}R_{x,y}(v)h_x.$$
There is a unique polynomial $\widetilde{R}_{x,y}(v)\in \mathbb{N}[v]$ such that $\widetilde{R}_{x,y}(v-v^{-1})=R_{x,y}(v).$ The exponents in the polynomial  $\widetilde{R}_{x,y}(v)$ are all even or all odd.
Brenti conjectured   that  the sequence of non-zero coefficients in $\widetilde{R}_{x,y}(v)$ is log-concave (for the precise formulation see \cite[Conjecture 7.1]{Brentihistory}). Usually log-concavity is an indication of deep phenomena, usually a hidden Hodge theory \cite{Huh}. 

The idea I would propose to solve this is to use David Plaza's interpretation of these polynomials in terms of certain light leaves \cite{Plaza} (in Section 2 of that paper there is more information about the notation used here) and then mix it with Hodge theoretic methods \cite{EWhodge} for Soergel bimodules in the vein of June Huh's works \cite{Huh}.

\begin{question} {\color{blue}$(I90,WK100,D70)$}
Give a combinatorial description for the coefficients of Kazhdan-Lusztig polynomials. 
\end{question}

One of the fundamental open questions in algebraic combinatorics. There is an approach via ``canonical light leaves'' \cite{LWKL}, developed in an interesting case in \cite{LP}. One would essentially try to see what light leaves ``die'' when precomposed with the favorite projector. Hopefully there are some basic ``relations'' between light leaves that die.
It is a difficult approach it but seems to hide so many beautiful mysteries!

\begin{question} {\color{blue}$(I20,WK3,D10)$}
For $\widetilde{S}_n$, find explicitly the decomposition of the pre-canonical basis $N^{i+1}$ in terms of $N^i$.  
\end{question}

I explained a bit this set of bases in an example, in Section \ref{funcase}. For the  definition of the pre-canonical bases for any affine Weyl group and what is known about it see \cite{precan}.

\begin{question} \danger $\  $ {\color{blue}$(I30,WK10,D15)$}
Find formulas for all Kazhdan-Lusztig polynomials in the affine Weyl groups, assuming the knowledge of the formulas for the corresponding Weyl group.  
\end{question}
We are thinking this problem with Leonardo Patimo and David Plaza, although we are far from solving it. I think that we are going to be able to solve it for the ``lowest'' Kazhdan-Lusztig double cell $C$ (which is ``almost all the group'' in the sense that for any element $x\in W$ there is $s\in S$ such that $xs\in C$) using geometric Satake. But it seems that the higher the codimension of the cell, the more difficult it is to calculate the Kazhdan-Lusztig polynomials.

\begin{question} {\color{blue}$(I15,WK5,D9)$}
Prove the Forking path conjecture. 
\end{question}

For a precise statement see \cite[Question 6.2]{Gentle}, and for a proof for other type of groups see \cite{fp}.
This conjecture   roughly says the following. Let $x\in S_n$ and $s_1s_2\cdots s_n$ be a reduced expression of $x$. Apply maps of the form $\mathrm{id}\otimes f_{s,r} \otimes \mathrm{id}$ to $BS(s_1,s_2,\cdots,s_n)$ ``all the times you can'', passing through every possible reduced expression of $x$ and then come back to $BS(s_1,s_2,\cdots,s_n)$. The conjecture says that the map thus defined is independent of the path you chose. If so, this path morphism defines an idempotent in $BS(s_1,s_2,\cdots,s_n)$ whose image, at least in ranks $\leq 6$ is exactly the indecomposable $B_x$ but by general theory this can not be  so for every rank. So, if the conjecture was proved, the image in general would be something really mysterious. I hope that it is related with some kind of ``universal'' $p$-canonical basis, but this is highly speculative. 

\begin{question}. {\color{blue}$(I80,WK7,D100)$}
Find a formula for the favorite projector for all elements of all Coxeter groups. 
\end{question}
This question is really hard. Recall that in Section \ref{sbim} we mentioned the favorite projector, while the full definition is in \cite[Section 4.1]{LLL} The only cases known are (in characteristic zero): dihedral groups \cite{Bendihedral}, universal Coxeter groups (i.e. all $m(s,r)=\infty$) \cite{BenNico}, the longest element in $S_n$ \cite{Ben2} and  $\widetilde{S}_3$ \cite{LP}. The only case where this is done in positive characteristic is for $\widetilde{S}_2$ in \cite{pjones}.

\begin{question} {\color{blue}$(I200,WK100,D300)$}
Compute the $p$-canonical basis for all affine Weyl groups. 
\end{question}
This is really close to the last question, to find the favorite projector in positive characteristic. It seems to me that this question will be a guiding problem for the next decade at least, maybe even for the next 30 years. The closest we are now is with the billiards conjecture \cite{LusWil} (and its correction \cite{Jensen}) for a small part of $\widetilde{S}_3$. 

\vspace{0.3cm}

\hspace{-0.4cm}\textrm{\textbf{Sub-Question 15.a.}} {\color{blue}$(I20,WK50,D20)$}
\textit{Prove the billiards conjecture.}  

\vspace{0.3cm}
Idea: Find the favorite projectors for the relevant elements. This might be done with a mixture of the $p$-Jones-Wenzl projectors \cite{pjones} (as the $SL_2$ ingredient that might appear in every projector), Ben Elias' clasp conjectures \cite{clasp} and the Quantum Geometric Satake equivalence \cite{QSat}. For more details on this approach see Section $1.3$ of \cite{pjones}.

\vspace{0.3cm}

\hspace{-0.4cm}\textrm{\textbf{Sub-sub-Question 15.a.i}} {\color{blue}$(I20,WK50,D20)$}
\textit{Prove Ben Elias' clasp conjectures.}

\vspace{0.3cm}

\hspace{-0.4cm}\textrm{\textbf{Sub-Question 15.b.}} {\color{blue}$(I35,WK30,D40)$}
\textit{Determine the exact set of primes for which the $p$-canonical basis is different from the canonical basis for each affine Weyl group. }

\vspace{0.3cm}
Idea: it seems to me that one needs more data to predict this curve (I mean, for $SL_n(\mathbb{F}_p)$, the curve $n$ against $p$). For this, it might help to categorify the spherical module as we did with Geordie Williamson for the anti-spherical module \cite{AS}.

\vspace{0.3cm}

\begin{question} {\color{blue}$(I50,WK15,D40)$}
Generalize the definition of the category of Soergel bimodules to any complex reflection group.  
\end{question}

The only  result I know in this direction is the rank $1$ construction in \cite{Gobet}, although this question is as old as I can remember (I start remembering from 2008).

\begin{question} \danger $\ $ {\color{blue}$(I70,WK100,D20)$}
Produce a singular Soergel calculus. 
\end{question}
By this I mean an analogue of Soergel calculus for singular Soergel bimodules. 
I have been working with Ben Elias and Hankyung Ko in producing singular light leaves. We have a construction but we have not yet proved that  it is a basis for singular Soergel bimodules. The next step would be to find all the relations and thus produce a singular Soergel calculus.

\begin{question}\danger $\ $ {\color{blue}$(I10,WK3,D5)$}
Prove a categorical Schur-Weyl duality.
\end{question}
I am working on this with Juan Camilo Arias. The ``action'' of a monoidal category or a 2-category on an additive category is hard to check because you have to check it on objects, on morphisms and also several other compatibilities. Here a major issue is to understand what is the precise definition of a category of endofunctors ``fixed'' by the action of some category.

\begin{question} {\color{blue}$(I40,WK5,D50)$}
Prove Relative hard Lefschetz for iterated tensor products of Soergel bimodules.
\end{question}
I mean  Conjecture 3.4 in \cite{Relativehl}. I will not define all the symbols and concepts involved, but just for the reader to get an idea, it says the following. If $x_1,x_2\ldots,x_n$ are elements of $W$ (any Coxeter group) and $\rho_1,\rho_2,\ldots,\rho_{n-1}\in \mathfrak{h}^*$ are such that $\langle \rho_i,\alpha_s^{\lor} \rangle\geq 0$ for all $i$ and $s\in S$, then the map 
$$B_{x_1}\otimes_RB_{x_2}\otimes_R\cdots \otimes_R B_{x_n}\longrightarrow B_{x_1}\otimes_RB_{x_2}\otimes_R\cdots \otimes_R B_{x_n}$$ 
$$\hspace{1.3cm}
b_1\otimes b_2\otimes \cdots b_{n-1}\otimes b_n \hspace{0.2cm} \mapsto \hspace{0.2cm} b_1 \rho_1\otimes b_2\rho_2\otimes \cdots b_{n-1}\rho_{n-1}  \otimes b_n
$$
induces an isomorphism for all $i$
$$H^{-i}(B_{x_1}\otimes_R\cdots \otimes_R B_{x_n})\longrightarrow H^{i}(B_{x_1}\otimes_R\cdots \otimes_R B_{x_n}).$$

\begin{question}\label{bri} \danger $\ $ {\color{blue}$(I55,WK4,D25)$}
Compute the space of Bridgeland stability conditions for  $K^b(\mathcal{B}(W))$.
\end{question}
It is quite rare that the paper in which a concept is introduced is still the best source to study the concept. In my opinion, this is the case with Bridgeland stability conditions, a manifold associated to any triangulated category \cite{Bridgeland}.

Question \ref{bri} seems absolutely fascinating to me. We have been thinking about it with Anthony Licata. It is a well established fact that the space of Bridgeland stability conditions  is a fundamental invariant of a triangulated category. 
I think that the answer to this question for $\mathcal{B}(W)$ in characteristic zero can give interesting new information about Kazhdan-Lusztig theory, and in positive characteristic it can give critical information about the $p$-canonical basis. That said, for me this is a fundamental question on its own. It seems that the whole space of stability conditions might be too massive, so maybe a variation of the construction of Bridgeland stability conditions might be in place, maybe related to the highest weight category structure of linear complexes in $K^b(\mathcal{B}(W))$ \cite{LibWil}, \cite[Theorem 26.17]{GBM} or maybe in the vein of the ``real variations of stability'' in \cite{Anno}. The braid group action (looking at Rouquier complexes as self-equivalences) is important.

\begin{song*}\label{desp} 
{\color{blue}$(I35,WK1,D50)$}  What is $\mathbb{F}_q$ when $q$ is not the power of a prime number, but any element of some field?
\end{song*}

 I would suggest to read this paper \cite{Deligne} by  Pierre Deligne. Instead of generalizing $S_t$ for $t$ not a natural number, he focuses on generalizing  the category $\mathrm{Rep}(S_n)$ of representations of $S_n$. In that paper, Deligne  defined the categories  $\mathrm{Rep}(S_t)$ and $\mathrm{Rep}(GL_t)$ for $t$  any element that is not a positive integer in any field $k$ of characteristic zero.  I would also encourage the reader to read for inspiration two very fun papers, one by Stephen Schanuel \cite{Schanuel} (this paper by James Propp explains it in easy terms \cite{Propp}) and one  by Daniel Loeb \cite{Loeb}, where the idea of sets with a negative number of elements is explored\footnote{If one uses Schanuel's improved version of the Euler characteristic, one could be tempted to admit that $\mathbb{F}_{-1}=\mathbb{R}.$}.

\bibliography{MAp2.bib} 
\bibliographystyle{alpha}

\end{document}